\documentclass[final,3p,times,authoryear,fleqn]{elsarticle}

\usepackage{amsmath,amssymb,amsthm}
\usepackage{geometry}%\geometry{left=2cm,right=2cm,top=2cm,bottom=2.4cm}
\usepackage[colorlinks,linkcolor=blue,urlcolor=blue,anchorcolor=red,citecolor=blue]{hyperref}
\usepackage{comment}
\usepackage{algorithm}
\usepackage{algorithmic}

\usepackage{ulem}

\usepackage{subcaption}
\usepackage{pxfonts,times}

\usepackage{amsmath,amssymb,bm}
\newcommand{\Vt}[1]{\bm{#1}}
\newcommand{\VtI}{\bm{I}}

\newcommand{\VtN}{\bm{N}}

\newcommand{\VtQ}{\bm{Q}}

\newcommand{\Rmd}{\mathrm{d}}

\newcommand{\ClN}{\mathcal{N}}

\newcommand{\ClT}{\mathcal{T}}

\providecommand{\inN}{\in\ClN}

\providecommand{\inT}{\in\ClT}

\usepackage{prettyref}
\newrefformat{eq}{Eq.~(\ref{#1})}
\newrefformat{fig}{Fig.~\ref{#1}}
\newrefformat{tab}{Table~\ref{#1}}
\newrefformat{assum}{{Assumption~\ref{#1}}}
\newrefformat{defn}{{Definition~\ref{#1}}}
\newrefformat{lemma}{{Lemma~\ref{#1}}}
\newrefformat{theo}{{Theorem~\ref{#1}}}
\newrefformat{prop}{{Proposition~\ref{#1}}}
\newrefformat{coro}{{Corollary~\ref{#1}}}
\newrefformat{app}{{\ref{#1}}}
\newrefformat{sec}{Section~\ref{#1}}
\newrefformat{conj}{{Conjecture~\ref{#1}}}

\newtheorem{definition}{Definition}[section]
\newtheorem{theorem}{Theorem}[section]
\newtheorem{lemma}{Lemma}[section]
\newtheorem{corollary}{Corollary}[section]
\newtheorem{example}{Example}[section]

\newtheorem{proposition}{Proposition}[section]
\numberwithin{equation}{section}
\newtheorem{conjecture}{Conjecture}[section]

\makeatletter
\def\th@plain{\upshape}
\makeatother

\makeatletter
\def\inN{\in\ClN}

\def\inT{\in\ClT}

\providecommand\Red[1]{#1}

\makeatother

\allowdisplaybreaks

\journal{Elsevier}

\begin{document}

\begin{frontmatter}



\title{Dynamic traffic assignment in a corridor network: Optimum versus Equilibrium}


\author[0]{Haoran Fu}
%[orcid=0000-0002-7204-0672]
\ead{haoranfu@xjtu.edu.cn}

%\credit{Conceptualization, Methodology, Software, Writing - original draft, Writing - review \& editing}
\address[0]{School of Management, Xi'an Jiaotong University, Xi'an, 710049, China}

\cortext[cor1]{Corresponding author.}
\pdfstringdefDisableCommands{%
  \def\corref#1{}%
}

\author[1]{Takashi Akamatsu\corref{cor1}}

\ead{akamatsu@plan.civil.tohoku.ac.jp}

%\credit{Conceptualization, Formal analysis, Supervision, Writing - review \& editing}

\address[1]{Graduate School of Information Sciences, Tohoku University, 6-6 Aramaki Aoba, Aoba-ku, Sendai, Miyagi, 980-8579, Japan}

\author[1]{Koki Satsukawa\corref{cor1}}

%\address[2]{New Industry Creation Hatchery Center, Tohoku University, 6-6 Aramaki Aoba, Aoba-ku, Sendai, Miyagi, 980-8579, Japan}

\ead{satsukawa@tohoku.ac.jp}

%\credit{Formal analysis, Validation, Writing - review \& editing}

\author[2]{Kentaro Wada}

\address[2]{Faculty of Engineering, Information and Systems, University of Tsukuba, 1-1-1 Tennodai, Tsukuba, Ibaraki, 305-8573, Japan}

\ead{wadaken@sk.tsukuba.ac.jp}

%\credit{Software, Writing - review \& editing}

\begin{abstract}
This study investigates dynamic system-optimal (DSO) and dynamic user equilibrium (DUE) traffic assignment of departure/arrival-time choices in a corridor network. 
The morning commute problems with a many-to-one pattern of origin-destination demand and the evening commute problems with a one-to-many pattern are considered. 
Specifically, a novel approach to derive closed-form solutions for both DSO and DUE problems is developed.
We first derive a closed-form solution to the DSO problem based on the regularities of the cost and flow variables at an optimal state.  
By utilizing this solution, we prove that the queuing delay at a bottleneck in a DUE solution is equal to an optimal toll that eliminates the queue in a DSO solution under certain conditions of a schedule delay function.
This enables us to derive a closed-form DUE solution by using the DSO solution.
We also show the theoretical relationship between the DSO and DUE assignment. 
Numerical examples are provided to illustrate and verify the analytical results.
\end{abstract}

%Research highlights
%\item An analytical approach is proposed for dynamic system-optimum (DSO) and dynamic user equilibrium (DUE) traffic assignment of departure/arrival-time choices.
%
%\item Closed-form solutions are obtained for DSO and DUE problems in corridor networks.
%
%
%\item The existence of DSO and DUE solutions with the same patterns of optimal tolls and queuing delays is revealed.
%
%
%\item Morning and evening commute problems are studied to clarify the differences in the relationship between DSO and DUE.


\end{frontmatter}


\setcounter{tocdepth}{2}
%\tableofcontents
% no subsubsection


%%%%%%%%%%%%%%%%%%%%%%%%%%%%%%%%%%%%%%%%%%%%%%%%%%%%%%%%%%%%%%%%%%%%
% Section 1
\section{Introduction}
\subsection{Background and purpose}
Traffic congestion during rush hours is a common issue in metropolitan areas.
It forces commuters to select specific departure times and compete with other commuters for limited road capacity.
To understand and characterize this problem, prior research has focused on dynamic user equilibrium (DUE) assignment and dynamic system optimal (DSO) assignment of commuters' departure-time choices \citep[see][for a recent comprehensive review]{Li2020l}.
%The DUE assignment problem considers a pure equilibrium of departure-time choice, whereas the DSO assignment problem considers the optimal assignment of departure-time choices that minimizes the total transportation cost for the commuters. 

Single-bottleneck models for DUE assignment (with a point queue model) are concise and effective in characterizing the departure-time choice behavior of commuters and has been successfully analyzed theoretically \citep{vickrey1969congestion,hendrickson1981schedule,smith1984existence,daganzo1985uniqueness,newell1987morning,arnott1993departure,lindsey2004existence,akamatsu2021new}.
The DSO assignment and its relation to the DUE assignment have also been extensively studied, 
%which lays the 
and the findings have established a theoretical foundation of a (first-best) travel demand management (TDM) policy to eliminate congestion, such as congestion pricing, tradable permits, and tradable credit schemes \citep[e.g.,][]{arnott1993departure,akamatsu2006tradable,Nie2013b}.
%For a recent comprehensive review, refer to \cite{Li2020l}. 
The models, however, have intrinsic limitations in terms of representing the spatial dynamics of congestion. 
Congestion may exist at multiple bottlenecks in a network, and the effects of mutually interacting queues must be accounted for.
%\textcolor{red}{and the effects of interactions between the congestion and traffic flow are required to be accounted for.} 

%There are two streams of research for analyzing the spatial dynamics of congestion in the departure-time choice framework.
% The first stream of research extends the model (or network) to ones with multiple bottlenecks for the DUE assignment  \citep{kuwahara1990equilibrium,arnott1993properties,Daniel2009,Lago2007,akamatsu2015corridor} and for the DSO assignment \citep{yodoshi2008pareto,shen2009morning,osawa2018first}. 
Many studies have focused on extensions to corridor networks with multiple tandem bottlenecks for DUE assignment  \citep{kuwahara1990equilibrium,arnott1993properties,Daniel2009,Lago2007,akamatsu2015corridor} and DSO assignment \citep{yodoshi2008pareto,shen2009morning,osawa2018first}. 
These studies have provided useful insights into the distribution and interaction of congestion along a corridor.
However, the {\it proof by cases} approach employed in most studies faces the challenge in which the number of cases rapidly increasing with an increase in the number of bottlenecks, as noted in \cite{arnott2001} and \cite{arnott2011corridor}, that is, it is difficult to obtain closed-form solutions for networks with an arbitrary number of bottlenecks using this approach.  
Furthermore, none of these studies have clarified the relationship between DUE and DSO assignment in a general manner. 
Note that some recent studies employed the Lighthill-Whitham-Richards (LWR) model rather than the classical point queue model in a corridor network primarily to analyze the spatial dynamics of the ``flow congestion" \citep{arnott2011corridor,DePalma2012,Wang2016g,Li2017a}.
The models in these studies, however, cannot handle queuing/bottleneck congestion, because of a uniform road width assumption (i.e., there is no bottleneck), and are analytically intractable.  

%The second stream of research employs the LWR traffic flow model primarily for modeling the spatial dynamics of the ``flow congestion" within the DUE and DSO assignments in a corridor network \citep{arnott2011corridor,DePalma2012,Wang2016g,Li2017a}. 
% The LWR model could describe the spatial traffic dynamics more realistically than the classical point queue model. 
% Unfortunately, however, queuing congestion, which is ubiquitous, never occurs because the models assume a corridor network with uniform width. 
%these models cannot capture the effects of mutually interacting queues at multiple bottlenecks because they assume a corridor network with uniform width and queues never occur. 
% Furthermore, the latter stream is more analytical intractable than the former. 

The purpose of this research is to investigate the theoretical properties of the DSO and DUE assignment of departure-time choices in a corridor network with multiple bottlenecks, using the point queue model. 
% Specifically, we aim to derive closed-form solutions and clarify the relationships between these solutions for the two problems.
To this end, we develop an analytical approach to solve the DSO and DUE problems sequentially, through the following steps. 
First, we formulate the DSO assignment without queues for a morning commute with a many-to-one origin-destination (OD) demand pattern as an infinite-dimensional linear programming (LP) problem. 
We show (i) an inclusion relationship between the destination arrival-time windows\footnote{Similar regularities were observed in the departure time choice equilibrium problems in corridor networks with many-to-one OD demands, but with the different technologies of congestion \citep[e.g.,][]{Tian2007,arnott2011corridor}.} and (ii) the regularities of cost and flow variables at the optimal state.
%an equivalent equilibrium state under a first-best TDM policy.
%First, we model the DSO problem in the morning commute with a many-to-one origin-destination demand pattern as an equilibrium state under a tradable network permit (TNP) scheme.
%This equilibrium under the TNP scheme (TNP equilibrium) is equivalent to a DSO solution.
%At the TNP equilibrium, we reveal an inclusion relationship of departure/arrival-time windows, and regularities of cost and flow variables.
These properties enable us to construct a reduced network and derive a closed-form solution for the DSO problem. 
Second, we prove that the optimal prices for bottlenecks in the DSO solution are equal to the queuing delays at bottlenecks in a DUE solution under certain conditions of a schedule delay function. 
By exploiting this fact and the DSO solution obtained in the first step, we obtain a closed-form solution to the DUE problem formulated as an infinite-dimensional linear complementarity problem (LCP). 
Then, we establish the relationships between flow/cost patterns under the DSO and DUE assignment and perform welfare analyses of first-best/second-best TDM policies. 
Furthermore, we show that a similar approach is applicable to the DSO and DUE assignment in the evening commute with a one-to-many OD demand pattern. 
The analysis shows that the DSO assignment in the morning and evening commute are symmetric, which means that they are mathematically equivalent; contrarily, the formulations and solutions of the DUE assignment differ for two problems.
% The analysis shows that the DSO problems in the morning and evening commute are symmetric, which means they are mathematically equivalent.
% On the contrary, the formulations and solutions of the DUE problem differ for the morning and evening commute although the analytical approach for the morning commute is applicable on the evening commute problem.

%Given the above existing studies, 
The contributions of this study are as follows:  
\begin{itemize}
    \item[(1)] A closed-form solution for the DSO problem is derived. 
    \item[(2)] The closed-form solution to the DUE problem is derived under certain conditions of a schedule delay function.
    \item[(3)] It is proved that the queuing delay at each bottleneck in the DUE assignment is exactly the same as the optimal (toll/permit) price that eliminates the queue in the DSO assignment (under the same conditions as mentioned in (2)). 
    \item[(4)] The (a)symmetric assignment properties of morning and evening commute problems are identified.      
\end{itemize}
Several remarks regarding each contribution are provided herein. 
(1) To circumvent the difficulty in the conventional approach based on the proof by cases, we introduce the notion of a reduced network, which enables us to derive the DSO closed-form solution.
%We successfully obtain the DSO closed-form solution by introducing a notion of a reduced network that circumvents . 
%reducing many cases in the original corridor network to only one case in a reduced network. 
The DSO closed-form solution is linked with the other contributions and also leads to a useful extension of the model, including more choice dimensions \citep{osawa2018first}.
(2) This is a theoretically and mathematically remarkable result. 
%Obtaining the closed-form solution to the DUE problem that is formulated as (infinite-dimensional) LCP is a theoretically/mathematically remarkable result.
Indeed, the DUE assignment was considered to be difficult to solve analytically \citep{arnott2001,arnott2011corridor,akamatsu2015corridor}; it is not trivial at all that the solution of an LCP (for the DUE assignment) is obtained from the solution of an LP problem (for the DSO assignment). 
(3) The relationship between the price variables in the DSO and DUE assignment not only enables us to derive the DUE solution from the DSO solution but also provides a powerful basis for the welfare analysis of several important (second-best as well as first-best) TDM policies. 
(4) While \cite{akamatsu2015corridor} suggested with a numerical example that the properties of the DUE assignment for morning and evening commute problems are asymmetric, the present study is analytic and provides more concrete insight into this point.

The remainder of this paper is organized as follows. 
After reviewing the related literature in the remainder of this section, Section \ref{Sec:MorningCommute} introduces the DSO problem for the morning commute in a corridor network and presents the derivation of a closed-form solution. 
Section \ref{Sec:MC_DUE} provides a closed-form DUE solution based on the results presented in the previous section. 
We here also show the relationship between the DSO and DUE assignment. 
Section \ref{Sec:EveningCommute} further studies the DSO and DUE assignment for the evening commute.
\prettyref{sec:Num} provides numerical examples. 
Finally, \prettyref{sec:Conclusion} concludes this paper.

\subsection{Related literature}
\subsubsection{Corridor networks with multiple bottlenecks}
\cite{kuwahara1990equilibrium} firstly examined a DUE assignment in more than one bottleneck (two-tandem bottleneck) networks with a many-to-one OD demand pattern, and showed that the bottleneck service priority differs for different origins. 
In a similar two-tandem bottleneck network, \cite{arnott1993properties} identified a capacity-increasing paradox; \cite{Daniel2009} later supported this theoretical hypothesis through a large-group laboratory experiment.
\cite{Lago2007} extended the model to incorporate queue spillovers and merging effects and provided policy insights accordingly.
However, the {\it proof by cases} approach adopted in these studies has a major limitation, as mentioned earlier.  
%came to face the challenge in which the number of cases increasing rapidly in the number of bottlenecks as noted in \cite{arnott2001}. 
Therefore, as an alternative, \cite{akamatsu2015corridor} developed a mathematical programming approach (an infinite-dimensional LCP formulation) for models in a corridor network with an arbitrary number of bottlenecks and obtained some fundamental theoretical results (i.e., existence and uniqueness of equilibrium). 
Although the study did not obtain closed-form solutions, the LCP formulation is an essential building block for characterizing the equilibrium flow and cost patterns in the present study.
%Although this study cannot obtain closed-form solutions \textcolor{red}{since the problem is formulated as an infinite-dimensional complementarity problem}, it is an essential building block for the current study. 

%In contrast to the DUE problem, 
Several studies have been conducted on the departure time choice DSO assignment without queues or the equivalent equilibrium assignment under a first-best TDM policy in a tandem bottleneck network (with a many-to-one OD demand pattern). %but they have treated an arbitrary number of bottlenecks.
\cite{yodoshi2008pareto} examined whether the Pareto improvement property is achieved at an equilibrium assignment under the {\it tradable network permits} (TNP) system \citep{akamatsu2006tradable,akamatsu2007system,akamatsu2017tradable} in a two-tandem bottleneck network.  
\cite{shen2009morning} considered a corridor network in which a freeway with multiple bottlenecks connects to a surface network with large capacity through on-/off-ramps and showed some of the features of the optimal traffic flow and toll patterns as well as a graphical solution procedure. 
\cite{osawa2018first} studied a model of DSO assignment that integrates the short-term problem (departure time choice with tolling) and the long-term problem (job and residential location choice). 
For the short-term problem (the departure time choice DSO assignment), they derived a generalized result by utilizing the findings obtained in the earlier versions of the present study \citep{Fu2016,Fu2018} and \cite{akamatsu2021new} as theoretical building blocks. 
% It should be noted that the result was obtained based on that for homogeneous commuters in the earlier versions of the current paper \textcolor{red}{\citep{Fu2015,Fu2018}} and 
However, their primary objective is fundamentally different from that of ours (i.e., their focus was to investigate the properties and policy implications of {\it long-term} equilibria under the first-best short-term TDM policies)\footnote{Unlike the discussion on the short-term/departure-time choice DSO problem in \cite{osawa2018first}, the present study provides complete, detailed arguments and is also more accessible with several illustrations and numerical examples.}. 
Furthermore, none of the aforementioned studies clarified the relationship between the DSO and DUE assignment with an arbitrary number of bottlenecks. % that provides a powerful basis for the welfare analysis of the second-best as well as the first dynamic control policies. 

\subsubsection{Corridor networks with flow congestion}
%In recent years, several studies employ the LWR traffic flow model for the spatial dynamic modeling of the ``flow congestion" within the DUE and DSO assignments in a corridor network. 
% In recent years, several studies employ the LWR traffic flow model rather than the classical point queue model for the spatial dynamic modeling of the ``flow congestion" within the DUE and DSO assignments in a corridor network \citep{arnott2011corridor,DePalma2012,Wang2016g,Li2017a}.
\cite{arnott2001} and \cite{arnott2011corridor} introduced the so-called {\it corridor problem}, in which continuum-entry points and a single destination exist along a corridor with uniform width, and the the LWR traffic flow model is assumed. 
The authors proposed a solution that meets a departure time choice equilibrium condition under restricted assumptions (e.g., a restricted demand pattern and piece-wise linear schedule delay function with no late arrival), but they could not provide a complete solution. 
\cite{Wang2016g} formulated the corridor problem into a partial differential complementarity system, which was then numerically solved using the Godunov scheme.
However, as mentioned earlier, the primary motivation for employing the LWR model is to capture the effects of spatial dynamics of ``flow congestion" rather than the queuing congestion at bottlenecks.
%focus is more on capturing the effects of spatial dynamics of ``flow congestion" than the queuing congestion. 
%Although the LWR model could describe the spatial traffic dynamics more realistically than the classical point queue model, queuing congestion, which is ubiquitous, never occurs in the corridor problem due to the uniform width assumption. 
% Therefore, which model is preferable depends on context \citep{DePalma2012}. 
Furthermore, because of its analytical intractability, % even for a uniform width corridor, 
\cite{DePalma2012} \citep[and][]{Li2017a} has turned their focus back on a simpler problem termed the ``single-entry corridor problem," which was first examined by \cite{Newell1988}.   
It is considered as a generalization of the single-bottleneck (entry) model than one with multiple bottlenecks (entries), which is  analyzed in the present study. 
Note that \cite{Tian2007} investigated a similar system of a many-to-one public transit with the in-vehicle crowding effect; however, the congestion technology adopted was almost the same as the flow congestion in static traffic equilibrium assignment.

% \prettyref{Sec:MorningCommute} introduces the DSO problem for the morning commute in a corridor network with the formulation of a TNP equilibrium.
% The analytical solution to the DSO problem is derived in \prettyref{sec:DSOsolution}, and the DUE problem was analyzed using these results to provide a closed-form DUE solution in \prettyref{sec:DSOvsDUE}.
% In addition, the DSO and DUE problems for the evening commute are further studied in \prettyref{sec:Eve}.
% \prettyref{sec:Num} provides numerical examples. 
% \prettyref{sec:Conclusion} concludes this paper.

%%%%%%%%%%%%%%%%%%%%%%%%%%%%%%%%%%%%%%%%%%%%%%%%%%%%%%%%%%%%%%%%%%%%
% Section 2

\begin{figure}[t]
	\centering
	\hspace{0mm}
    \includegraphics[width=80mm,clip]{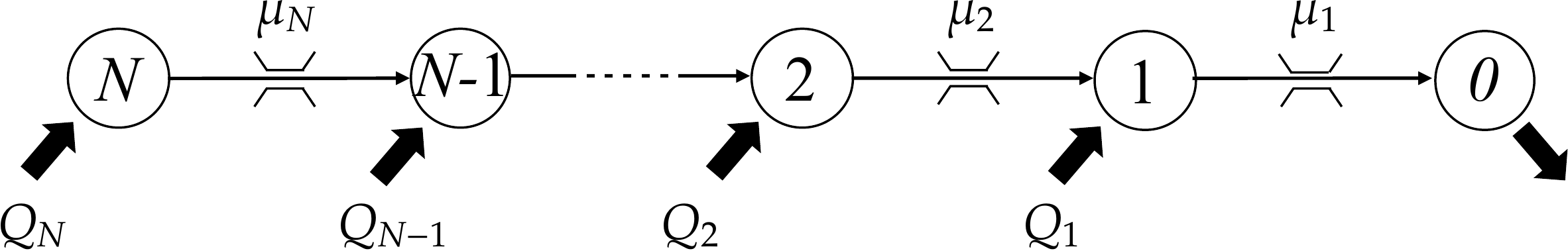}
    \vspace{-2mm}
	\caption{Corridor network with $N$ on-ramps and $N$ bottlenecks}
    \vspace{-0mm}
    \label{Fig:Corridor}
\end{figure}

\section{Dynamic system optimal assignment for the morning commute}\label{Sec:MorningCommute}

\subsection{Networks}\label{Sec:MC_Networks}

We consider a freeway corridor network consisting of $N$ on-ramps (origin nodes) and a single off-ramp (destination node), as shown in Figure~\ref{Fig:Corridor}.
These nodes are numbered sequentially from Destination node $0$ to the most distant origin node $N$.
A set of origin nodes is denoted by $\mathcal{N}\equiv\{1,\ldots,N\}$.
There is a single bottleneck with a finite capacity $\mu_{i}$ immediately downstream of the origin $i\in\mathcal{N}$, which is referred to as the bottleneck $i$.
At each bottleneck, a queue is formed when the arrival flow rate exceeds the capacity.
The queue evolution and associated queuing delay are modeled using the standard point queue model in accordance with the first-in-first-out (FIFO) principle.
The cumulative arrival and departure flows for the bottleneck $i$ by time $t$ are denoted by $A_{i}(t)$ and $D_{i}(t)$, respectively.
The arrival and departure flow rates at bottleneck $i$ at time $t$ are expressed as follows:
\begin{align}
\lambda_{i}(t) = \cfrac{\mathrm{d}A_{i}(t)}{\mathrm{d}t},\quad x_{i}(t) = \cfrac{\mathrm{d}D_{i}(t)}{\mathrm{d}t}.
\end{align}
The free-flow travel time from origin (and bottleneck) $i$ to the destination is denoted by $c_{i}$.

From each origin $i$, $Q_{i}$ commuters enter the network and reach the destination during the morning rush-hour $\mathcal{T}\equiv [0,T]$.
Commuters are considered as a continuum, and the total mass $Q_i$ at each origin $i$ is constant.
Each commuter incurs a schedule delay cost that is associated with the deviation from the desired arrival time (e.g., working start time) to the destination.
We assume that all commuters are homogeneous; thus, they have the same desired arrival time $t_{d}$, same value of time, and same penalty function $s(t)$ for commuters arriving at the destination at time $t$.
Moreover, we assume that the function $s(t)$ is strictly quasi-convex and piecewise differentiable for time periods $[0,t_{d}]$ and $[t_{d}, T]$.

\subsection{Formulation of dynamic system optimal assignment}\label{Sec:MC_DSOFormulation}
Under the abovementioned conditions, we first consider the dynamic system optimal (DSO) assignment.
Here, we define a DSO state as a state in which the total transport cost in the network is minimized without queues, i.e., congestion externalities are completely eliminated.
Mathematically, a DSO state is defined as a solution to the following infinite-dimensional linear programming problem: 
\begin{align}
\text{\textbf{[DSO-LP]}}\quad
\min_{\mathbf{q}\geq \mathbf{0}}
&\sum_{i\in \mathcal{N}}\int_{t\inT}\left(s(t) + c_{i}\right) q_i(t)\Rmd t  \label{Eq:DSOLP_Obj} \\
\text{s.t.}\quad
&\sum_{j = i}^{N}q_j(t)\leq\mu_i &\forall i\in \mathcal{N},\forall t\inT,\label{Eq:DSOLP_Const1}\\
&\int_{t\inT}q_i(t)\Rmd t=Q_i  &\forall i\in \mathcal{N},\label{Eq:DSOLP_Const2}
\end{align}
where $q_{i}(t)$ is the destination arrival flow rate of commuters at time $t$ departing from origin $i\in \mathcal{N}$ (hereinafter, the commuters are referred to as $i$\textit{-commuters}).
The objective function~\eqref{Eq:DSOLP_Obj} is the total schedule and travel costs when there is no queuing delay.
The first constraint~\eqref{Eq:DSOLP_Const1} is the capacity constraint on each link that ensures no queuing delay, and the second constraint~\eqref{Eq:DSOLP_Const2} is the flow conservation condition for each origin-destination (OD) pair.
It should be noted that the arrival flow rate $q_{i}(t)$ has the following relationship with departure flows from bottlenecks:
\begin{align}
&q_{i}(t) = 
\begin{cases}
x_{i}(\sigma_{i}(t)) - x_{i+1}(\sigma_{i+1}(t)),	&\quad \forall i\in\mathcal{N}\setminus \{ N\},\\
x_{N}(\sigma_{N}(t)), 									&\quad i = N
\end{cases}
\quad \forall t\in\mathcal{T},\label{Eq:DSO_FlowConservation}\\
&\text{where}\quad \sigma_{i}(t) = t - c_{i},\label{Eq:DSO_Sigma}
\end{align}
where $\sigma_{i}(t)$ represents the departure time from bottleneck $i$ of commuters, with a destination arrival time $t$.
The optimality conditions of the problem are given by the following conditions~\citep[for example,][]{Luenberger1997,akamatsu2021new}: 
\begin{align}
&\begin{cases}
s(t) + c_{i} + \sum_{j = 1}^{i}p_{j}(t) = \rho_{i}	\quad	&\text{if}\quad q_{i}(t) > 0,\\
s(t) + c_{i} + \sum_{j = 1}^{i}p_{j}(t) \geq \rho_{i}	\quad	&\text{if}\quad q_{i}(t) = 0,
\end{cases}
&\quad \forall i\in\mathcal{N},\forall t\in\mathcal{T},\label{Eq:DSO_KKT1}\\
&\begin{cases}
\sum_{j = i}^{N}q_{j}(t) = 		\mu_{i}	\quad	&\text{if}\quad p_{i}(t) > 0,\\
\sum_{j = i}^{N}q_{j}(t) \leq \mu_{i}	\quad	&\text{if}\quad p_{i}(t) = 0,
\end{cases}
&\quad \forall i\in\mathcal{N},\forall t\in \mathcal{T},\label{Eq:DSO_KKT2}\\
&\int_{t\inT}q_i(t)\Rmd t=Q_i  &\forall i\in \mathcal{N},\label{Eq:DSO_KKT3}
\end{align}
where $p_{i}(t)$ and $\rho_{i}$ are the Lagrange multipliers for Constraints~\eqref{Eq:DSOLP_Const1} and \eqref{Eq:DSOLP_Const2}, respectively.

%It is worth noting that these optimality conditions can be \textcolor{red}{economically interpreted} as \textit{equilibrium conditions under a first-best TDM policy} (we refer to such equilibrium as \textit{pricing equilibrium}) that completely eliminates bottleneck congestion.
It is worth noting that these optimality conditions have several economic interpretations, such as \textit{equilibrium conditions under a first-best TDM policy} (such equilibrium is referred to as \textit{pricing equilibrium}) that completely eliminates bottleneck congestion.
For example, we can interpret the conditions as equilibrium conditions under an optimal dynamic congestion pricing scheme.
This indicates that the Lagrange multipliers $p_{i}(t)$ can be regarded as the congestion price charged to the commuters with destination-arrival time $t$, and $\rho_{i}$ can be regarded as the equilibrium commuting cost of $i$-commuters.
Eq.~\eqref{Eq:DSO_KKT1} can then be interpreted as the departure time equilibrium condition for commuters, where the commuting cost of $i$-commuters with destination-arrival time $t$ is expressed as the sum of the schedule delay, free-flow travel time, and prices charged to the commuters, as follows: 
\begin{align}
v_{i}(t) = s(t)  +  c_{i} + \sum_{j = 1}^{i}p_{j}(t).
\end{align}
Eqs.~\eqref{Eq:DSO_KKT2} and \eqref{Eq:DSO_KKT3} can be interpreted as the conditions for optimal pricing required to eliminate congestion externalities and flow conservation conditions, respectively.

An alternative interpretation is the equilibrium conditions under a \textit{tradable network permit} (TNP) scheme~\citep[][]{akamatsu2006tradable,akamatsu2007system,akamatsu2017tradable}.
Under this scheme, the road manager issues a permit that allows for a permit holder to pass through a bottleneck at a pre-specified time period, such that the number of permits does not exceed the capacity to eliminate the congestion.
Each commuter is required to purchase a set of permits corresponding to a set of links passed by the commuter from trading markets where the road manager sells the permits.
The Lagrange multiplier $p_{i}(t)$ is then interpreted as the market-clearing price under the TNP scheme, and Condition~\eqref{Eq:DSO_KKT2} is interpreted as the demand-supply equilibrium condition for each permit trading market under the TNP scheme, where the demand of time period $t$ is equal to the flow $\sum_{j = 1}^{N}q_{j}(t)$, and the maximum supply of the permit is given by the bottleneck capacity $\mu_{i}$.

It should be noted that such an interpretation from the perspective of equilibrium conditions allows us to obtain the closed-form solution of the DSO assignment problem.
In the following section, we present an analytical approach for solving the DSO problem from the equilibrium conditions and its closed-form solutions.

\subsection{Closed-form solution of the DSO assignment problem}
In this section, we first introduce definitions to simplify the presentation of the pricing equilibrium.
Thereafter, we establish the inclusion relationship between the arrival time windows of commuters, which allows us to derive the useful properties for analyzing the flow and cost patterns at the pricing equilibrium.
We then obtain the closed-form solution of the pricing equilibrium based on these properties. 

\subsubsection{Definitions}
We first define an \textit{arrival time window} of commuters from each origin, which is the range of time within which the commuters arrive at the destination.
The earliest and latest arrival times of the $i$-commuters at the destination are denoted by $t_{i}^{-}$ and $t_{i}^{+}$, respectively.
The arrival time window of $i$-commuters is denoted by $\mathcal{T}_{i}\equiv [t_{i}^{-}, t_{i}^{+}]$, and its length is denoted by $T_{i}\equiv |\mathcal{T}_{i}|$.
We denote by $\bar{s}(T_{i})$ the schedule delay as a function of $T_{i}$: 
\begin{align}
\bar{s}(T_i) = s(t) = s(t+T_i),\quad \forall i\in\mathcal{N}.\label{Eq:Defi_ScheduleDelay}
\end{align}
Given that the schedule delay function is strictly quasi-convex, the schedule delay $\bar{s}(T_i)$ is uniquely determined for a given length of the arrival time window $\mathcal{T}_{i}$.

We then introduce the following concept of the \textit{false bottleneck}: 
\begin{definition}(\textbf{False bottleneck})
A bottleneck $i\in\mathcal{N}$ is defined as a \textit{false bottleneck} when optimal prices on the bottleneck are always zero in the pricing equilibrium, i.e., $p_{i}(t) = 0$, $\forall t\in\mathcal{T}$.
\end{definition}
\noindent From the definition, if bottleneck $i$ is a false bottleneck, the following relationship between the equilibrium costs holds: 
\begin{align}
\rho_i - c_i = \rho_{i-1} -c_{i-1}.
\end{align}
Thus, the false bottleneck is always under free-flow conditions, i.e., passing the bottleneck at any time does not influence the commuting cost for any commuter.
This indicates that a false bottleneck does not influence the essential properties of pricing equilibrium.
However, it also causes arbitrariness of the equilibrium flow patterns of the closed-form solution, and brings some unnecessary complications for the analysis.

\begin{figure}[t]
	\centering
	\hspace{0mm}
    \includegraphics[width=80mm,clip]{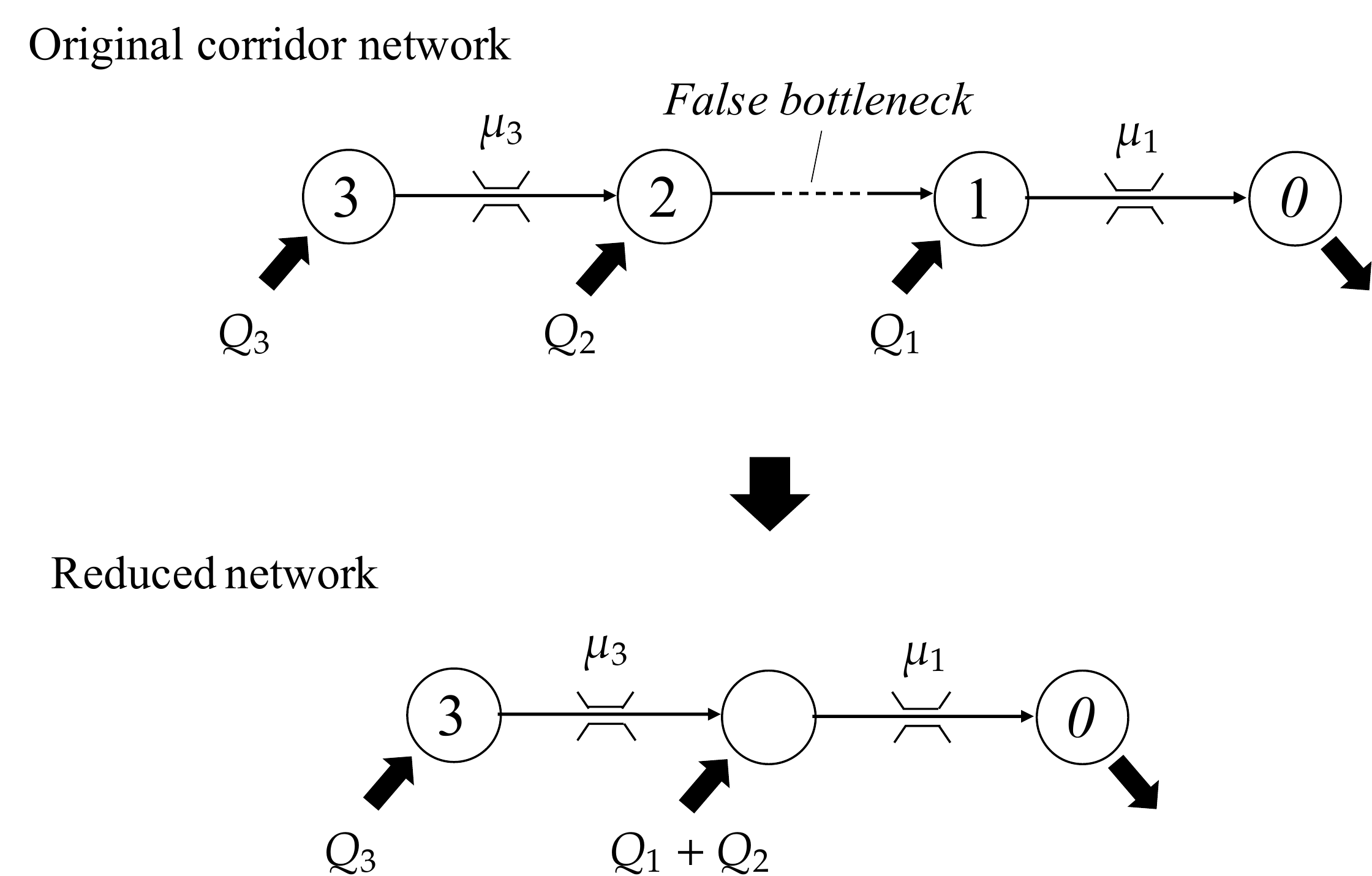}
    \vspace{-2mm}
	\caption{Example of constructing a reduced network in the case where $N = 3$ and bottleneck $2$ is a false bottleneck}
    \vspace{-0mm}
    \label{Fig:ReducedNetwork}
\end{figure}

To avoid the complications, we introduce the concept of a \textit{reduced network}, as follows:

\begin{definition}(\textbf{Reduced network})\label{Defi:ReducedNet}
Consider a corridor network with patterns of free-flow travel times, bottleneck capacities, and inflow demands $\{\mathbf{c}, \boldsymbol{\mu}, \mathbf{Q} \}$.
A reduced network of the (original) corridor network is the new corridor network constructed by the following methods. 
\begin{enumerate}
\item Unifying the initial and terminal nodes of each false bottleneck into a single node.
\item Aggregating the inflow demands from the initial and terminal nodes (origins), such that the aggregated demand is the inflow demand from the unified node (see, for example, Figure~\ref{Fig:ReducedNetwork}).
\end{enumerate}
If the (original) corridor network does not have false bottlenecks, then the corridor network itself is a reduced network.
\end{definition}
\noindent By definition, a reduced network has no false bottlenecks.
Thus,we can avoid the arbitrariness and complications due to the presence of false bottlenecks.

There is a simple criterion to detect false bottlenecks in a corridor network by utilizing the following concept of \textit{normalized demand}: 
\begin{definition}(\textbf{Normalized demand})
At each bottleneck $i$, the normalized demand $\overline{Q}_{i}$ is defined as follows: 
\begin{align}
&\overline{Q}_{i}\equiv\cfrac{\sum_{j = i}^{n-1}Q_j}{\mu_i-\mu_n}, \label{eq:NorDemand}
\end{align}
where bottleneck $n$ is the closest non-false bottleneck upstream of bottleneck $i$, and $\bar Q_i\equiv \sum_{j = i}^{N}Q_j/\mu_i$ if bottleneck $i$ is the most upstream non-false bottleneck in the corridor network.
\end{definition}
Based on this definition, we can detect false bottlenecks according to the following lemma:

\begin{lemma}\label{Lemma:DetectFBN}
In a corridor network, bottleneck $i$ is non-false iff all downstream bottlenecks of $i$ have a \textit{smaller} normalized demand than that of the bottleneck $i$, i.e., $\forall m<i$, subject to the following:
\begin{align}
&\bar Q_m< \bar Q_i.\label{eq:nonFBcriterion}
\end{align}
\end{lemma}
\noindent This lemma indicates that whether a bottleneck is false or not depends only on the normalized demands of downstream bottlenecks. 
Hence, we can construct a reduced network from any corridor network with arbitrary capacity patterns by using an appropriate algorithm (see \prettyref{app:algo}); this algorithm checks all the downstream bottlenecks from the upstream to downstream of each bottleneck and screens out false bottlenecks, and aggregates travel demands from the upstream and downstream origins of each false bottleneck.
Hereafter, we consider a reduced corridor network constructed from an original corridor network and redefine Index $i\inN$ for bottlenecks and locations in the reduced networks.

\subsubsection{Closed-form solution of the DSO assignment problem}
The analytical solution of the pricing equilibrium of the (reduced) corridor network is presented below.
We first obtain the following lemma that states the property of the arrival time windows: 
\begin{lemma}\label{lemma:CoT}
Consider the pricing equilibrium in a reduced corridor network.
The arrival time window of commuters from an origin includes that of commuters from the downstream origin, as follows:
\begin{align}
\mathcal{T}_{i}\subset \mathcal{T}_{i+1},\quad \forall i\in \mathcal{N}\setminus\{N\}.
\end{align}
\end{lemma}
\noindent This lemma indicates that the arrival time window of commuters departing from the upstream origin is longer than that of commuters departing from the downstream origin.
That is, the arrival time windows have nested structures\footnote{We will later show that this inclusion relationship also holds in dynamic user equilibrium, with the related previous studies.}.

The following lemma can then be obtained, which provides useful insights into the flow pattern of the pricing equilibrium.

\begin{figure*}[t]
	\centering
	\begin{minipage}[t]{0.5\textwidth}
	\centering
		\includegraphics[clip, width=1.0\columnwidth]{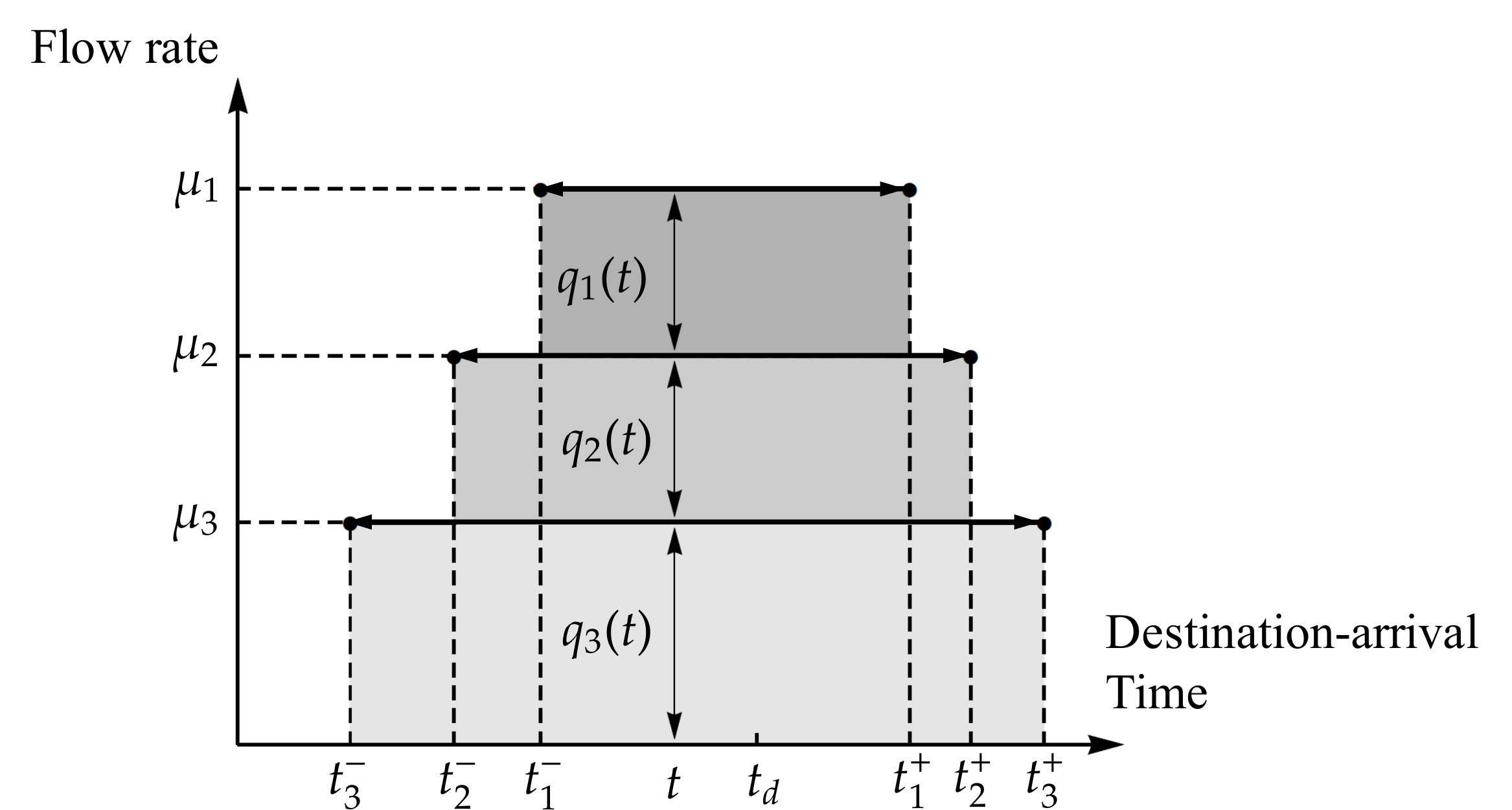}
		\subcaption{Arrival-flow rates at the destination}
		\label{Fig:TNPFlowinCorridor}
	\end{minipage}%
	\begin{minipage}[t]{0.5\textwidth}
	\centering
		\includegraphics[clip, width=1.0\columnwidth]{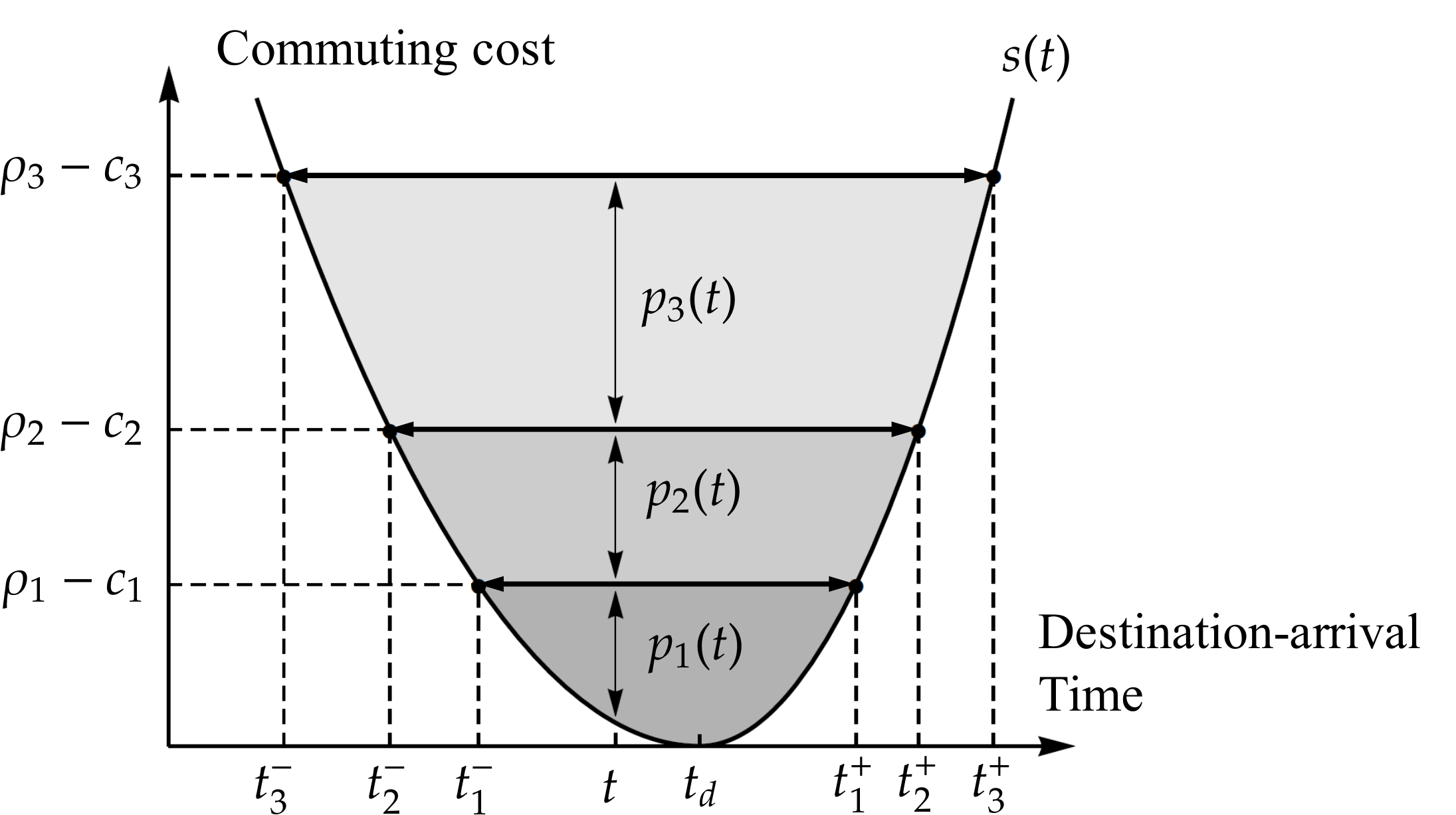}
		\subcaption{Equilibrium commuting costs and optimal prices}
		\label{Fig:TNPCostinCorridor}
	\end{minipage}
	\vspace{0mm}
	\caption{Flow and cost patterns in the pricing equilibrium ($N=3$)}
	\label{Fig:TNPinCorridor}
	\vspace{-2mm}
\end{figure*}

\begin{lemma}\label{lemma:NonnegPrice}
Consider the pricing equilibrium in a reduced corridor network.
The optimum prices on bottleneck $i$ satisfy the following relationship: 
\begin{align}
p_i(t)
\begin{cases}
>0  \quad   &  \quad \text{if}\quad t\in\mathcal{T}_{i},\\
=0  \quad   &  \quad \text{otherwise.}
\end{cases}
\end{align}
\end{lemma}
\noindent This lemma implies that when the destination arrival flow from origin $i$ is positive, the corresponding optimal prices on bottleneck $i$ are positive. 
By combining these lemmas and the capacity constraint condition~\eqref{Eq:DSO_KKT2}, we obtain the following arrival flow pattern in the pricing equilibrium:
\begin{align}
&q_i(t)=
\begin{cases}
\hat{\mu}_i &    \text{if}\quad t\in\mathcal{T}_i,\\
0           &   \text{otherwise},
\end{cases}
 &&    \forall i\in\mathcal{N},	\label{Eq:TNPFlowinCorridor}\\
&\text{where}\quad \hat{\mu}_i\equiv
\begin{cases}
\mu_i-\mu_{i+1}	& \text{for}\quad i\in\mathcal{N}\setminus\{N\},\\
\mu_N					&\text{for}\quad i=N.
\end{cases}
\end{align}
\noindent These equations indicate that the destination arrival flow from each origin becomes an \textit{all-or-nothing pattern} in the pricing equilibrium, as shown in Figure~\ref{Fig:TNPFlowinCorridor}; that is, the inflow $q_{i}(t)$ enters at a constant rate from each origin $i$, such that the departure flow from the immediate downstream bottleneck equals the capacity, or the inflow becomes zero.

After the derivation of the flow pattern in the pricing equilibrium, the remaining problems are the analytical derivation of the earliest and latest arrival times of commuters from each origin and equilibrium commuting costs. 
To achieve this, we employ the length of the arrival time window. 
By combining the flow conservation condition~\eqref{Eq:DSO_KKT3} and \eqref{Eq:TNPFlowinCorridor}, we obtain the following: 
\begin{align}
T_{i}=\cfrac{Q_i}{\hat{\mu}_i},  &&\forall i\in \mathcal{N}. \label{eq:ATWLengthinCorridor}
\end{align}
\noindent By employing the length, we obtain the arrival time window and the equilibrium commuting cost of each origin, as shown in Figure~\ref{Fig:TNPCostinCorridor}.
Given that the optimal prices paid by $i$-commuters arriving at the destination at $t_{i}^{-}$ or $t_{i}^{+}$ (i.e. $t_{i}^{-} + T_i$) are zero, their equilibrium commuting costs, which consist of the free-flow travel times and schedule delay costs, should be the same. 
This indicates that $s(t_{i}^{-}) =  s(t_{i}^{-}+ T_i)$.
Using the definition of $\bar{s}(\cdot)$ in Eq.~\eqref{Eq:Defi_ScheduleDelay}, such a schedule delay is uniquely determined for a given length of the arrival time window $T_i$.
Thus, the equilibrium commuting cost and the earliest/latest arrival times are determined by solving the following equation of the schedule delay: 
\begin{align}
\rho_{i} = \bar{s}(T_{i}) + c_{i},\quad \forall i\in\mathcal{N}.
\end{align}

From the equilibrium commuting costs, the patterns of optimal prices are derived using the equilibrium conditions~\eqref{Eq:DSO_KKT1}, as follows:
\begin{align}
&\sum_{j = 1}^{i}p_j(t) = \rho_i -s(t) -c_i,&\forall i\inN,t\inT_i.\label{Eq:TNPPermit}
\end{align}
This equation suggests that $p_i(t)$ can be derived sequentially from the downstream bottleneck 1 to the upstream bottleneck $N$.

In summary, the following proposition of the closed-form solution is obtained for a reduced corridor network in the pricing equilibrium:
\begin{proposition}\label{prop:AnaSol}
Consider the pricing equilibrium in a reduced corridor network.
The closed-form solution for $i$-commuters $(\forall i\in\mathcal{N})$ is given as follows:
\begin{subequations}
\begin{align}
&\bullet \textrm{arrival-flow rate}: &&q_i(t)=
 \begin{cases}
    \hat{\mu}_i &    \text{if}\quad t\inT_i\\
    0           &   \text{otherwise}
 \end{cases}\label{Eq:TESolution_Flow}\\
&\bullet \textrm{commuting cost}:&&\rho_i=\bar s(T_i)+c_i    
\label{Eq:TESolution_Cost}\\
&\bullet \textrm{optimal prices}:&&p_i(t)=
\begin{cases}
\rho_i-s(t)-c_{i} - \sum_{j = 1}^{i-1}p_j(t)&\text{if}\quad t\inT_i\\
0&\text{otherwise}
\end{cases}\label{Eq:TESolution_Tolls}
\end{align}
\end{subequations}
\end{proposition}
\noindent This is the solution of \textbf{[DSO-LP]}: Eq.~\eqref{Eq:TESolution_Flow} is the optimal arrival flow pattern, and Eqs.~\eqref{Eq:TESolution_Cost} and \eqref{Eq:TESolution_Tolls} correspond to Lagrange multipliers.

\begin{example}
We demonstrate the approach and closed-form solution of the DSO assignment problem in a network with two tandem bottlenecks $(N = 2)$ as an example.
The approach starts with the detection of false bottlenecks based on \textbf{Lemma~\ref{Lemma:DetectFBN}}.
Since the most downstream bottleneck (i.e., bottleneck $1$) is always a non-false bottleneck, the following two cases are exhaustive, wherein (i) bottleneck $2$ is a non-false bottleneck (Case 1), and (ii) bottleneck $2$ is a false bottleneck (Case 2).

In Case 1, since the original corridor network is the reduced network, we can obtain the flow and cost patterns in the DSO assignment from the proposition, as follows:
\begin{align*}
&\bullet \textrm{arrival-flow rate}: 
&&q_{1}(t)=
 \begin{cases}
    \hat{\mu}_{1} &    \text{if}\quad t\in\mathcal{T}_{1}\\
    0           &   \text{otherwise}
 \end{cases},\quad 
 q_{2}(t)=
 \begin{cases}
    \mu_{2} &    \text{if}\quad t\in\mathcal{T}_{2}\\
    0           &   \text{otherwise}
 \end{cases}
 \\
&\bullet \textrm{commuting cost}:
&&\rho_{1}=\bar s(T_{1})+c_1,\quad     \rho_{2}=\bar s(T_{2})+c_1
\\
&\bullet \textrm{optimal prices}:&&
p_{1}(t)=
\begin{cases}
\rho_{1} - s(t) - c_{1}  &\text{if}\quad t\inT_{1}\\
0&\text{otherwise}
\end{cases}, \quad 
p_{2}(t)=
\begin{cases}
\rho_{2} - s(t) - c_{2} - p_{1}  &\text{if}\quad t\inT_{2}\\
0&\text{otherwise}
\end{cases}
\end{align*}

In Case 2, since bottleneck $2$ is a false bottleneck, we first construct the reduced network based on \textbf{Definition~\ref{Defi:ReducedNet}}.
Therefore, the DSO problem reduces to that of the single bottleneck, where $Q_{1} + Q_{2}$ commuters depart from a single origin and pass through the bottleneck with capacity $\mu_{1}$.
Here, we refer to the single origin and bottleneck as origin $r$ and bottleneck $r$, respectively.
We obtain the following optimal flow and cost patterns in the reduced network:
\begin{align*}
&q_{r}(t) = 
\begin{cases}
    \mu_{1} &    \text{if}\quad t\in\mathcal{T}_{r},\\
    0           &   \text{otherwise}.
 \end{cases}\quad 
 \rho_{r} = \bar{s}(T_{r}) + c_{1},\quad 
 p_{r}(t) = 
 \begin{cases}
 	\rho_{r} - s(t) - c_{1}, 	&    	\text{if}\quad t\in\mathcal{T}_{r},\\
 	0									&		\text{otherwise}.
 \end{cases}
 \\
&\text{where}\quad T_{r} = (Q_{1} + Q_{2})/\mu_{1},\quad s(t_{r}^{-}) = s(t_{r}^{+}).
\end{align*}

We then obtain the complete closed-form solution in the original corridor network by disaggregating the solution obtained in the reduced network.
Specifically, we regard the arrival flow rate $q_{r}(t)$ from the reduced origin $r$ as the sum of the arrival flow rates from the aggregated origins in the original corridor network, i.e., $q_{r}(t) = q_{1}(t) + q_{2}(t)$.
In addition, we can regard the equilibrium commuting cost, excluding the free-flow travel time (i.e., $\rho_{r} - c_{1}$) of the commuters departing from the reduced origin as the costs of the commuters from the aggregated origins in the original corridor network.
We obtain the equilibrium commuting cost of commuters from each origin by adding the free-flow travel time (to the destination) to $\rho_{r} - c_{1}$.
This is summarized as follows: 
\begin{align*}
&\bullet \textrm{arrival-flow rate}: 
&&q_{1}(t) + q_{2}(t)=
\begin{cases}
	\mu_{1} &    \text{if}\quad t\in\mathcal{T}_{r}\\
	0           &   \text{otherwise}
\end{cases}\\
& &&
\text{s.t.}\quad q_{2}(t)\leq \mu_{2},\quad \int_{t\in\mathcal{T}_{r}}q_{2}(t) = Q_{2}.
 \\
&\bullet \textrm{commuting cost}:
&&\rho_{1}=\bar s(T_{r})+c_{1},\quad     \rho_{2}=\bar s(T_{r})+c_{2}
\\
&\bullet \textrm{optimal prices}:&&
p_{1}(t)=
\begin{cases}
\rho_{1} - s(t) - c_{1}  &\text{if}\quad t\inT_{r}\\
0&\text{otherwise}
\end{cases}, \quad 
p_{2}(t)=0,\quad \forall t\in\mathcal{T}
\end{align*}
It should be noted that the flow patterns have arbitrariness due to the free-flow condition of bottleneck $2$, whereas the equilibrium commuting costs are uniquely determined, as mentioned in Section~\ref{Sec:MC_DSOFormulation}.
\end{example}

%%%%%%%%%%%%%%%%%%%%%%%%%%%%%%%%%%%%%%%%%%%%%%%%%%%%%%%%%%%%%%%%%%%%
% Section 4
\section{Dynamic user equilibrium assignment for the morning commute}\label{sec:DSOvsDUE}\label{Sec:MC_DUE}
The closed-form solution of the DSO assignment problem is more than a simple contribution toward clarifying characteristics of optimal states.
It theoretically clarifies the relationship between the optimal states and dynamic user equilibrium (DUE) \textit{with queues}.
This relationship allows for the derivation of the closed-form solution in the DUE, and provides a powerful basis for the welfare analysis of several important optimal congestion pricing schemes.

%\sout{This section derives the closed-form solution of the DUE assignment problem for the morning commute with clarifying the relationship of the DUE state with the DSO state.}
In this section, we first describe the settings and formulation of the DUE assignment problem with queues.
Thereafter, we present an analytical approach to solve the DUE assignment problem based on the solution of the DSO assignment derived in the previous section.
Finally, we compare the flow and cost patterns of the DSO and DUE assignment, and derive insights about the first-best and second-best dynamic TDM policies.
%\sout{dynamic optimal congestion pricing schemes} 

\subsection{Settings}\label{Sec:MC_DUE_Settings}
We consider a DUE assignment problem for the morning commute in the corridor network described in Section~\ref{Sec:MC_Networks}.
In this case, each commuter chooses the arrival time of the trip from the origin node, to minimize his/her disutility.
This disutility is composed of the schedule delay cost, free-flow travel time, and queuing delay on every bottleneck that the commuter passes.
The queuing delay is modeled by the point queue model based on the FIFO principle, as mentioned in Section~\ref{Sec:MC_Networks}.

In modeling the DUE assignment accompanying the queuing delay, we describe traffic variables mainly in a \textit{Lagrangian coordinate system}~\citep[][]{akamatsu2015corridor}.
In this system, traffic variables at each origin and bottleneck are expressed as functions of the arrival time at the destination, and not at the origin or bottleneck.
Such an expression is suitable for considering the \textit{ex-post} travel time of each commuter during his/her trip.
Hence, we can easily trace the time-space paths of commuters, which allows for the representation of the equilibrium conditions in a simple manner\footnote{For more details of the Lagrangian coordinate system and the standard Eulerian coordinate system, please refer to \cite{akamatsu2015corridor}}.

\begin{figure}[t]
	\center
	\includegraphics[clip, width=0.85\columnwidth]{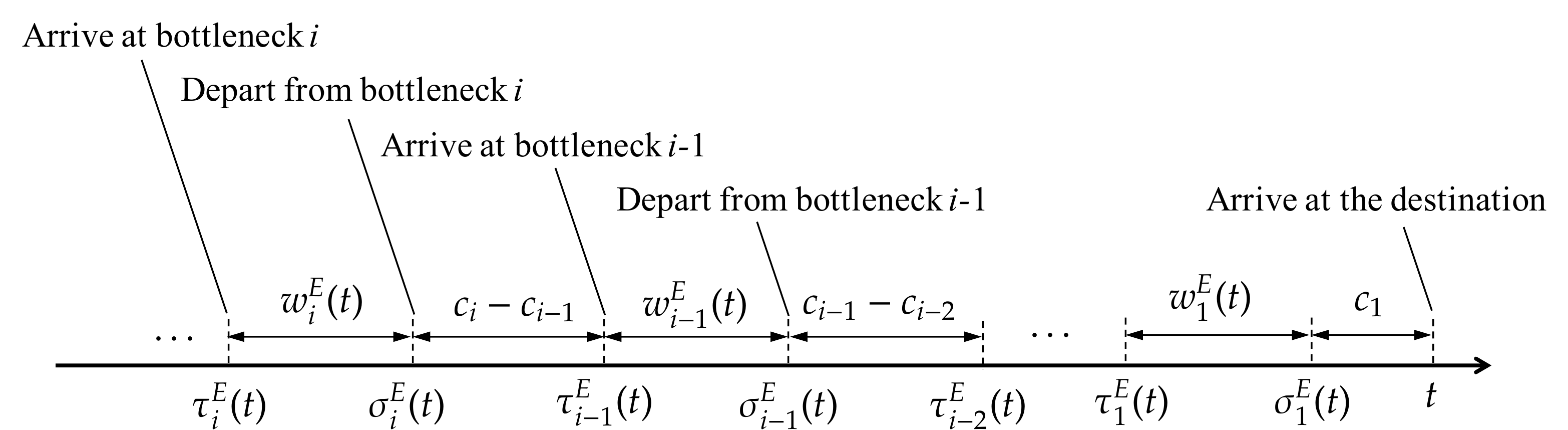}
	\caption{The time-space path of commuters arriving at the destination at time $t$}
	\label{fig:ADTRelation_Morning}
\end{figure}

Several variables are introduced for the expression of the DUE in the Lagrangian coordinate system.
We denote by $\tau_{i}^{E}(t)$ and $\sigma_{i}^{E}(t)$ the equilibrium arrival and departure times at bottleneck $i$ for commuters arriving at the destination at $t$, respectively.
Here, the superscript E indicates that the variables are defined in the DUE problem.
We denote by $w^E_{i}(t)$ the queuing delay at bottleneck $i$ for the commuter.
These variables satisfy the following relationships (see \prettyref{fig:ADTRelation_Morning}):
\begin{align}
&\sigma^E_i(t)\equiv t-\sum_{j = 1}^{i-1}w^E_j(t)-c_i,&\forall i\in \mathcal{N},\forall t\inT,\label{eq:DUE_Sigma}\\
&\tau^E_i(t)\equiv \sigma^E_i(t) - w^E_i(t),&\forall i\in \mathcal{N},\forall t\inT.\label{eq:DUE_Tau}
\end{align}
In addition to these variables representing the time-space paths of commuters, the arrival flow rate $y^E_{i}(t)$ at bottleneck $i$ for commuters arriving at the destination at time $t$ is defined as follows: 
\begin{align}
&y^E_{i}(t)\equiv \dfrac{\mathrm{d}A_i^E(\tau^E_i(t))}{\mathrm{d} t} = \dfrac{\mathrm{d}D_i^E(\sigma^E_i(t))}{\mathrm{d} t} = x^{E}_{i}(\sigma^E_{i}(t))\cdot \dot{\sigma}^E_{i}(t),\label{eq:Definition_DUE_Y}
\end{align}
where $A^E_i(\tau^E_i(t))$ [$D^E_i(\sigma^E_i(t))$] is the cumulative arrival [departure] flow of bottleneck $i$ by time $\tau^{E}_{i}(t)$ [$\sigma^E_i(t)$] under DUE, $x_{i}^{E}(\sigma_{i}^{E}(t))$ is the departure flow rate from bottleneck $i$ at time $\sigma_{i}^{E}(t)$, and $\dot{\sigma}^E_{i}(t)\equiv \Rmd \sigma^E_{i}(t)/\Rmd t$.
Under the FIFO condition, $y_{i}^{E}(t)$ has the following relationship with the destination arrival flow rate of the $i$-commuters at time $t$, $q_{i}^{E}(t)$: 
\begin{align}
&y^E_{i}(t) = \sum_{j = i}^{N} q^E_j(t),  &\forall i\in\mathcal{N},\forall t\inT,\label{eq:DUECond_FC_1}
\end{align}
The generalized transport cost of the $i$-commuters arriving at the destination at time $t$ is expressed as follows:
\begin{align}
v^{E}_{i}(t) = s(t) + c_{i} + \sum_{j = 1}^{i}w_{j}(t).
\end{align}

\subsection{Formulation of dynamic user equilibrium assignment}\label{Sec:MC_DUE_Formulation}
Here, we define the DUE state, in which no commuter can reduce his/her commuting cost by unilaterally changing his/her arrival time.
This DUE assignment problem in a corridor network can be formulated as a linear complementarity problem (LCP) under the following three conditions~\citep[][]{akamatsu2015corridor}.
First, the equilibrium condition for the commuter arrival time choice is expressed as follows: 
\begin{align}
\begin{cases}
\sum_{j = 1}^{i}w_j^E(t)+s(t)+c_i=\rho_i^E&\text{if}\quad  q_{i}^E(t)>0,\\
\sum_{j = 1}^{i}w_j^E(t)+s(t)+c_i\geq\rho_i^E&\text{if}\quad  q_{i}^E(t)=0,
\end{cases}
&&\forall i\in\mathcal{N},\forall t\inT.\label{eq:DUECond_DTCE}
\end{align}
Second, the following flow conservation conditions for all the travel demands should be satisfied: 
\begin{align}
&\int_{t\inT} q_{i}^E(t)\mathrm{d} t = Q_i, &&\forall i\in\mathcal{N}. \label{eq:DUECond_FCMass}
\end{align}
Third, the queueing delay condition at each bottleneck $i$ is expressed as follows: 
\begin{align}
\begin{cases}
\dot{w}^E_{i}(t) = y^E_{i}(t)/\mu_{i} - \dot{\tau}^E_{i}(t) & \text{if}\quad w^E_{i}(t)>0,\\
\dot{w}^E_{i}(t) \geq y^E_{i}(t)/\mu_{i} - \dot{\tau}^E_{i}(t) &\text{if}\quad w^E_{i}(t)=0,
\end{cases}
&&\forall i\in\mathcal{N},\forall t\inT.\label{eq:DUECond_QueueDelay_1proto}
\end{align}
\prettyref{eq:DUECond_FC_1} can be utilized to express this condition, as follows: 
\begin{align}
\begin{cases}
\sum_{j = i}^{N} q_j^E(t)=\mu_{i} \cdot \dot{\sigma}^E_{i}(t) & \text{if}\quad w^E_{i}(t)>0,\\
\sum_{j = i}^{N} q_j^E(t)\leq\mu_{i} \cdot \dot{\sigma}^E_{i}(t) &\text{if}\quad w^E_{i}(t)=0,
\end{cases}
&&\forall i\in\mathcal{N},\forall t\inT.\label{eq:DUECond_QueueDelay_1}
\end{align}

In summary, the DUE in the corridor network can be defined as follows.
\begin{definition}(\textbf{Dynamic user equilibrium})\label{defn:DUE}
Dynamic user equilibrium is a collection of variables $\{\boldsymbol{\rho}^E,\mathbf{q}^E(t),\mathbf{w}^E(t)\}_{t\inT}$ that satisfy Eqs.~\eqref{eq:DUECond_DTCE}, \eqref{eq:DUECond_FCMass}, and \eqref{eq:DUECond_QueueDelay_1}.
\end{definition}

\subsection{Closed-form solution of the DUE assignment problem}\label{Sec:MC_DUE_AnalyticalSolution}

\subsubsection{Underlying principle of the analytical derivation of DUE solution}
The proposed approach for deriving the DUE solution in this study is based on the following conjecture (as will be proved later) that claims the equality between the queuing delays in DUE and the Lagrange multiplier of the DSO assignment problem (i.e., optimal prices in the pricing equilibrium):  
\begin{conjecture}\label{conj:EqualityQD_PP}
There is a DUE solution $\{\boldsymbol{\rho}^E,\mathbf{q}^E(t),\mathbf{w}^E(t)\}_{t\inT}$ in a corridor network of morning commute that satisfies the following equation:
\begin{align}
&w^E_i(t)=p_i(t),&\forall i\in\mathcal{N},\forall t\inT,\label{eq:Equality_QD_PP}
\end{align}
where $\{\mathbf{p}(t)\}_{t\inT}$ is the Lagrange multiplier of \textbf{[DSO-LP]}.
\end{conjecture}
\noindent This indicates that the equilibrium commuting costs are equal in the pricing equilibrium and DUE, as can be observed from the comparison between Eqs.~\eqref{Eq:DSO_KKT1} and \eqref{eq:DUECond_DTCE}; i.e., $\rho_{i} = \rho_{i}^{E}$, $\forall i\in\mathcal{N}$.

This conjecture is actually true in a single-bottleneck model: queuing delays are equal to optimal prices, and equilibrium commuting costs are the same~\citep[for example,][]{akamatsu2021new}.
This implies that the solution of the DUE can be constructed from that of the DSO assignment problem (pricing equilibrium) by appropriately coordinating the arrival flows to the bottleneck, such that queuing delays and the Lagrange multipliers [optimal prices] are the same.
Therefore, the solution of the DUE can be derived from that of the DSO assignment problem in a constructive manner if the conjecture is true in the corridor network.

In the following, we analytically derive the DUE solution by assuming that the conjecture is true; the derivation also serves as a proof of the conjecture.
Specifically, we first derive the departure flows of each bottleneck and destination-arrival flows from each origin by using \prettyref{eq:Equality_QD_PP}.
Thereafter, the derived variables are confirmed as consistent with the DUE conditions.

\subsubsection{Proof of the conjecture}
We first clarify the relationship between the departure flow from bottlenecks under the DSO and DUE states by comparing the optimal pricing condition~\eqref{Eq:DSO_KKT2} and queuing delay condition~\eqref{eq:DUECond_QueueDelay_1}.
By combining~Eqs.~\eqref{Eq:DSO_FlowConservation} and \eqref{Eq:DSO_KKT2}, the optimal pricing condition can be expressed as follows:
\begin{align}
&\begin{cases}
x_i(\sigma_{i}(t))=\mu_i		&\text{if}\quad p_{i}(t)>0,\\
x_i(\sigma_{i}(t))\leq\mu_i	&\text{if}\quad p_{i}(t)=0,
\end{cases}
&\forall i\in\mathcal{N},\forall t\inT,\label{eq:TNP_PPE_2}
\end{align}
\noindent In contrast, the queuing delay condition can be expressed by combining Eqs.~\eqref{eq:DUE_Tau}, \eqref{eq:Definition_DUE_Y}, and \eqref{eq:DUECond_QueueDelay_1proto} as follows: 
\begin{align}
\begin{cases}
x^{E}_{i}(\sigma_i^E(t))=\mu_i      &\text{if}\quad w^E_{i}(t)>0,\\
x^{E}_{i}(\sigma_i^E(t))\leq\mu_i   &\text{if}\quad w^E_{i}(t)=0,
\end{cases}
&&\forall i\in \mathcal{N},\forall t\inT.\label{eq:DUECond_QueueDelay_2}
\end{align}
These conditions imply that the \textit{departure flow patterns of the bottlenecks in both problems are the same} when the queuing delays are equivalent to the optimal prices.
By combining this relationship and the flow conservation condition~\eqref{eq:DUECond_FCMass}, we can formally obtain the following lemma: 
\begin{lemma}\label{lemma:Equality_OutFlow}
If Conjecture~\ref{conj:EqualityQD_PP} is true, the departure flows of bottlenecks in the DUE and DSO assignment have the following relationship: 
\begin{align}
&x_i^E(\sigma^E_i(t))=x_i(\sigma_{i}(t)),&\forall i\in\mathcal{N},\forall t\inT.\label{eq:Equality_OutFlow}
\end{align}
\end{lemma}
\noindent
\begin{proof}
See \ref{app:proofx=x}.
\end{proof}

This relationship enables us to derive the destination arrival flow rates in DUE from the solution of the DSO assignment.
From Eqs.~\eqref{eq:Definition_DUE_Y} and \eqref{eq:DUECond_FC_1}, we obtain the following:
\begin{align}
q^E_{i}(t)&=y^E_{i}(t) - y^E_{i+1}(t) \notag\\
&= x^{E}_{i}(\sigma^E_{i}(t))\cdot \dot{\sigma}^E_{i}(t) - x^{E}_{i+1}(\sigma^E_{i+1}(t))\cdot \dot{\sigma}^E_{i+1}(t),	\quad\forall i\in\mathcal{N}\setminus\{N\},\forall t\in\mathcal{T},\label{eq:DUEqUsingx}\\
q^E_{N}(t) &=   x^{E}_{N}(\sigma^E_{N}(t))\cdot \dot{\sigma}^E_{N}(t),\qquad \forall t\in\mathcal{T}.\label{eq:DUEqNUsingx}
\end{align}
Furthermore, Eqs.~\eqref{Eq:TESolution_Tolls}, \eqref{eq:DUE_Sigma}, and \eqref{eq:Equality_QD_PP} are combined to derive the following:
\begin{align}
\dot{\sigma}_{i}^{E}(t) &= 1 - \sum_{j = 1}^{i-1}\dot{w}_{j}^{E}(t) = 1 - \sum_{j = 1}^{i-1}\dot{p}_{j}^{E}(t) \\
& = \begin{cases}
1+\dot s(t)&\mathrm{if}\quad t\inT_{i-1},\\
1&\mathrm{otherwise},
\end{cases}&\forall i\inN.\label{eq:dotsigma}
\end{align}
Subsequently, by substituting this equation into Eqs.~\eqref{eq:DUEqUsingx} and \eqref{eq:DUEqNUsingx}, $q_{i}^{E}(t)$ is derived as follows.
\begin{align}
& q_i^E(t) =
\begin{cases}
(1 + \dot s(t) )  \cdot\hat\mu_i			& \text{if}\quad t\in\mathcal{T}_{i-1}, \\
\hat\mu_i - \dot s(t)\cdot\mu_{i+1}	& \text{if}\quad t\in\mathcal{T}_i\setminus\mathcal{T}_{i-1},\\
0															& \text{otherwise,}
\end{cases}
&\forall i\inN\setminus\{N\},\label{eq:DUE_ArrivalFlow_1}\\
& q_{N}^{E}(t)  = 
\begin{cases}
(1 + \dot s(t) )  \cdot \mu_N			& \text{if}\quad t\in\mathcal{T}_{N-1}, \\
\mu_N	& \text{if}\quad t\in\mathcal{T}_N\setminus\mathcal{T}_{N-1},\\
0															& \text{otherwise.}
\end{cases}\label{eq:DUE_ArrivalFlow_2}
\end{align}
This relationship means that the destination arrival time windows $\{\ClT_i\}$ in the pricing equilibrium and those in the DUE are the same.
Letting $\mathcal{T}_{i}^{E}$ be the arrival time window of $i$-commuters in the DUE; then $\mathcal{T}_{i}=\mathcal{T}_{i}^{E}$, $\forall i\in\mathcal{N}$.
As a result, all the variables in the DUE can be analytically derived from the solution of the DSO assignment based on the assumption that the conjecture is true.

The derived variables satisfy almost all the equilibrium conditions in the DUE because the variables are induced by the conditions themselves.
However, the derivation may not guarantee the \textit{non-negativity constraint} on the arrival flow $q_{i}^{E}(t)$ from each origin: the cumulative arrival flow curve at an origin may be backward bending.
Thus, the conditions in the following lemma should be satisfied to ensure a physically appropriate arrival flow.

\begin{lemma}\label{lemma:DUE_Feasible}
A DUE solution under \textbf{\prettyref{conj:EqualityQD_PP}} exists if the following conditions are satisfied:
\begin{subequations}\label{Eq:DUE_Feaibility}
\begin{align}
& -1	\leq \dot{s}(t), & \forall t\inT_{N},\label{eq:DUE_Feasible_1}\\
&\dot{s}(t) \leq \cfrac{\mu_i}{\mu_{i+1}}-1, & \forall i\in\mathcal{N}\setminus\{N\},\forall t\inT_{i}\setminus\mathcal{T}_{i-1},\label{eq:DUE_Feasible_2}
\end{align}
\end{subequations}
where $\ClT_0\equiv\emptyset$.
\end{lemma}
\noindent Condition~\eqref{eq:DUE_Feasible_1} corresponds to the existence condition of DUE in a corridor network~\citep{akamatsu2015corridor}, whereas \prettyref{eq:DUE_Feasible_2} guarantees the existence of the DUE solution, where the pattern of queuing delays is equal to that of dynamic prices in pricing equilibrium.
This indicates that a DUE solution exists when Condition~\eqref{eq:DUE_Feasible_1} is satisfied.
However, the bottleneck departure flows and cost patterns are different from those of the pricing equilibrium when Condition~\eqref{eq:DUE_Feasible_2} is not satisfied (as demonstrated later in numerical examples).

The following proposition summarizes the above discussion. 
\begin{proposition}\label{prop:DUEAnaSol}
Consider the DUE problem for the morning commute in a corridor network.
Assuming that Eqs.~\eqref{eq:DUE_Feasible_1} and \eqref{eq:DUE_Feasible_2} are satisfied, the closed-form solution for $i$-commuters $(\forall i\in\mathcal{N})$ is expressed as follows:
\begin{subequations}
\begin{align}
&\bullet \textrm{arrival-flow rate}: && q_i^E(t) =
\begin{cases}
(1 + \dot s(t) )  \cdot\hat\mu_i			& \text{if}\quad t\in\mathcal{T}_{i-1}, \\
\hat\mu_i - \dot s(t)\cdot\mu_{i+1}	& \text{if}\quad t\in\mathcal{T}_i\setminus\mathcal{T}_{i-1},\\
0															& \text{otherwise,}
\end{cases}
&\forall i\inN\setminus\{N\},\label{eq:DUESolution_Flow}\\
& && q_{N}^{E}(t)  = 
\begin{cases}
(1 + \dot s(t) )  \cdot \mu_N			& \text{if}\quad t\in\mathcal{T}_{N-1}, \\
\mu_N	& \text{if}\quad t\in\mathcal{T}_N\setminus\mathcal{T}_{N-1},\\
0															& \text{otherwise.}
\end{cases}\label{eq:DUESolution_Flow_2}\\
&\bullet \textrm{commuting cost}: &&\rho^E_i=\bar s(\dfrac{Q_i}{\hat\mu_i})+c_i&&\label{eq:DUESolution_Cost}\\
&\bullet \textrm{queuing delay}: &&w_i^E(t)=
\left\{
\begin{array}{ll}
\rho^E_i-s(t)- c_{i} - \sum_{j = 1}^{i-1}w_j^E(t)&\mathrm{if}\quad t\inT_i\\
0&\mathrm{otherwise}
\end{array}
\right.&&\label{eq:DUESolution_Delay}
\end{align}
where $\ClT_0\equiv \emptyset$.
\end{subequations}
\end{proposition}

\begin{figure}[t]
\center
\includegraphics[width=0.8\columnwidth,clip]{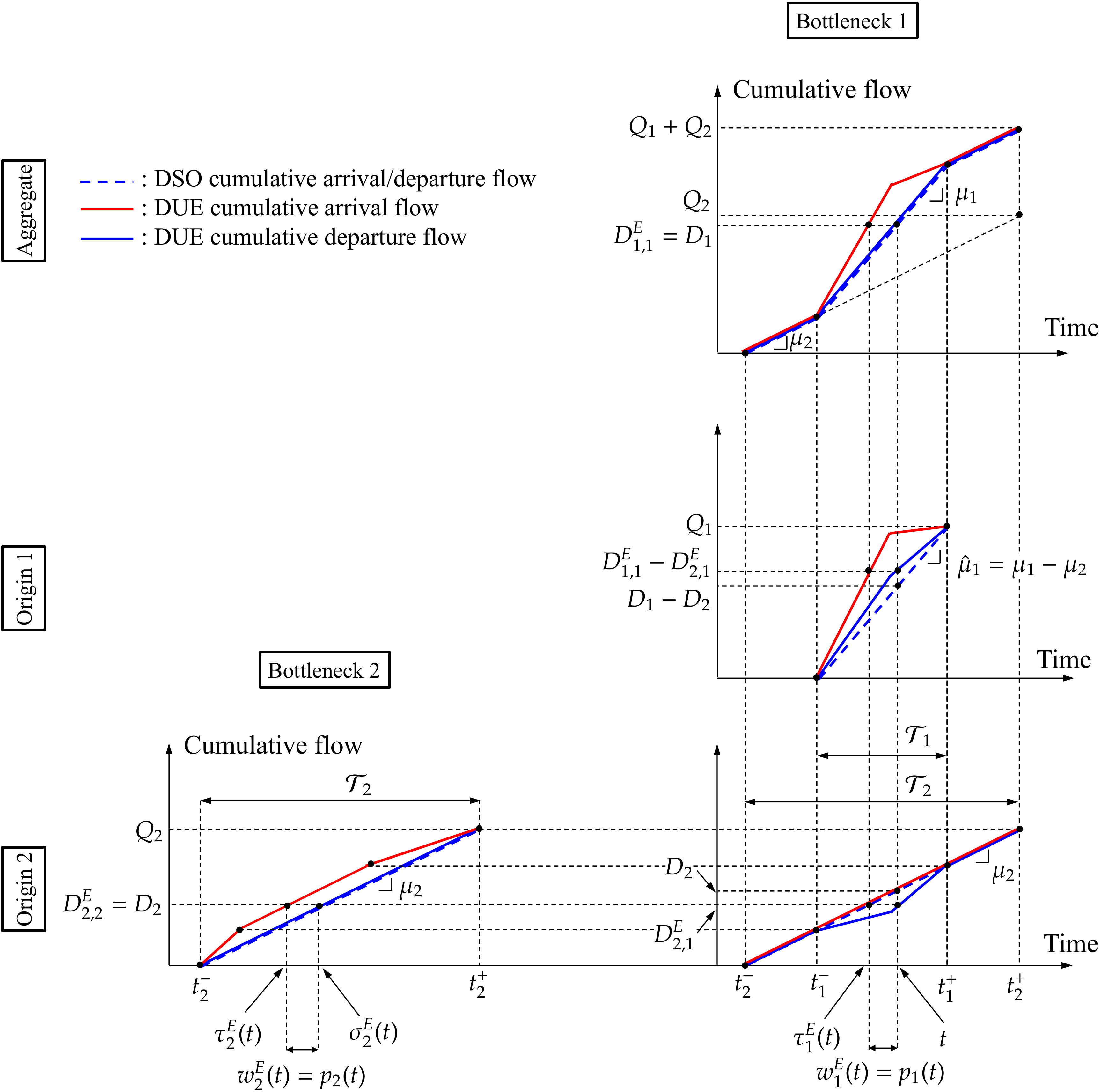}
\caption{DUE cumulative flow curves of bottlenecks for the morning commute in a corridor network ($N=2$). Note that $D_{i,m}^E\equiv D_i^E(\sigma^E_m(t))$.}
\label{Fig:CumulativeCurves_DUEDSO}
\end{figure}

The derivation of the DUE solution from the DSO solution in a corridor network with a simple setting is demonstrated.
Here, we depict the cumulative flow curves of the bottlenecks in a corridor network with $N = 2$, a piecewise linear schedule delay function, and zero free-flow travel times (i.e., $c_{1} = c_{2} = 0$) in Figure~\ref{Fig:CumulativeCurves_DUEDSO}.
The cumulative flow curves are obtained from the solution of \textbf{[DSO-LP]} by the following steps:
(i) The cumulative departure curve for each bottleneck $i$ in the DSO state $D_{i}(t)$ is plotted.
We obtain the disaggregate cumulative departure curves for each origin $i$, which represent the cumulative number of departures of $i$-commuters at each bottleneck, by $D_{i}(t)-D_{i+1}(t)$.
(2) The cumulative departure flow curve of bottleneck $i$ in the DUE $D^E_i(\sigma_i^E(t))$, which is the same as that of $D_i(t)$, as expressed by \textbf{Lemma~\ref{lemma:Equality_OutFlow}}, si plotted.
We also obtain the disaggregate cumulative departure curves for origin $i$ in the DUE by $D^E_i(t)-D^E_{i+1}(t)$.
(3) The cumulative arrival curve $A^E_i(\tau_i^E(t))$ of each bottleneck $i$ in the DUE is plotted by shifting the departure curve, such that the horizontal distance between the departure and arrival curves is equal to the Lagrange multiplier $p_i(t)$ (as queuing delays $w_i^E(t)$ are equal to the Lagrange multiplier).

In this manner, the cumulative arrival and departure curves can be derived based on the DSO solution using a systematic approach.
Moreover, as mentioned above, we can also derive a closed-form solution of the arrival time windows of commuters from each origin, as follows: 
\begin{align}
\mathcal{T}_{i}^{E} = \mathcal{T}_{i} = \cfrac{Q_{i}}{\hat{\mu}_{i}},\quad \forall i\in\mathcal{N}.
\end{align}
This equation means that each arrival time window in the DUE state is the same as that in the DSO state.
Thus, the arrival time windows in the DUE state have the nested structures similar to those in the DSO state.
The more distant the origin of commuters is, the earlier the commuters start to arrive at the destination and the later the commuters finish arriving at the destination.

It should be noted that the properties of such nested structures of destination-arrival time windows seems to be robust in departure/arrival time choice problems in single-destination corridor networks, even if the traffic flow model (i.e., congestion technology) is different.
For example, \cite{Tian2007} considered the departure time choice DUE assignment in a single-destination public transit system with an in-vehicle crowding effect, and revealed nested departure time windows in a discrete time setting for discrete train schedules.
In addition, \cite{arnott2011corridor} employed the Lighthill-Whitham-Richards (LWR) traffic flow model in a corridor network with continuum-entry points, and proposed a horn-shaped equilibrium departure-time set from origins, which may result in the nested structures of arrival time windows at the destination.

\begin{figure}[t]
	\begin{minipage}[t]{0.49\columnwidth}
	\centering
		\includegraphics[clip, width=1.0\columnwidth]{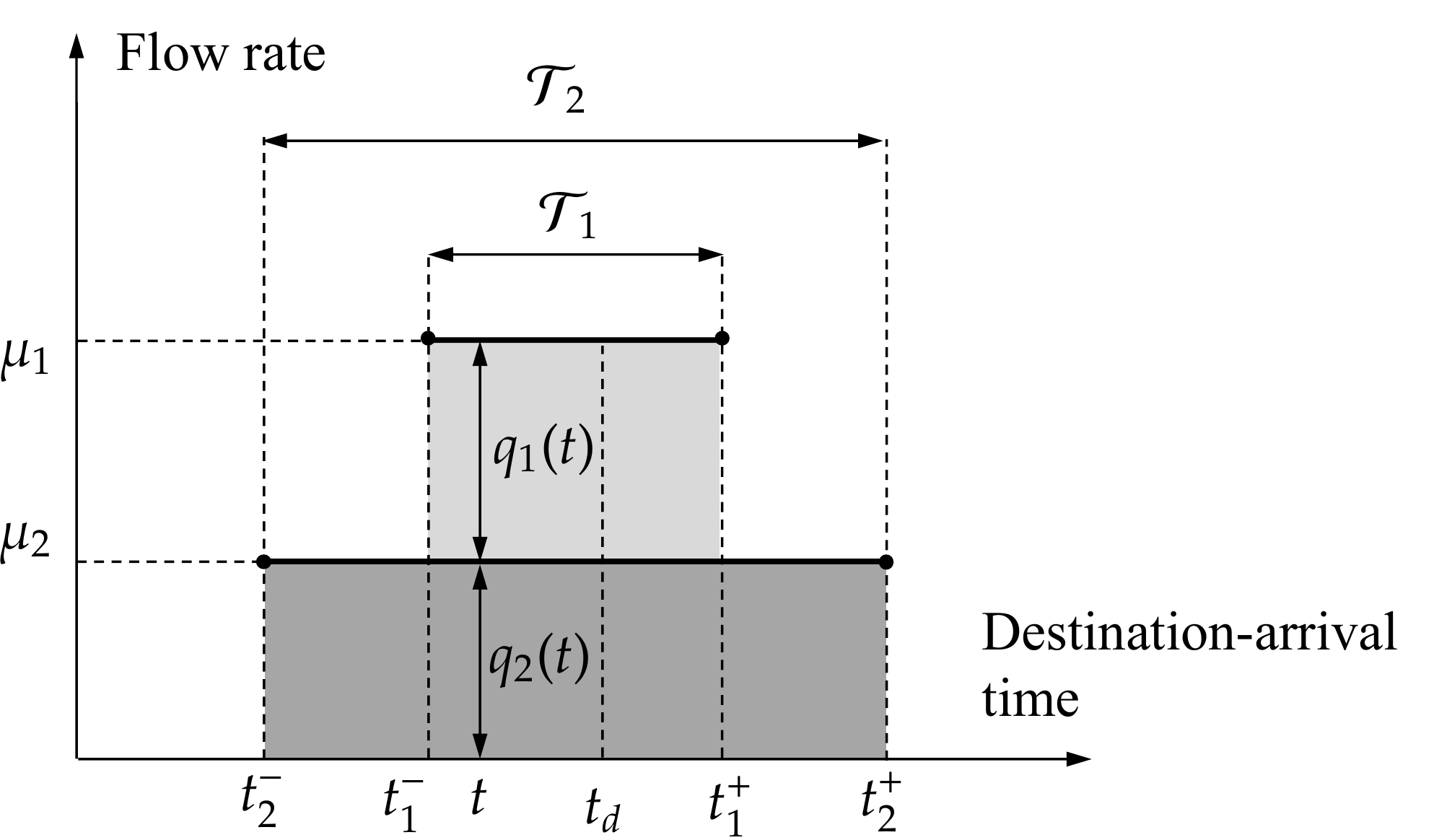}
		\subcaption{Arrival flow rate in the DSO state}
		\label{Fig:CompareMC_DSOFP}
	\end{minipage}
	\begin{minipage}[t]{0.49\columnwidth}
	\centering
		\includegraphics[clip, width=1.0\columnwidth]{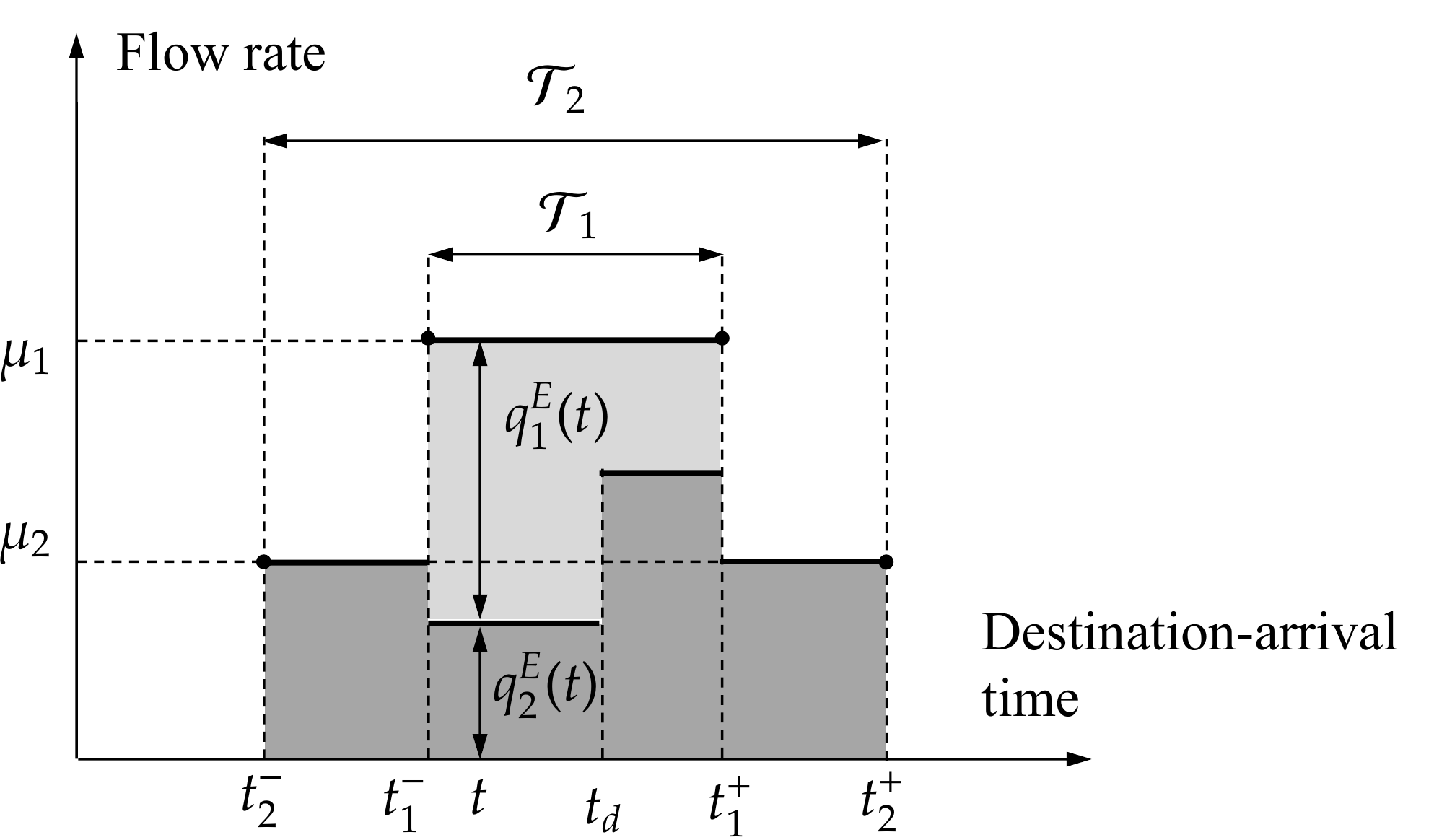}
		\subcaption{Arrival flow rate in the DUE state}
		\label{Fig:CompareMC_DUEFP}
	\end{minipage}\vspace{0mm}
 \caption{Destination-arrival flow rates of the DSO and DUE solutions for the morning commute}
\label{Fig:CompareMC_FlowPatterns}
\end{figure}

\subsection{Optimum vs. Equilibrium}
Now, we compare the flow and cost patterns of the DSO and DUE states, and present several remarks on the relationships.
First, Figure~\ref{Fig:CompareMC_FlowPatterns} presents the difference between the destination arrival flow rates in the DSO and DUE states in a corridor network with tandem-bottlenecks ($N = 2$) as an example. 
As can be seen from the figure, the total arrival flows in the DSO state and those in the DUE state are the same, i.e., $\sum_{i\in\mathcal{N}}q_{i}(t) = \sum_{i\in\mathcal{N}}q^{E}_{i}(t)$.
Meanwhile, the disaggregate (origin-specific) arrival flows in the DUE state $q_{i}^{E}(t)$ are different from those in the DSO state due to the congestion effect.
Specifically, when $t\in[t_{1}^{-}, t_{d}]$ (i.e., commuters arrive at the destination earlier than the desired arrival time), the arrival flow rates from the downstream origin $1$ in the DUE state are larger than those in the DSO state.
However, when $t\in[t_{d}, t_{1}^{+}]$, the DUE arrival flow rates are smaller than the DSO arrival flow rates.
Such differences in the disaggregate arrival flows in the optimal and (user) equilibrium states are characteristic of corridor problems, which do not occur in a single bottleneck problem.

Despite the differences in the disaggregate arrival flow rates, the equilibrium commuting costs from every origin in the DUE and DSO states are the same (when Condition~\eqref{Eq:DUE_Feaibility} is satisfied), as confirmed in the previous section. 
Considering the equivalence of the DSO assignment and the pricing equilibrium analyzed in Section~\ref{Sec:MorningCommute}, this implies that a \textit{Pareto improvement can be achieved} by an appropriate first-best TDM policy (such as dynamic pricing and tradable permits schemes) that imposes optimal prices equal to queueing delays.
Specifically, such first-best TDM policy does not increase the equilibrium costs of all the commuters.
In addition, the road manager who implements the policy can gain an income equal to the sum of the queuing delays in the DUE.
%Thus, the social transport cost, which composes of the free-flow travel times and schedule delay costs of all the commuters, decreases without increasing the cost of each commuter or road manager.
This is formally summarized in the following theorem:
\begin{theorem}
Consider a DUE state in a corridor network for the morning commute that satisfies Eqs.~\eqref{eq:DUE_Feasible_1} and \eqref{eq:DUE_Feasible_2}.
By imposing dynamic prices equal to queuing delays in DUE on every bottleneck, the social transport cost, which is composed of the travel times and schedule delay costs, is minimized.
Moreover, a Pareto improvement is achieved.
\end{theorem}
%\noindent Note that this social transport cost does not include the prices paid by commuters although these prices are included in commuting costs of commuters.
%This is because they are just income transfers between the commuters and a road manager who collects the prices.
%Hence, when the existence conditions are satisfied, the road manager can gain the income equal to the sum of the queuing delays in the DUE with achieving a Pareto improvement,.
%This implies that we can interpret Eq.~\eqref{eq:DUE_Feasible_2} as a \textit{sufficient condition} for a Pareto improvement in corridor problems.

The equivalence between the patterns of queuing delays and optimal prices implies that a Pareto improvement can be achieved by a \textit{second-best} TDM policy; i.e., the social transport cost decreases even if we \textit{partially} implement a TDM policy on several bottlenecks, as follows:
\begin{corollary}
Consider a DUE state in a corridor network for the morning commute that satisfies Eqs.~\eqref{eq:DUE_Feasible_1} and \eqref{eq:DUE_Feasible_2}.
By imposing dynamic prices equal to queuing delays in the DUE on \textit{several bottlenecks}, the social transport cost is always improved without increasing the commuting cost of every commuter.
\end{corollary}
\noindent 
This implies that (i) we can implement optimal TDM policies on each bottleneck individually, and (ii) the dynamic prices in the first-best TDM policy correspond to those in the second-best policy.
Thus, we can improve the social transport cost in the corridor network by partially implementing dynamic pricing on bottlenecks without considering the influence of such implementation on the other bottlenecks.
In summary, we can interpret Condition~\eqref{Eq:DUE_Feaibility} as a sufficient condition for a Pareto improvement by first-best/second-best TDM policies in corridor problems.

%%%%%%%%%%%%%%%%%%%%%%%%%%%%%%%%%%%%%%%%%%%%%%%%%%%%%%%%%%%%%%%%%%%%
% Section 5

\section{Corridor problems for the evening commute}\label{Sec:EveningCommute}
This section considers the DSO and DUE corridor problems for the evening commute, in which the corridor network defined in Section~\ref{Sec:MC_Networks} is used by commuters in reverse.
The specific settings are described in Section~\ref{Sec:EC_Networks}.
Section~\ref{Sec:EC_DSO} shows a DSO assignment problem for the evening commute and derives the analytical solution of the problem.
In Section~\ref{Sec:EC_DUE}, we derive the analytical solution of a DUE assignment problem for the evening commute from that of the DSO assignment using the same approach used in Section~\ref{Sec:MC_DUE_AnalyticalSolution}.

\begin{figure}[t]
 \centering
 \includegraphics[width=80mm,clip]{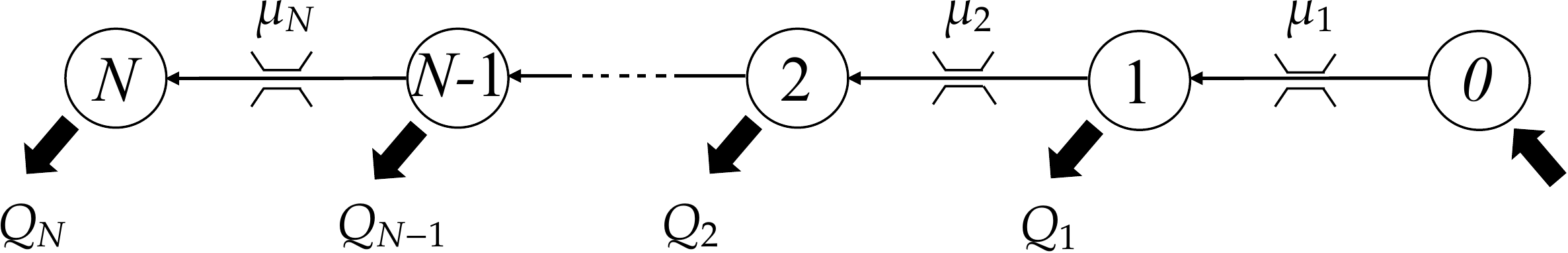}
 \caption{The evening commute in a corridor network with $N$ off-ramps and $N$ bottlenecks}
\label{Fig:EC_CorridorNetwork}
\end{figure}

\subsection{Networks and settings}\label{Sec:EC_Networks}
In corridor problems for the evening commute, we consider a corridor network with a single on-ramp (origin node) and $N$ off-ramps (destination nodes) (see Figure~\ref{Fig:EC_CorridorNetwork}).
These nodes are numbered sequentially from origin node $0$ to the most distant destination node $N$. 
Here, the set of destination nodes is denoted by $\mathcal{N}$ here.
We refer to the bottleneck immediate upstream of the destination $i\in\mathcal{N}$ as bottleneck $i$.
For the bottleneck and origin $i$, the notations of the finite capacity, free-flow travel time, cumulative arrival (departure) flow, and arrival (departure) flow rates are the same as the corridor problems for the morning commute.

From the unique origin, $Q_{i}$ commuters enter the network and reach the destination $i$ during the evening rush hour $\mathcal{T}\equiv [0, T]$ (we refer to the commuters as `$i$-commuters').
Each commuter receives the disutility associated with the travel time from the origin to the destination and the schedule delay cost in the same manner as the morning commute problems.
The only difference is that the evening commuters experience a schedule delay cost based on the deviation between the desired departure time and actual departure time at the origin.
We denote by $t_d$ the desired departure time of all commuters, and denote by $s(t)$ the schedule delay function for commuters departing from the origin at time $t$.

\subsection{Dynamic system optimal assignment}\label{Sec:EC_DSO}
Here, we consider a DSO assignment problem for the evening commute, where the total transport cost is minimized without any queuing delay.
We first re-denote by $q_{i}(t)$ the departure flow rate of the $i$-commuters from the origin at (departure) time $t$.
It is clear that the DSO assignment problem has the same mathematical structure as [DSO-LP] based on the description of the variables in Section~\ref{Sec:MC_DSOFormulation}: the total schedule and (free-flow) travel costs are expressed by Eq.~\eqref{Eq:DSOLP_Obj}, the flow conservation conditions are expressed by Eq.~\eqref{Eq:DSOLP_Const1}, and the capacity constraints are expressed by Eq.~\eqref{Eq:DSOLP_Const2}.

Since the optimization problem has the same structure, the optimality conditions, namely, the corresponding pricing equilibrium conditions under the optimal dynamic congestion pricing scheme, are the same as follows: 
\begin{align}
&\begin{cases}
s(t) + c_{i} + \sum_{j = 1}^{i}p_{j}(t) = \rho_{i}	\quad	&\text{if}\quad q_{i}(t) > 0,\\
s(t) + c_{i} + \sum_{j = 1}^{i}p_{j}(t) \geq \rho_{i}	\quad	&\text{if}\quad q_{i}(t) = 0,
\end{cases}
&\quad \forall i\in\mathcal{N},\forall t\in\mathcal{T},\label{Eq:DSO_EC_KKT1}\\
&\begin{cases}
\sum_{j = i}^{N}q_{j}(t) = 		\mu_{i}	\quad	&\text{if}\quad p_{i}(t) > 0,\\
\sum_{j = i}^{N}q_{j}(t) \leq \mu_{i}	\quad	&\text{if}\quad p_{i}(t) = 0,
\end{cases}
&\quad \forall i\in\mathcal{N},\forall t\in \mathcal{T},\label{Eq:DSO_EC_KKT2}\\
&\int_{t\inT}q_i(t)\Rmd t=Q_i  &\forall i\in \mathcal{N}.\label{Eq:DSO_EC_KKT3}
\end{align}
Thus, Eqs.~\eqref{Eq:TESolution_Flow}-\eqref{Eq:TESolution_Tolls} for the morning commute yield the closed-form solution of the DSO (and the pricing equilibrium) for the evening commute, although the definitions of the variables between the two problems are different.
This result is summarized by the following proposition.
\begin{proposition}\label{prop:DSOAnaSolEve}
Consider the DSO problem for the evening commute in a reduced corridor network. The analytical solution for $i$-commuters $(\forall i\inN)$ is obtained as Eqs.\eqref{Eq:TESolution_Flow}-\eqref{Eq:TESolution_Tolls}.
\end{proposition}

Before proceeding to the following section, it should be noted that there is a slight difference from the morning commute problem in the relationship between $q_{i}(t)$ and $x_{i}(\sigma_{i}(t))$, departure flow rates on bottleneck $i$ at time $\sigma_{i}(t)$, as follows: 
\begin{align}
&q_{i}(t) = 
\begin{cases}
x_{i}(\sigma_{i}(t)) - x_{i+1}(\sigma_{i+1}(t)),	&\quad \forall i\in\mathcal{N}\setminus \{ N\},\\
x_{N}(\sigma_{N}(t)), 									&\quad i = N
\end{cases}
\quad \forall t\in\mathcal{T},\label{Eq:EC_DSO_FC}	\\
&\text{where}\quad \sigma_{i}(t) = t + c_{i}.
\end{align}
Here, $\sigma_{i}(t)$ represents the departure time from the bottleneck $i$ of commuters with an \textit{origin departure time of} $t$, and this is not represented in the same manner as Eq.~\eqref{Eq:DSO_Sigma}.
In other words, the relationship between $q_{i}(t)$ and $x_{i}(\sigma_{i}(t))$ exhibits a symmetric structure to the morning commute problem.

\begin{figure}[t]
\center
 \includegraphics[width=1.0\columnwidth,clip]{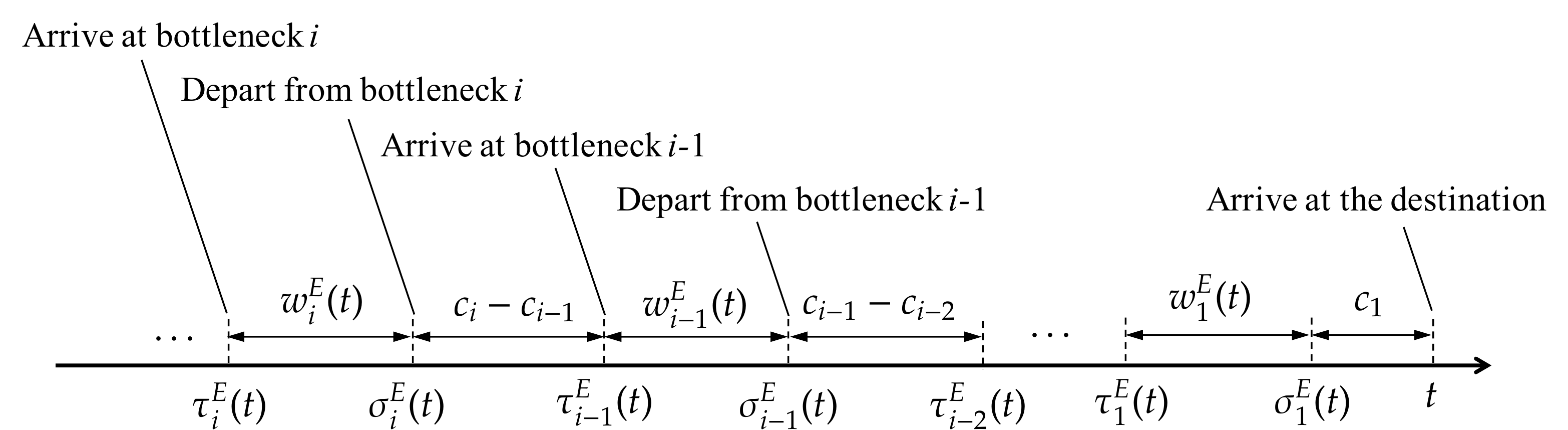}
\caption{Time-space path of commuters departing from the origin at $t$}
\label{Fig:EC_TimeSpacePath}
\end{figure}

\subsection{Dynamic user equilibrium assignment}\label{Sec:EC_DUE}

\subsubsection{Settings and formulation}
We next consider a DUE assignment problem (with queuing delay) for the evening commute, where each commuter chooses the departure time from the origin to minimize his/her disutility.
In modeling the DUE assignment problem, we employ the Lagrangian coordinate system, wherein traffic variables are expressed as functions of the origin departure time.
Specifically, we denote by $\tau_{i}^{E}(t)$ and $\sigma^E_i(t)$ the equilibrium arrival and departure times at bottleneck $i$ for commuters departing from the origin at time $t$, respectively.
We denote by $w_{i}^{E}(t)$ the queuing delay of the commuters at bottleneck $i$.
In the evening commute problem, these variables satisfy the following relationship (see Figure~\ref{Fig:EC_TimeSpacePath}):
\begin{align}
&\sigma^E_i(t)\equiv t+\sum_{j = 1}^{i}w^E_j(t)+c_i, &\forall i\inN,\forall t\inT.\label{eq:sigmaEve}
\end{align}
The arrival flow rate at bottleneck $i$ for commuters with origin departure time $t$ is denoted by $y_{i}^{E}(t)$, and expressed in the same manner as the morning commute problem, as follows:
\begin{align}
&y^E_{i}(t)\equiv \dfrac{\mathrm{d}A_i^E(\tau^E_i(t))}{\mathrm{d} t} = \dfrac{\mathrm{d}D_i^E(\sigma^E_i(t))}{\mathrm{d} t} = x^{E}_{i}(\sigma^E_{i}(t))\cdot \dot{\sigma}^E_{i}(t).\label{eq:Definition_DUEEve_Y}
\end{align}
Note that $y_{i}^{E}(t)$ satisfies the following relationship under the FIFO condition:
\begin{align}
y^E_{i}(t) = \sum_{j = i}^{N} q^E_j(t),  &\forall i\in\mathcal{N},\forall t\inT,\label{Eq:EC_DUE_FC1}
\end{align}
where $q_{i}^{E}(t)$ is the origin departure flow rate of the $i$-commuters at $t$.

The DUE problem can be formulated as a linear complementarity problem (LCP) consisting of the following three conditions~\citep{akamatsu2015corridor}. 
First, the equilibrium condition for commuters departure time choice is expressed as follows
\begin{align}
&\begin{cases}
\sum_{j = 1}^{i}w_j^E(t)+s(t)+c_i=\rho_i^E&\text{if}\quad  q_{i}^E(t)>0,\\
\sum_{j = 1}^{i}w_j^E(t)+s(t)+c_i\geq\rho_i^E&\text{if}\quad  q_{i}^E(t)=0,
\end{cases}
&&\forall i\inN,\forall t\inT.\label{eq:DUEEveDTCE}
\end{align}
Second, the following flow conservation conditions regarding all travel demands should be satisfied:
\begin{align}
&\int_{t\inT} q_{i}^E(t)\mathrm{d} t = Q_i, &&\forall i\in\mathcal{N}. \label{eq:DUEEveFCMass}
\end{align}
Third, the queuing delay condition at each bottleneck $i$ is expressed in the same manner as the morning commute problem, as follows:
\begin{align}
\begin{cases}
\sum_{j = i}^{N} q_j^E(t)=\mu_{i} \cdot \dot{\sigma}^E_{i}(t) & \text{if}\quad w^E_{i}(t)>0,\\
\sum_{j = i}^{N} q_j^E(t)\leq\mu_{i} \cdot \dot{\sigma}^E_{i}(t) &\text{if}\quad w^E_{i}(t)=0,
\end{cases}
&&\forall i\in\mathcal{N},\forall t\inT.\label{eq:DUEEveQueue}
\end{align}

The DUE for the evening commute is summarized as follows.
\begin{definition}(\textbf{Dynamic user equilibrium for the evening commute})
The dynamic user equilibrium for the \textit{evening commute} is a collection of variables $\{\boldsymbol{\rho}^E,\mathbf{q}^E(t),\mathbf{w}^E(t)\}_{t\inT}$ that satisfy Eqs.~\eqref{eq:DUEEveDTCE}, \eqref{eq:DUEEveFCMass}, and \eqref{eq:DUEEveQueue}.
\end{definition}

\subsubsection{Closed-form solution for the evening commute}
Although the DUE problems for the morning and evening commute exhibit different formulations, the same approach can be applied to obtain the solutions.
Specifically, we suppose that queuing delays of the DUE solution are equal to the optimal prices of the pricing equilibrium (i.e., the Lagrange multipliers of the DSO assignment) for the evening commute:
\begin{conjecture}\label{conj:EveEqualityQD_PP}
In the evening commute, there is a DUE solution $\{\boldsymbol{\rho}^E,\mathbf{q}^E(t),\mathbf{w}^E(t)\}_{t\inT}$ in a corridor network that satisfies the following: 
\begin{align}
&w^E_i(t)=p_i(t),&\forall i\in\mathcal{N},\forall t\inT,\label{eq:EveEquality_QD_PP}
\end{align}
where $\{\mathbf{p}(t)\}_{t\inT}$ is the Lagrange multiplier of the DSO assignment.
\end{conjecture}

\noindent Note that the equality between equilibrium commuting costs in the pricing equilibrium and DUE (i.e., $\boldsymbol{\rho}=\boldsymbol{\rho}^E$) is obtained by comparing Eq.~\eqref{Eq:DSO_EC_KKT1} and Eq.~\eqref{eq:DUEEveDTCE}.
In the remainder of this section, we analytically derive the DUE solution assuming that the conjecture is true, and verify the conjecture in a similar manner to the morning commute problem.

First, we clarify the relationship between the bottleneck departure flows in the DSO and DUE states.
We obtain the following relationship by combining Eqs.~\eqref{Eq:DSO_EC_KKT2} and \eqref{Eq:EC_DSO_FC}:
\begin{align}
&\begin{cases}
x_i(\sigma_{i}(t))=\mu_i		&\text{if}\quad p_{i}(t)>0,\\
x_i(\sigma_{i}(t))\leq\mu_i	&\text{if}\quad p_{i}(t)=0,
\end{cases}
&\forall i\in\mathcal{N},\forall t\inT,\label{Eq:EC_DSOOutflow}
\end{align}
Eqs.~\eqref{eq:Definition_DUEEve_Y}, \eqref{Eq:EC_DUE_FC1}, and \eqref{eq:DUEEveQueue} are combined to yield the queuing delay condition in the DUE problem, as follows: 
\begin{align}
&\begin{cases}
x^{E}_{i}(\sigma_i^E(t))=\mu_i      &\text{if}\quad w^E_{i}(t)>0,\\
x^{E}_{i}(\sigma_i^E(t))\leq\mu_i   &\text{if}\quad w^E_{i}(t)=0,
\end{cases}
&&\forall i\in \mathcal{N},\forall t\inT.\label{eq:DUEEveQueue_2}
\end{align}
Substituting Eqs.~\eqref{Eq:TESolution_Tolls} and \eqref{eq:EveEquality_QD_PP} into \eqref{eq:DUEEveQueue_2}, the bottleneck departure flow is obtained as follows:
\begin{align}
&x_i^E(\sigma_i^E(t))=\mu_i,&\forall i\inN,t\inT_i.\label{eq:x=hatmu}
\end{align}
It should be noted that this equation only provides the profile of the bottleneck departure flows for a certain time period of $t\inT_i$; $x_i^E(t)$ for $t\notin\ClT_i$ can be obtained by deriving $q_i^E(t)$ for $t\inT_i$.
Similar to Eqs.~\eqref{eq:DUEqUsingx} and \eqref{eq:DUEqNUsingx} in the morning commute, we obtain the following:
\begin{align}
&q^E_{i}(t)= x^{E}_{i}(\sigma^E_{i}(t))\cdot \dot{\sigma}^E_{i}(t) - x^{E}_{i+1}(\sigma^E_{i+1}(t))\cdot \dot{\sigma}^E_{i+1}(t),	&&\forall i\in\mathcal{N}\setminus\{N\},\forall t\in\mathcal{T},\label{eq:DUEEveqUsingx}\\
&q^E_{N}(t) = x^{E}_{N}(\sigma^E_{N}(t))\cdot \dot{\sigma}^E_{N}(t),&& \forall t\in\mathcal{T}.\label{eq:DUEEveqNUsingx}
\end{align}
where $\dot\sigma_i^E(t)$ is derived by combining \eqref{Eq:TESolution_Tolls} for the DSO solution, \prettyref{eq:sigmaEve} from the definition of $\sigma_i^E(t)$, and \prettyref{eq:EveEquality_QD_PP} in \textbf{\prettyref{conj:EveEqualityQD_PP}},% separate sentence % refer to lemma and conjecture
\begin{align}
&\dot\sigma_i^E(t) = 1 + \sum_{j = 1}^{i}\dot{w}_{j}^{E}(t) = 1 +  \sum_{j = 1}^{i}\dot{p}_{j}^{E}(t) 
=\begin{cases}
1-\dot s(t)&\mathrm{if}\quad t\inT_i,\\
1&\mathrm{otherwise},
\end{cases}&\forall i\inN.\label{eq:EveDotsigma}
\end{align}
Substituting this equation and \prettyref{eq:x=hatmu} into \prettyref{eq:DUEEveqUsingx} and \prettyref{eq:DUEEveqNUsingx}, we obtain the following:
\begin{align}
&q_i^E(t)= (1-\dot s(t))\cdot\hat\mu_i,			&\forall i\in \mathcal{N}\setminus\{N\},\forall t \in \mathcal{T}_i,\label{eq:DUEEveqi}\\
&q_N^E(t)= (1-\dot s(t))\cdot\mu_N,			&\forall t \in \mathcal{T}_N,\label{eq:DUEEveqi_2}
\end{align}
Integrating this equation with respect to $t$, we obtain the following from the flow conservation condition \eqref{eq:DUEEveFCMass}:
\begin{align}
&q_i^E(t) = 0,&\forall i\inN,t\notin\ClT_i.\label{eq:DUEEveqN}
\end{align}
%Together, \prettyref{eq:DUEEveqi} and \prettyref{eq:DUEEveqN} provide the origin departure-flow rates $\{q_i^E(t)\}$. 
Substituting $\{q^E_i(t)\}$ into Eqs.~\eqref{eq:DUEEveqUsingx} and \eqref{eq:DUEEveqNUsingx}, the bottleneck departure flow $x_i^E(\sigma_i^E(t))$ can be obtained as follows:
\begin{align}
&x_i^E(\sigma_i^E(t))= 
\begin{cases}
\mu_i &\forall t\inT_i,\\
(1-\dot s(t))\cdot \mu_{i+1}&\forall t\inT_{i+1}\setminus\ClT_i,\\
0&\text{otherwise},
\end{cases}&\forall i\inN\setminus\{N\},\label{eq:DUEEvexi}\\
&x_N^E(\sigma_i^E(t))= 
\begin{cases}
\mu_N &\forall t\inT_N,\\
0&\text{otherwise}.\label{eq:DUEEvexN}
\end{cases}
\end{align}

From Eqs.~\eqref{eq:DUEEveqi}, \eqref{eq:DUEEveqi_2}, and \eqref{eq:DUEEveqN} of the origin departure-flow rates, the departure time windows in pricing equilibrium and DUE can be observed as the same, considering $\mathcal{T}_{i}^{E}$ as the departure time window of $i$-commuters in DUE yields $\mathcal{T}_{i}=\mathcal{T}^E_{i}$, $\forall i\in\mathcal{N}$, which is the same as that in the morning commute.
Consequently, all the variables in DUE can be analytically derived from the solution of the DSO assignment, assuming that the conjecture is true.

The derived flow variables satisfy all the equilibrium conditions in the DUE.
However, they may not satisfy certain physical constraints.
Specifically, the non-negativity constraints on the origin departure-flow rates (i.e. $q^E_i(t)\geq0$) and the \textit{capacity constraints} on the bottleneck departure flow rates (i.e. $x_i^E(\sigma_i^E(t))\leq\mu_i$).
Although the non-negativity constraints are the same as those in the morning commute problem, the capacity constraints are explicitly required in the evening commute problem, different from the morning commute problem.
This is because, in the morning commute problem, \prettyref{eq:Equality_OutFlow} in \textbf{\prettyref{lemma:Equality_OutFlow}} yields the bottleneck departure flows of the DUE for the entire time period from those of the DSO solution.
Thus, the bottleneck departure flows satisfies the capacity constraints.
For the evening commute, \prettyref{eq:x=hatmu} yields the bottleneck departure flows within the departure time window $\ClT_i$, and the capacity constraints are not guaranteed for $t\notin\ClT_i$. 
Thus, the conditions in the following lemma should be satisfied to ensure that the obtained flow variables are physically appropriate for the DUE in the evening commute.

\begin{lemma}\label{lemma:DUE_Feasible_Eve}
A DUE solution for the \textit{evening commute} under \textbf{\prettyref{conj:EveEqualityQD_PP}} exists if the following conditions hold true:
\begin{subequations}\label{Eq:DUE_EC_Feaibility}
\begin{align}
&\dot s(t) \leq 1,&& t\inT_N,\label{eq:DUEEve_Feasible_1}\\
&\dot s(t) \geq 1-\dfrac{\mu_{i}}{\mu_{i+1}}, && \forall i\inN\setminus\{N\},\forall t\inT_{i+1}\setminus\mathcal{T}_{i},\label{eq:DUEEve_Feasible_2}
\end{align}
\end{subequations}
where $\ClT_{N+1}\equiv\ClT_N$.
\end{lemma}
It should be noted that although Eqs.~\eqref{eq:DUEEve_Feasible_1} and \eqref{eq:DUEEve_Feasible_2} are similar to Conditions \eqref{eq:DUE_Feasible_1} and \eqref{eq:DUE_Feasible_2} for the morning commute problems, their meanings are slightly different.
Specifically, both conditions~\eqref{eq:DUEEve_Feasible_1} and \eqref{eq:DUEEve_Feasible_2} are interpreted as the conditions for the existence of the DUE solution, where the queueing delays are equal to the optimal prices in the pricing equilibrium.
In other words, even if the conditions are not satisfied, a DUE solution exists although its variables are different from those of the pricing equilibrium (as demonstrated in Section 6.4).

\begin{figure}[t]
\center
 \includegraphics[width=0.9\columnwidth,clip]{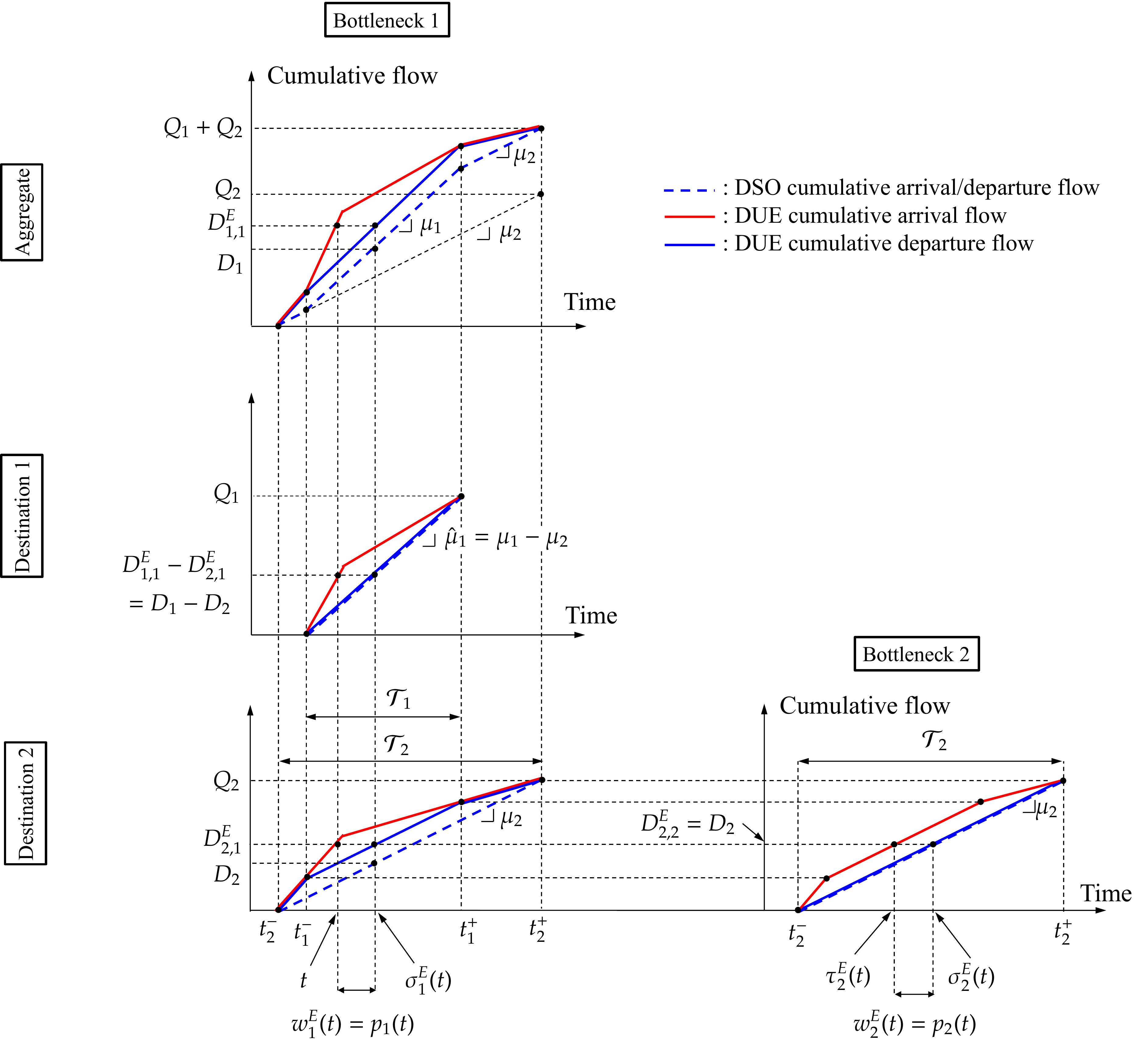}
 \caption{DUE cumulative flow curves of bottleneck $i$ for evening commute in a corridor network. Note that $c_1=c_2=0$, $D_{i,m}^E\equiv D_i^E(\sigma^E_m(t))$,  $D_{i}\equiv D_i(t)$.}
\label{fig:DUEEve_Example_CFC}
\end{figure}

The following proposition summarizes the above discussions of the DUE solution for the evening commute: 
\begin{proposition}\label{prop:DUEAnaSolEve}
Consider the DUE problem for the \textit{evening commute} in a reduced corridor network. 
Suppose that Eqs.~\eqref{eq:DUEEve_Feasible_1} and \eqref{eq:DUEEve_Feasible_2} are satisfied. 
The DUE solution for $i$-commuters ($\forall i\inN$) is given as follows:
\begin{subequations}
\begin{align}
&\bullet \textrm{departure-flow rate}: &&
q^E_i(t)=\left\{
\begin{array}{ll}
(1-\dot s(t))\cdot\hat\mu_i&\mathrm{if}\quad t\inT_i\\
0&\mathrm{otherwise}
\end{array}
\right.\label{eq:DUEqEve}\\
& &&q^E_N(t)=\left\{
\begin{array}{ll}
(1-\dot s(t))\cdot \mu_N&\mathrm{if}\quad t\inT_N\\
0&\mathrm{otherwise}
\end{array}
\right.\label{eq:DUEqEve2}\\
&\bullet \textrm{commuting cost}: &&\Red{\rho^E_i=\bar s\left(\dfrac{Q_i}{\hat{\mu}_i}\right)+c_i}\label{eq:DUErhoEve}\\
&\bullet \textrm{queuing delay}: &&w_i^E(t)=
\left\{
\begin{array}{ll}
\rho^E_i-s(t)- c_{i} - \sum_{j = 1}^{i-1}w_j^E(t)&\mathrm{if}\quad t\inT_i\\
0&\mathrm{otherwise}
\end{array}
\right.\label{eq:DUEwEve}
\end{align}
\end{subequations}
\end{proposition}
\noindent Note that the properties discussed in the morning commute problem (e.g., the nested arrival times and a Pareto improvement property) are true in this evening commute problem.

\begin{figure}[t]
\centering
\includegraphics[width=80mm,clip]{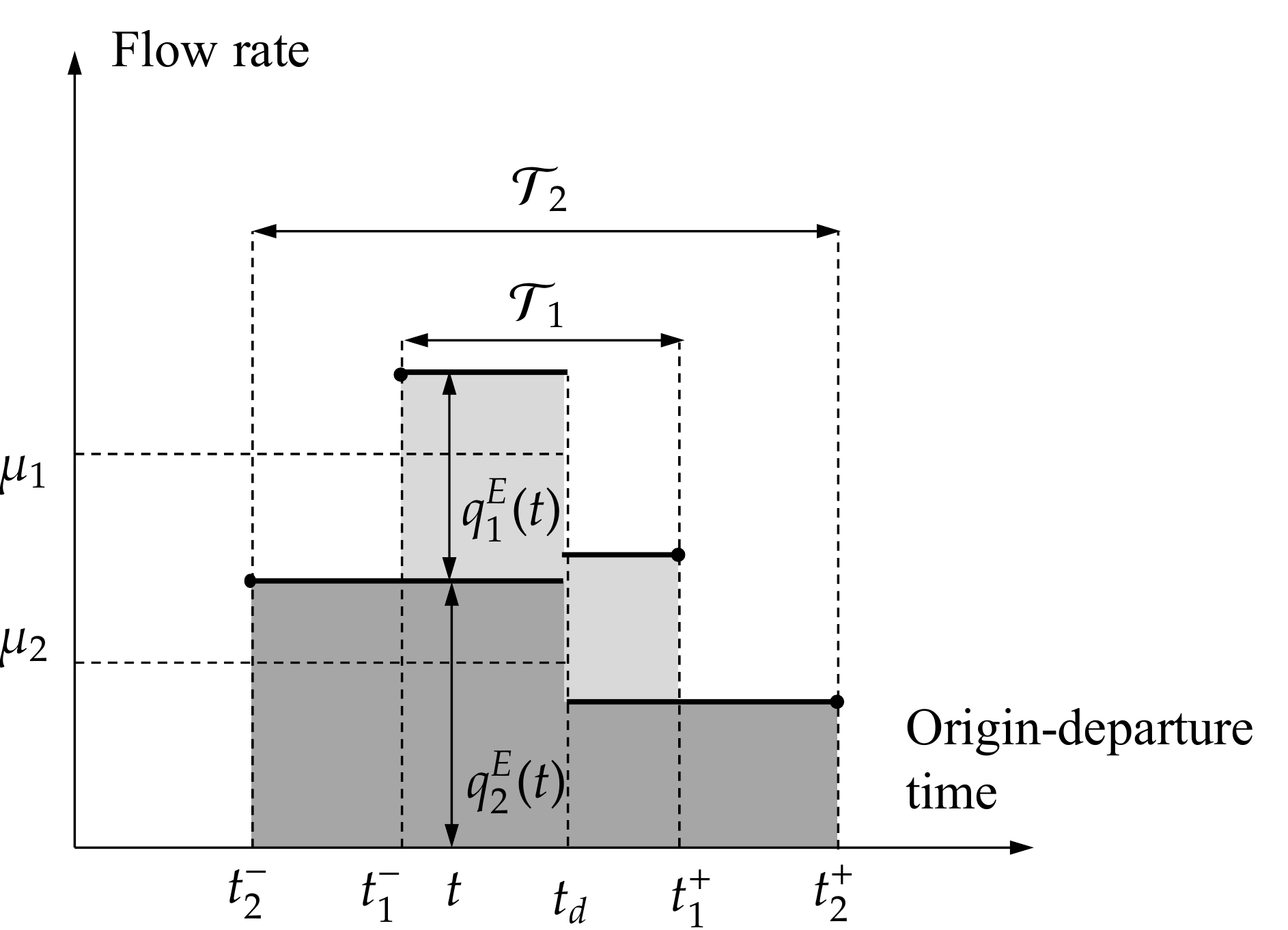}
\caption{Arrival-flow rate at destinations of the DUE solution in the evening commute}
\label{fig:DUEEveAnaSol}
\end{figure}

Figure~\ref{fig:DUEEve_Example_CFC} illustrates the cumulative flow curves of bottlenecks in a simple example with $I=2$ and a piecewise linear schedule delay function (for simplicity, we again set $c_1=c_2=0$). 
This figure reveals that the aggregate departure flows of bottlenecks in the DSO and DUE problems are not the same in the evening commute (i.e. $D_{1,1}^E\neq D_1$). 
However, the destination arrival flows are the same for both the problems, namely, $D_{1,1}^E-D_{2,1}^E= D_1-D_2$ and $D_{2,2}^E=D_2$, which is opposite to that in the morning commute. 
The origin departure-flow rates $\{q^E_i(t)\}$ are shown in Figure~\ref{fig:DUEEveAnaSol}.

%%%%%%%%%%%%%%%%%%%%%%%%%%%%%%%%%%%%%%%%%%%%%%%%%%%%%%%%%%%%%%%%%%%%
% Section 6

\section{Numerical examples}\label{sec:Num}
This section provides numerical examples to illustrate the theoretical results for the relationship between the DSO and DUE clarified so far.
We present several examples for which Condition~\eqref{Eq:DUE_Feaibility} for the morning commute and Condition~\eqref{Eq:DUE_EC_Feaibility} for the evening commute are satisfied.
In these examples, we verify our analyses and analytical solutions in the morning and evening commute problems.
Moreover, we present other examples for which Conditions Eqs.~\eqref{eq:DUE_Feasible_2} and \eqref{eq:DUEEve_Feasible_2} (i.e., conditions ensuring the existence of the DUE where the queueing delays are equal to the optimal prices in the pricing equilibrium) are not satisfied.
In these examples, we demonstrate the existence of DUE where the aggregate flow and queuing delay patterns do not correspond to the patterns of aggregate flow and the Lagrange multiplier of the DSO state.

\subsection{Experimental settings}
We consider a corridor network with three bottlenecks.
The parameters are set such that all bottlenecks are non-false.
The travel demands are $\mathbf{Q}=[100,350,250]^{T}$, the bottleneck capacities are $\boldsymbol{\mu}=[50,30,10]^{T}$, and the free-flow travel times are $\mathbf{c} = \mathbf{0}$ for simplicity.
Note that the time period $\mathcal{T}$ is divided into a finite number of intervals, $K$, labeled as $k\in\mathcal{K}\equiv\{1, 2, \ldots,K\}$, which has a finite duration $\Delta k$.

We employ two piecewise linear schedule delay functions for the morning and evening commute problems, respectively (four schedule functions in total).
%In each commute problem, one schedule delay function satisfies Conditions~\eqref{eq:DUE_Feasible_2} or \eqref{eq:DUEEve_Feasible_2}, and the other does not satisfy the condition.
Specifically, the schedule functions are set as follows: 
\begin{align*}
    &\text{Example 1 (for the morning commute)}: s(t)=\max\{0.5(30-t),0.5(t-30)\}\\
    &\text{Example 2 (for the morning commute)}: s(t)=\max\{0.5(30-t),8(t-30)\}\\
    &\text{Example 3 (for the evening commute)}: s(t)=\max\{0.5(30-t),0.5(t-30)\}\\
    &\text{Example 4 (for the evening commute)}: s(t)=\max\{8(30-t),0.5(t-30)\}
\end{align*}
The first schedule delay function satisfies Condition~\eqref{eq:DUE_Feasible_2}; the second function does not satisfy the condition.
The third function satisfies Condition~\eqref{eq:DUEEve_Feasible_2}, and the fourth function does not satisfy this condition.

To compute the DSO problems in discrete time, which are formulated as finite-dimensional LP problems, we apply the interior point algorithm.
Meanwhile, the DUE problems are formulated as finite-dimensional LCPs\footnote{For the details of the formulation of the finite-dimensional LCPs (including the problems for the evening commute), refer to~\citep[][]{akamatsu2015corridor}.}.
Specifically, the DUE problem for the morning commute is expressed as
\begin{align}
\mathbf{F}(\mathbf{X})\cdot\mathbf{X}=0,\mathbf{F}(\mathbf{X})\geq0, \mathbf{X}\geq  \mathbf{0},
\end{align}
where
\begin{align}
&\mathbf{X} = [\mathbf{q}^E, \mathbf{w}^E, \boldsymbol{\rho}^E]^{T}, \mathbf{F}( \mathbf{X} )=\mathbf{M}\mathbf{X}+\mathbf{b},\\
&\mathbf{M}=\left[
\begin{array}{ccc}
\mathbf{0}												&	\mathbf{I}_{K}	\otimes \mathbf{L}																							&	-\mathbf{1}_{K} \otimes	\mathbf{I}_{N}\\ 
-\mathbf{I}_{K}\otimes	\mathbf{L}^{T}		&	\boldsymbol{\Delta}_{K} \otimes \text{diag}(\boldsymbol{\mu})(\mathbf{I}_{N} - \mathbf{L})		&\mathbf{0}\\
\mathbf{1}_{K}^{T}\otimes \mathbf{I}_{N} &\mathbf{0}				&\mathbf{0}
\end{array}
\right],\quad  
\mathbf{b}=
\left[
\begin{array}{c}
\mathbf{s}\\
\mathbf{1}_{K}\otimes	\boldsymbol{\mu}\\
-\mathbf{Q}.
\end{array}\right].
\end{align}
Here, $\otimes$ is the operator of Kronecker product; $\mathbf{L}$ is a lower triangular $N\times N$ matrix where all non-zero entries equal to 1; 
$\mathbf{I}_k$ $(\forall k\in\mathbb{N})$ is an identity matrix of dimensions $k\times k$; 
$\mathbf{1}_k$ is a vector of dimensions $k\times 1$ with all entries equal to 1; 
$\mathbf{s}\equiv [s(t)\mathbf{1}]_{t = 1}^{K}$.
$\boldsymbol{\Delta}_{K}$ is the $K\times K$ matrix representing the first order backward difference operator:
\begin{align}
\boldsymbol{\Delta}_{K} = \left[
\begin{array}{cccc}
1 		& 				& 				& \\
-1 	& 	1			&				& \\
		& \ddots	& \ddots	& \\
		&				&		-1		& 1
\end{array} \right].
\end{align}

%All constant travel costs are set to be zero and thus is omitted for simplicity.
The LCP can be equivalently represented as the following quadratic programming (QP)\footnote{For the proof that the LCP is equivalent to the QP, refer to \cite{cottle1992linear}.}:
\begin{align}
\min \mathbf{F}(\mathbf{X})\cdot \mathbf{X}, \quad  \text{s.t.}\quad \mathbf{X}\geq \mathbf{0}, \mathbf{F}( \mathbf{X})\geq \mathbf{0}.  
\end{align}
We can apply the well-known Frank-Wolfe algorithm to solve it.

\begin{figure*}[t]
	\centering
	\begin{minipage}[t]{1.0\textwidth}
	\centering
		\includegraphics[clip, width=1.0\columnwidth]{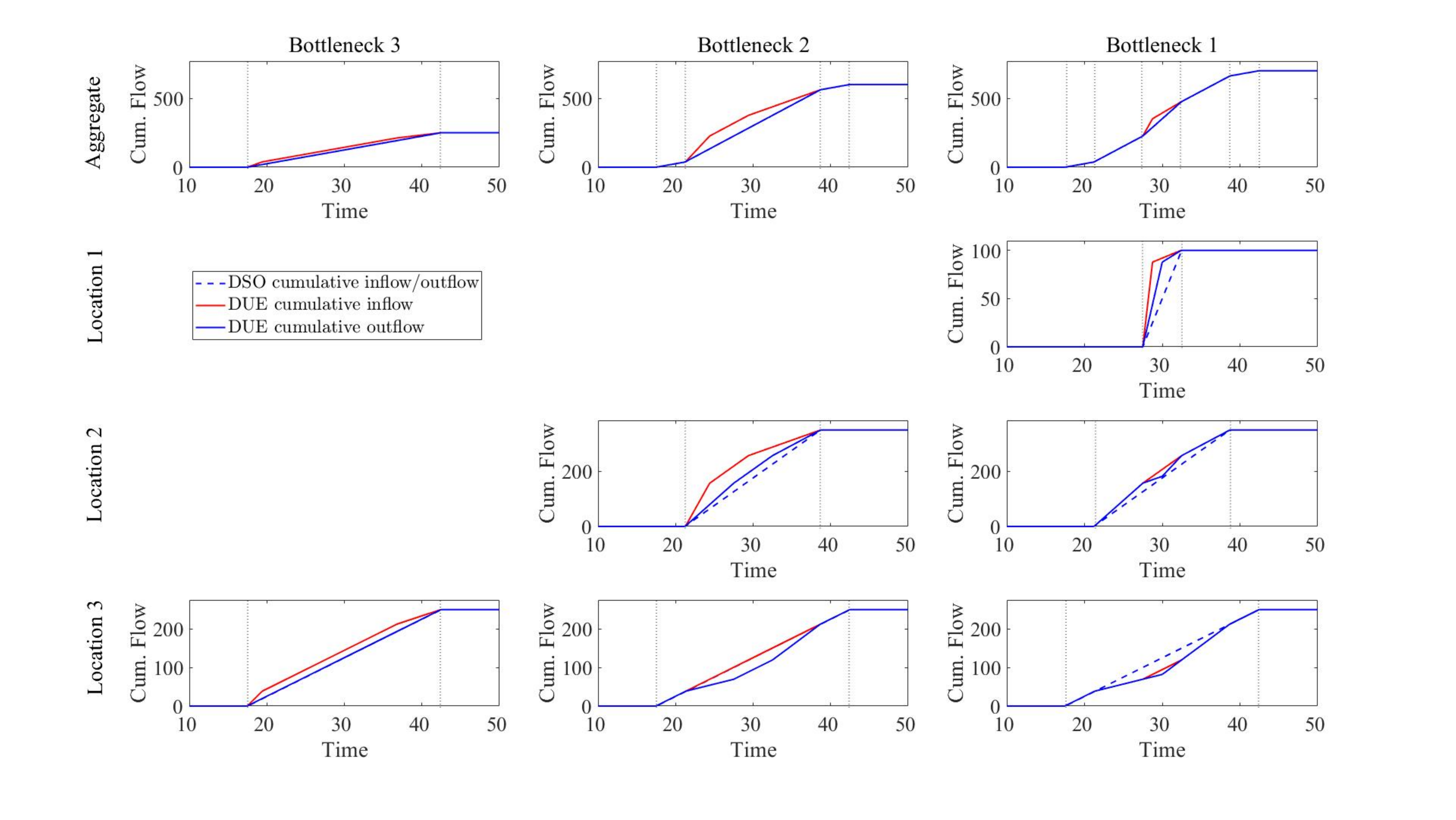}
		\vspace{-10mm}
		\subcaption{Cumulative flow curves of bottlenecks}
		\label{Fig:NE1_CumFlows}
	\end{minipage}\vspace{5mm}\\
	\begin{minipage}[t]{0.45\textwidth}
	\centering
		\includegraphics[clip, width=1.0\columnwidth]{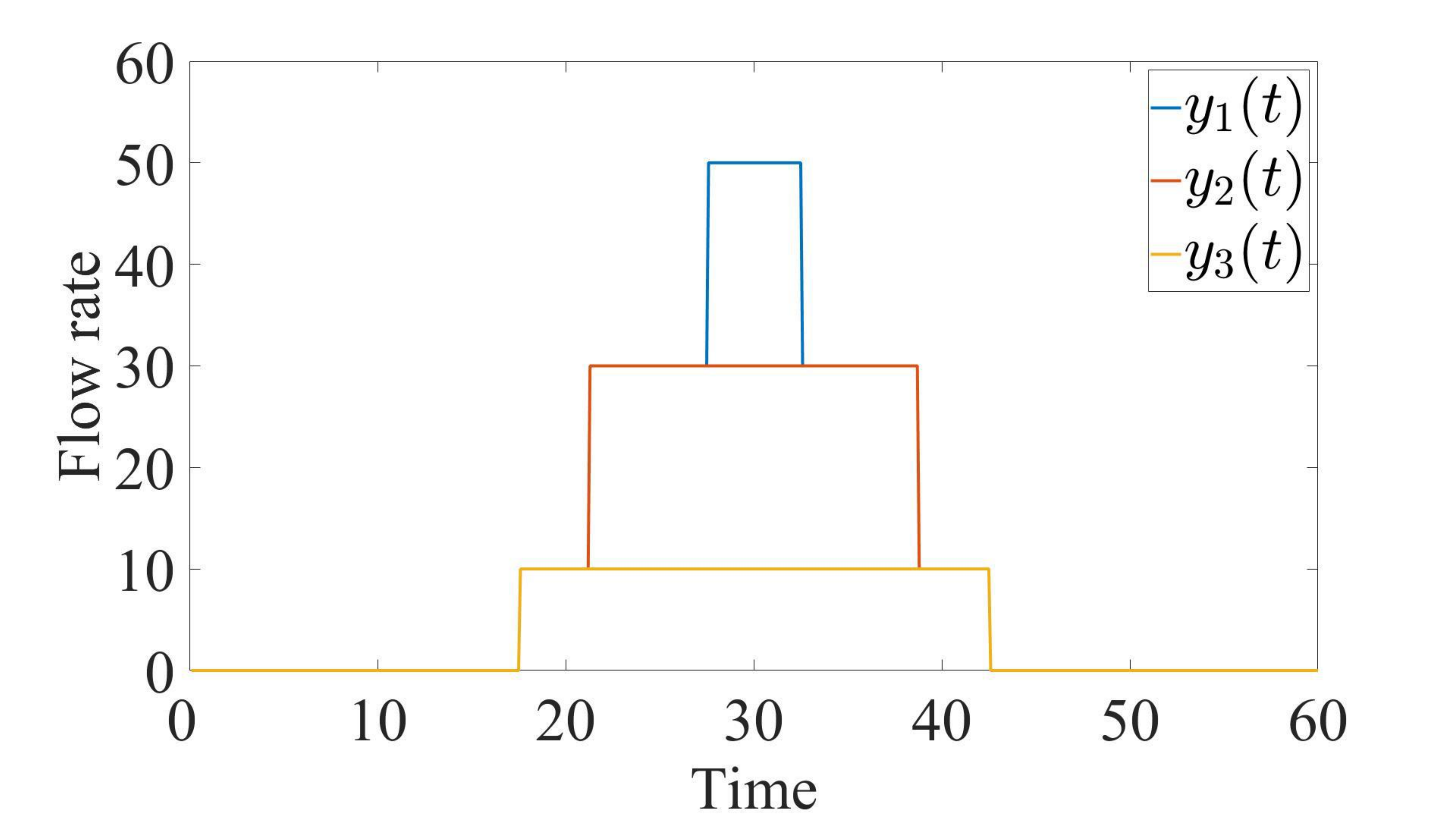}
		\subcaption{Arrival flow rate in the DSO state}
		\label{Fig:NE1_DSOFlows}
	\end{minipage}
	\begin{minipage}[t]{0.45\textwidth}
	\centering
		\includegraphics[clip, width=1.0\columnwidth]{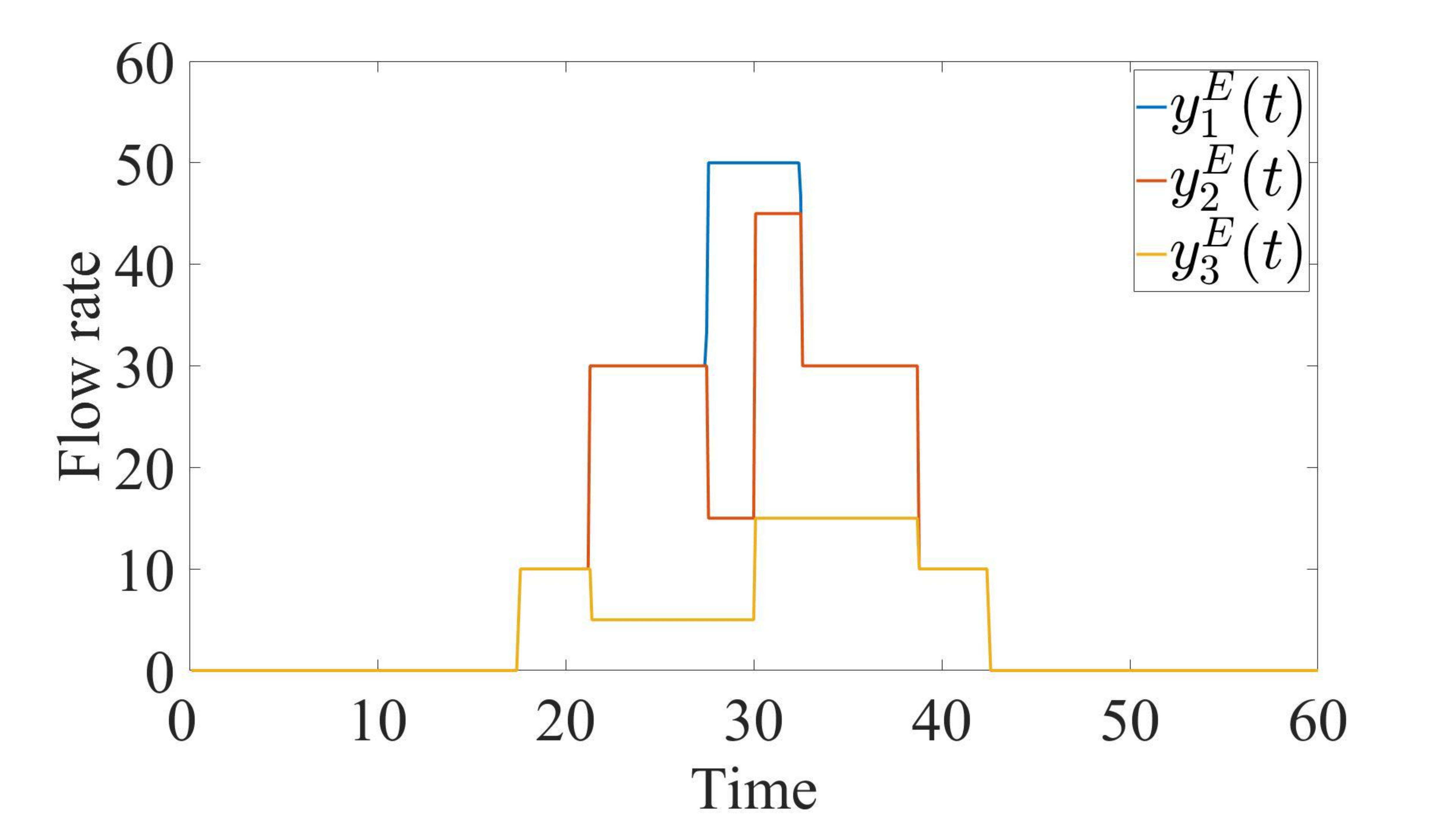}
		\subcaption{Arrival flow rate in the DUE state}
		\label{Fig:NE1_DUEFlows}
	\end{minipage}\vspace{5mm}\\
	\begin{minipage}[t]{0.45\textwidth}
	\centering
		\includegraphics[clip, width=1.0\columnwidth]{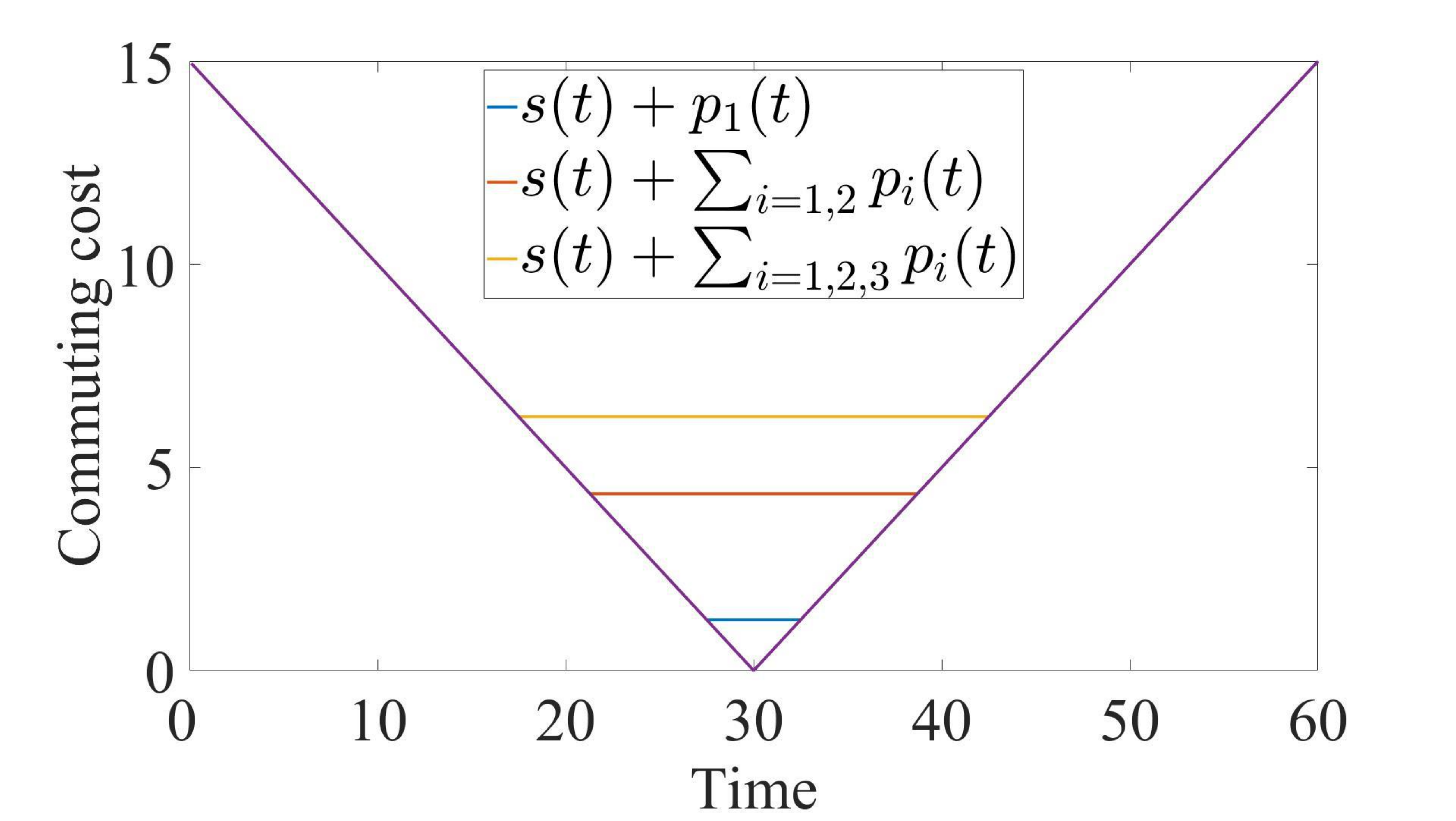}
		\subcaption{Equilibrium commuting costs in the DSO state}
		\label{Fig:NE1_DSOCosts}
	\end{minipage}
	\begin{minipage}[t]{0.45\textwidth}
	\centering
		\includegraphics[clip, width=1.0\columnwidth]{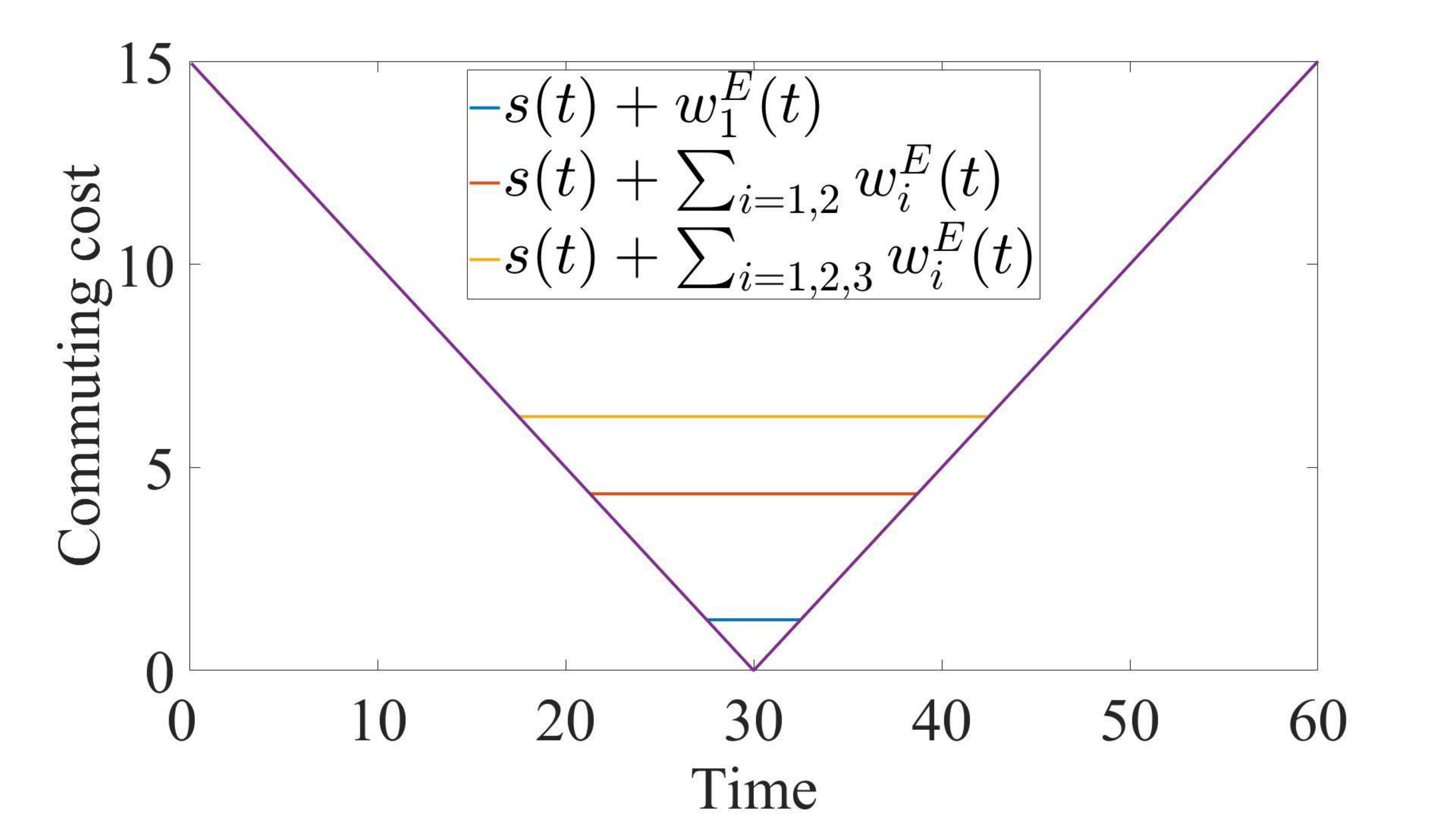}
		\subcaption{Equilibrium commuting costs in the DUE state}
		\label{Fig:NE1_DUECosts}
	\end{minipage}
	\vspace{0mm}
	\caption{Flow and cost patterns in DSO and DUE states in Example 1}
	\label{Fig:NE_Example1}
	\vspace{-2mm}
\end{figure*}

\begin{figure*}[t]
	\centering
	\begin{minipage}[t]{1.0\textwidth}
	\centering
		\includegraphics[clip, width=1.0\columnwidth]{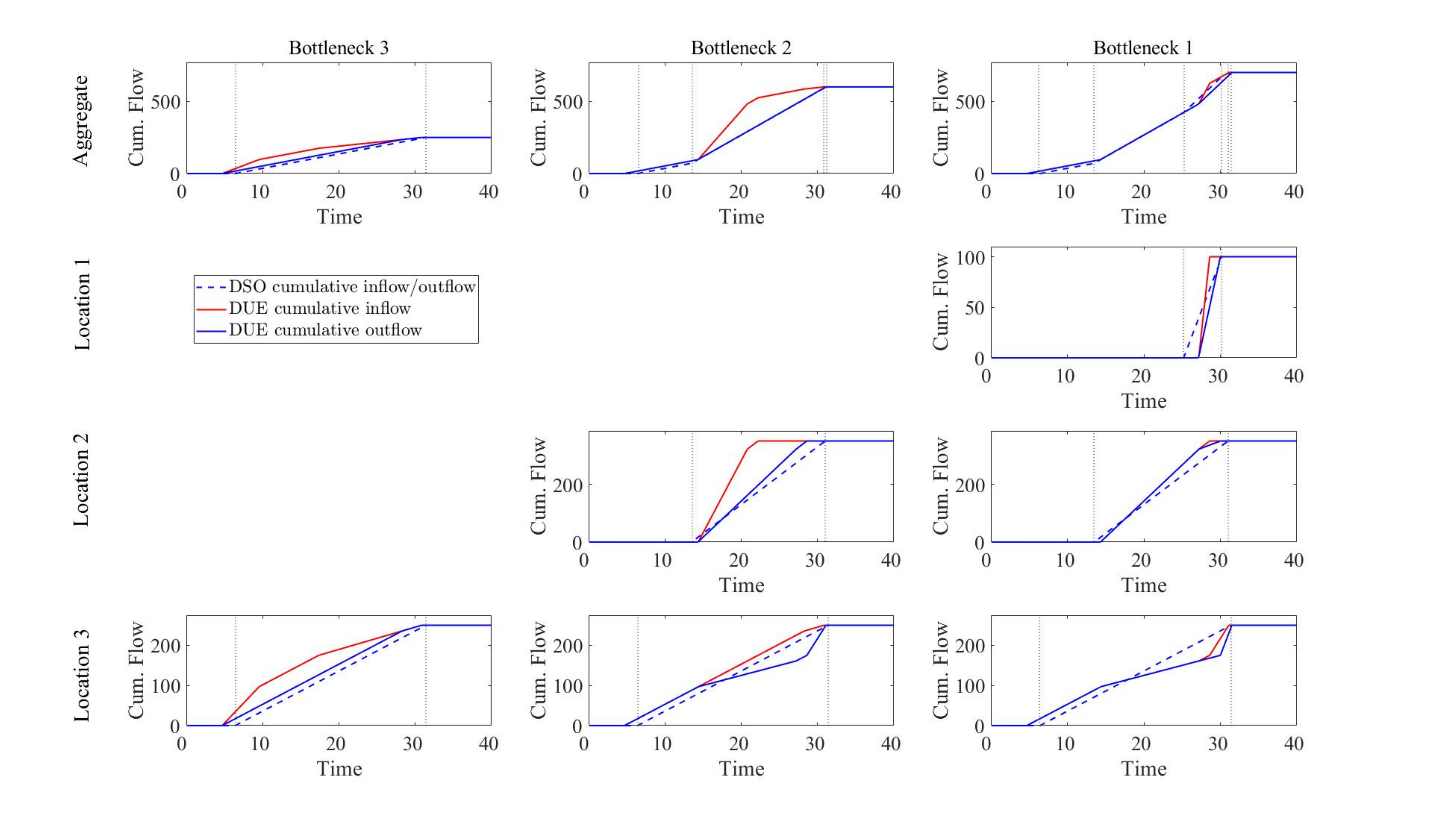}
		\vspace{-10mm}
		\subcaption{Cumulative flow curves of bottlenecks}
		\label{Fig:NE2_CumFlows}
	\end{minipage}\vspace{5mm}\\
	\begin{minipage}[t]{0.45\textwidth}
	\centering
		\includegraphics[clip, width=1.0\columnwidth]{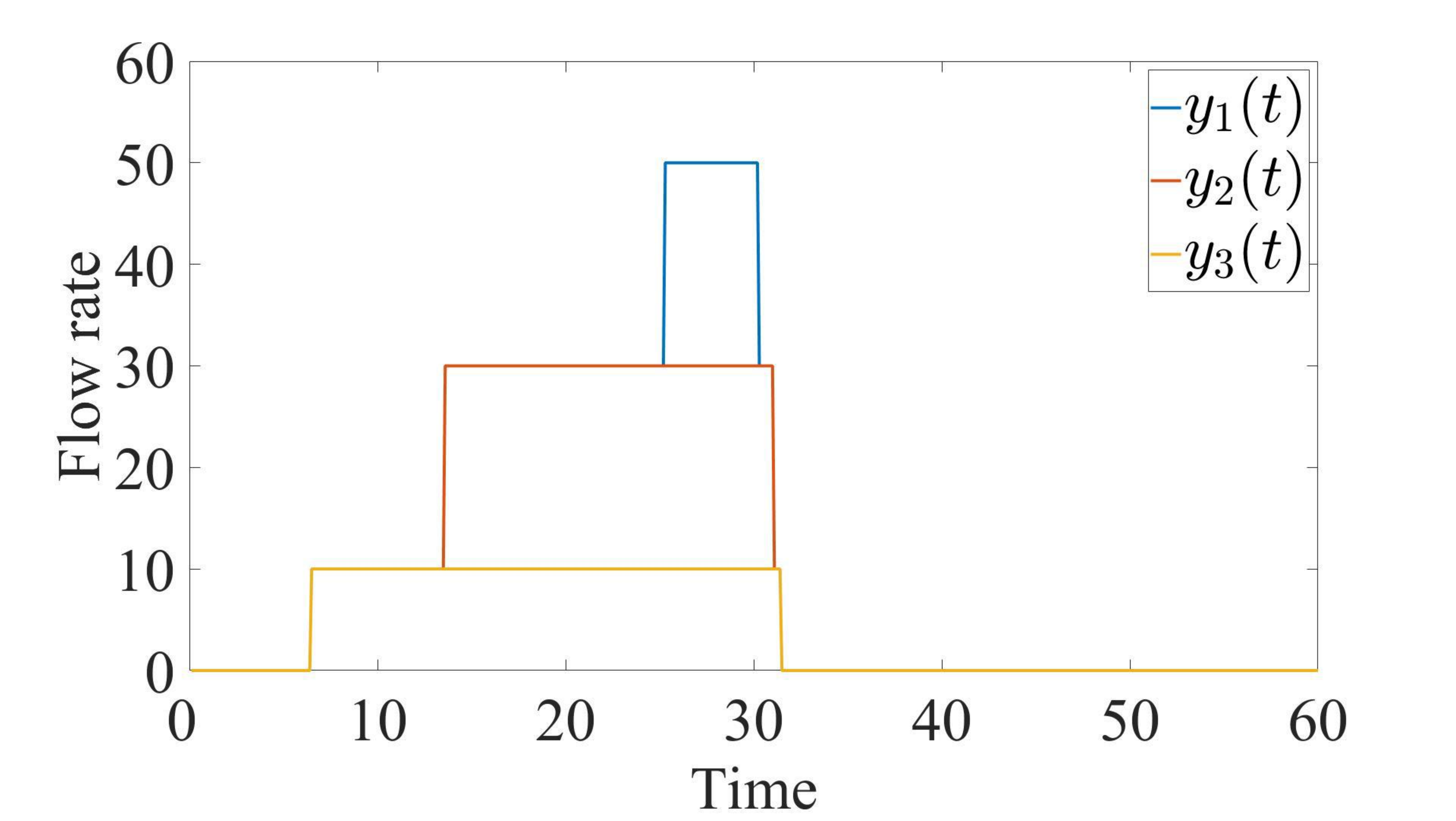}
		\subcaption{Arrival flow rate in the DSO state}
		\label{Fig:NE2_DSOFlows}
	\end{minipage}
	\begin{minipage}[t]{0.45\textwidth}
	\centering
		\includegraphics[clip, width=1.0\columnwidth]{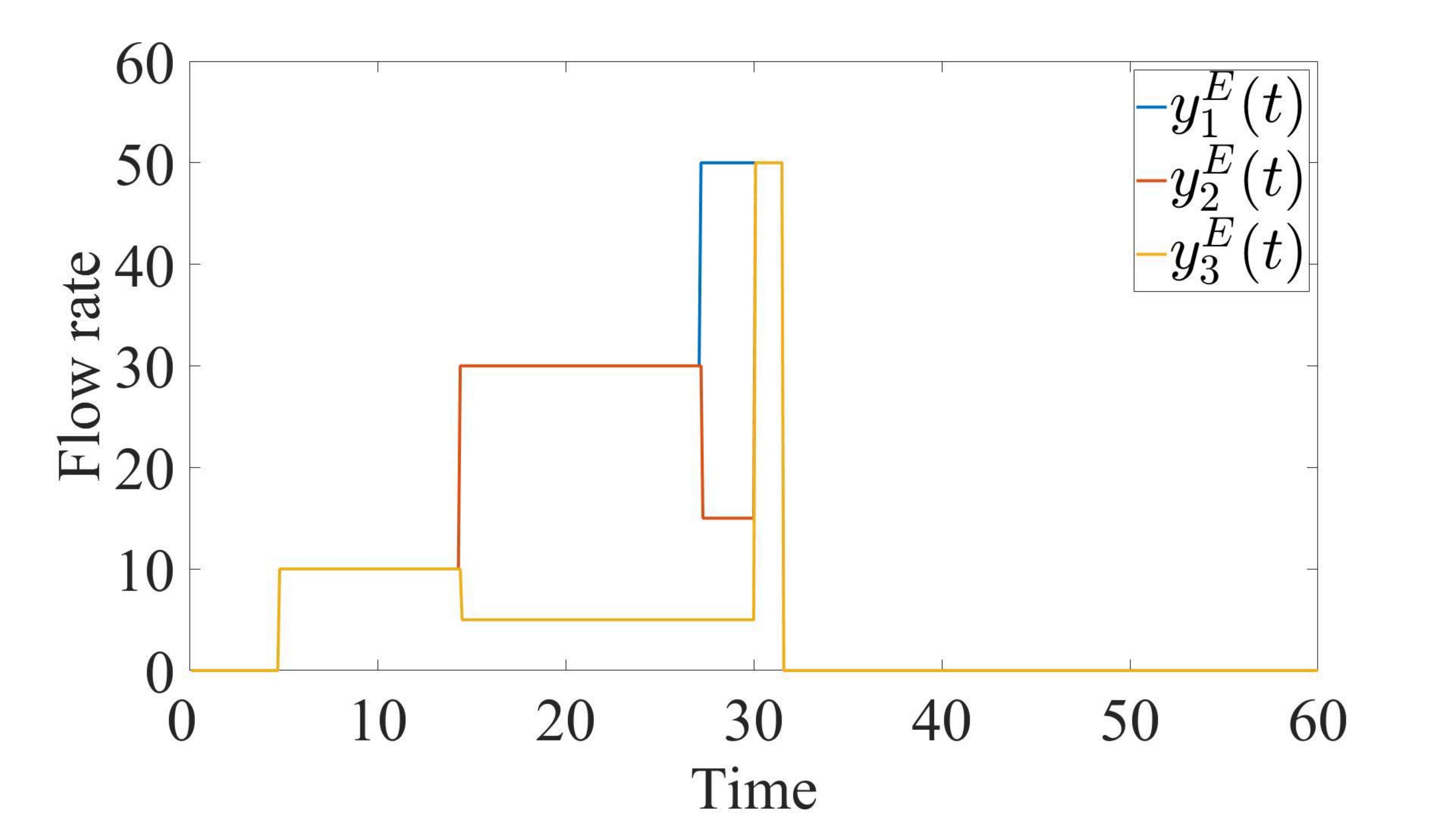}
		\subcaption{Arrival flow rate in the DUE state}
		\label{Fig:NE2_DUEFlows}
	\end{minipage}\vspace{5mm}\\
	\begin{minipage}[t]{0.45\textwidth}
	\centering
		\includegraphics[clip, width=1.0\columnwidth]{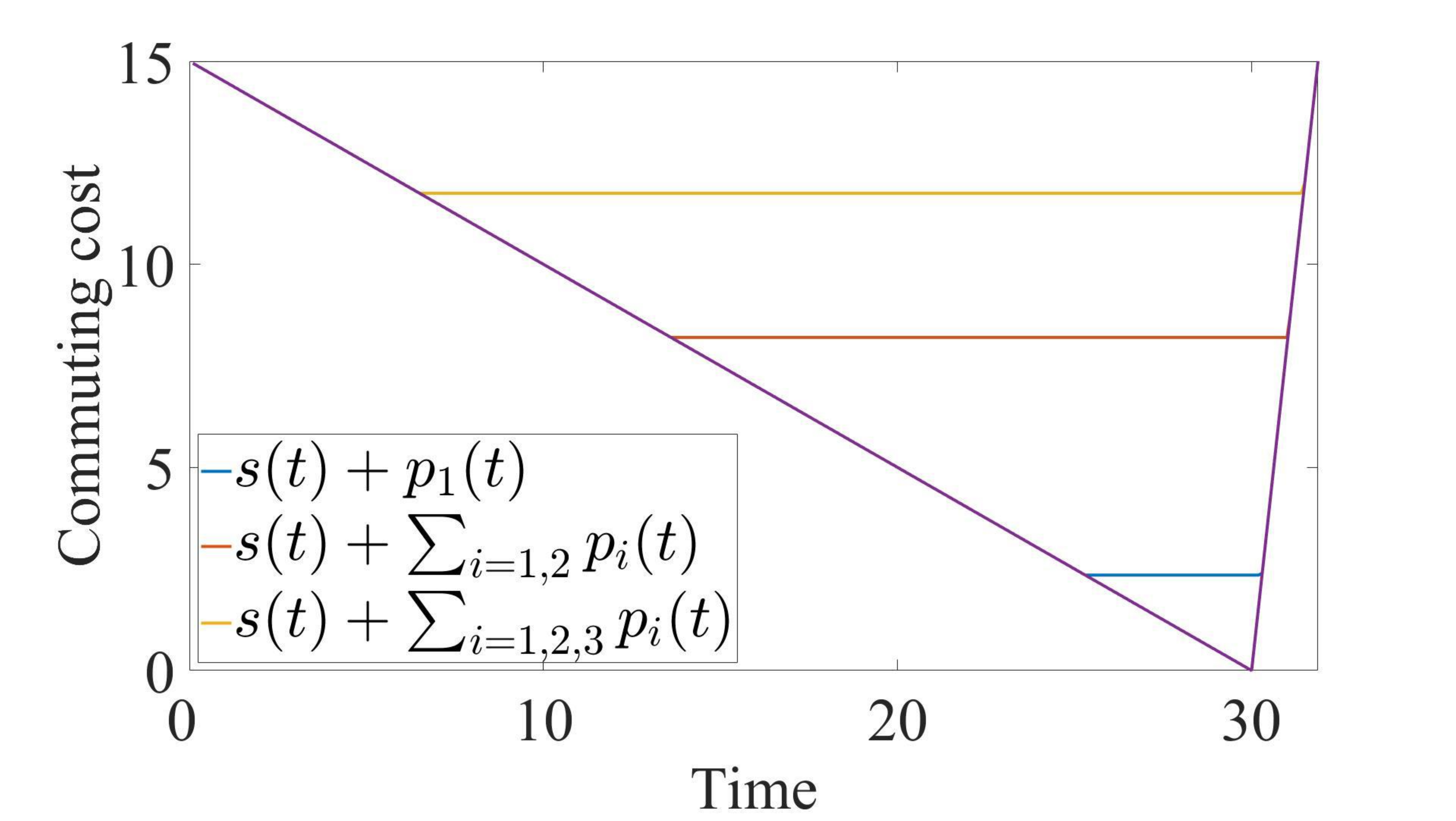}
		\subcaption{Equilibrium commuting costs in the DSO state}
		\label{Fig:NE2_DSOCosts}
	\end{minipage}
	\begin{minipage}[t]{0.45\textwidth}
	\centering
		\includegraphics[clip, width=1.0\columnwidth]{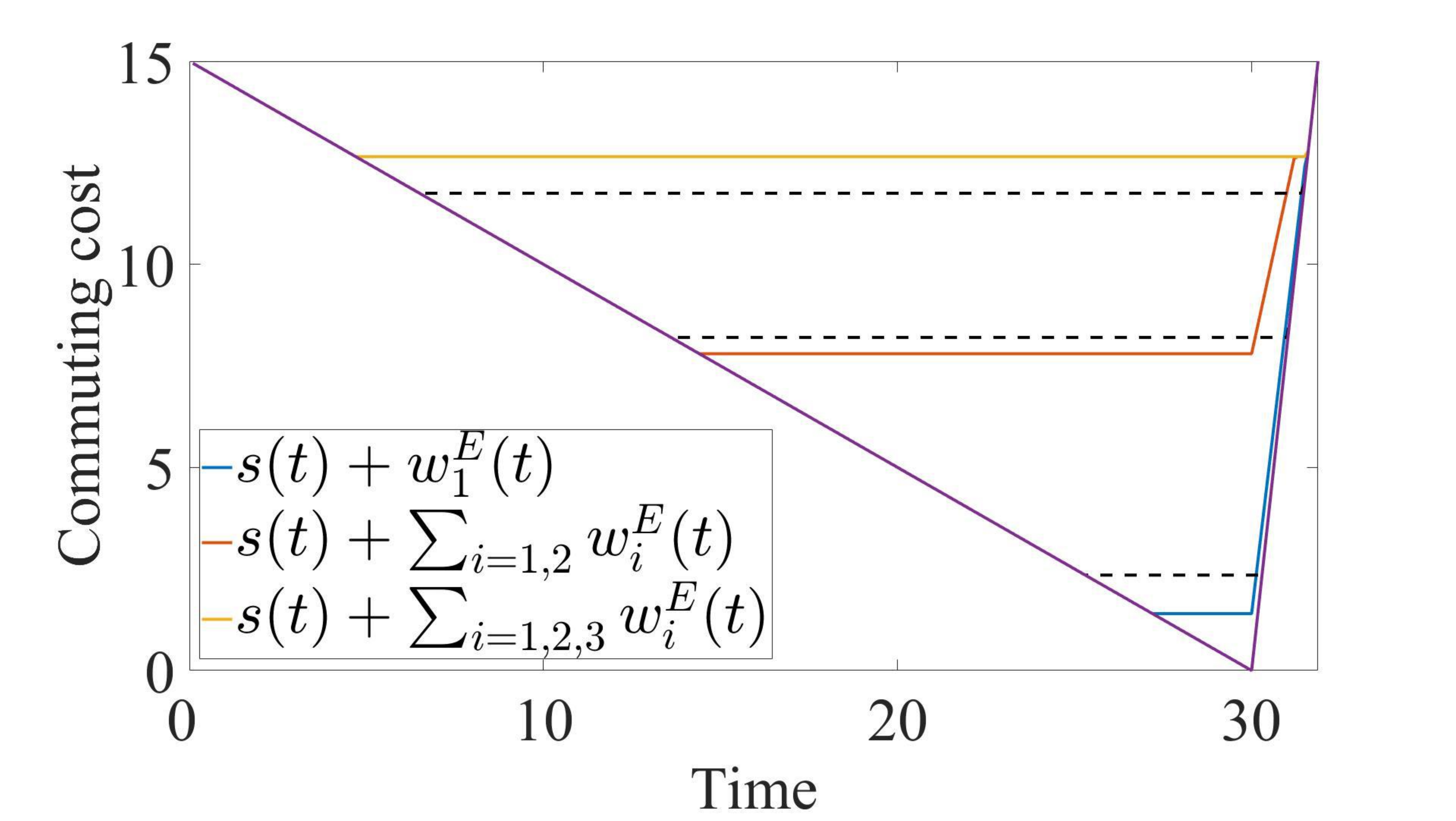}
		\subcaption{Equilibrium commuting costs in the DUE state}
		\label{Fig:NE2_DUECosts}
	\end{minipage}
	\vspace{0mm}
	\caption{Flow and cost patterns in DSO and DUE states in Example 2}
	\label{Fig:NE_Example2}
	\vspace{-2mm}
\end{figure*}

\begin{figure*}[t]
	\centering
	\begin{minipage}[t]{1.0\textwidth}
	\centering
		\includegraphics[clip, width=1.0\columnwidth]{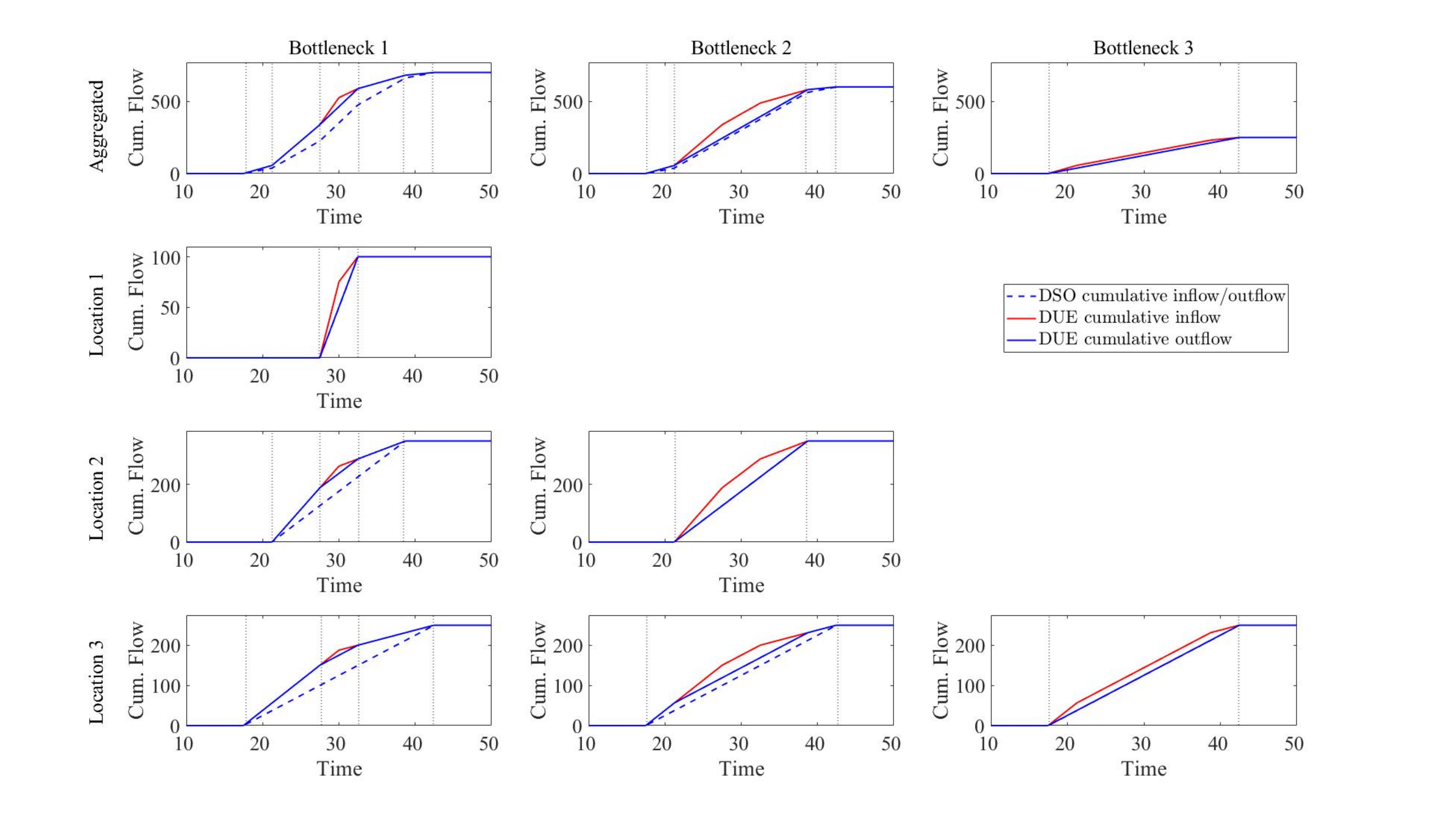}
		\vspace{-10mm}
		\subcaption{Cumulative flow curves of bottlenecks}
		\label{Fig:NE3_CumFlows}
	\end{minipage}\vspace{5mm}\\
	\begin{minipage}[t]{0.45\textwidth}
	\centering
		\includegraphics[clip, width=1.0\columnwidth]{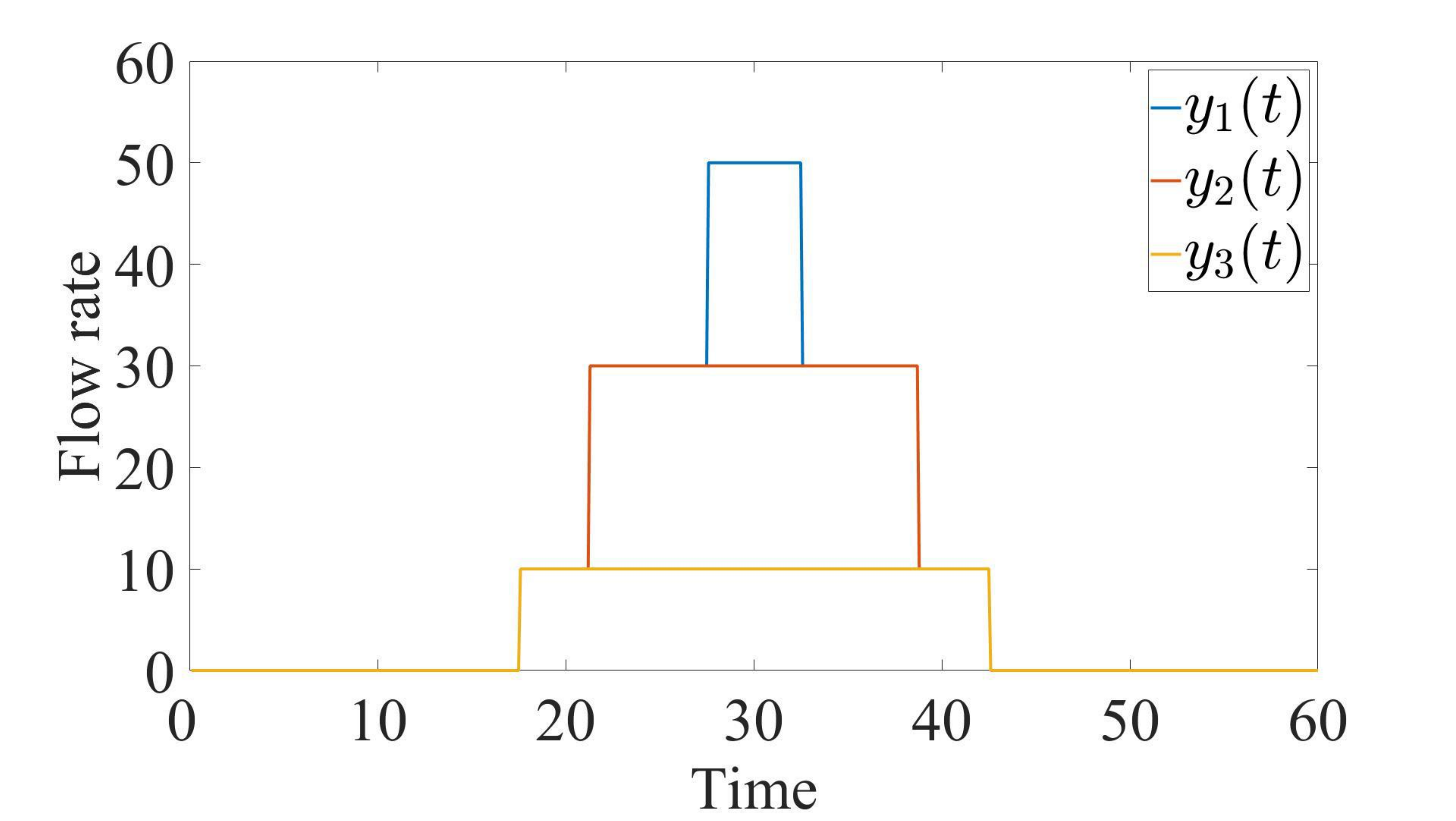}
		\subcaption{Arrival flow rate in the DSO state}
		\label{Fig:NE3_DSOFlows}
	\end{minipage}
	\begin{minipage}[t]{0.45\textwidth}
	\centering
		\includegraphics[clip, width=1.0\columnwidth]{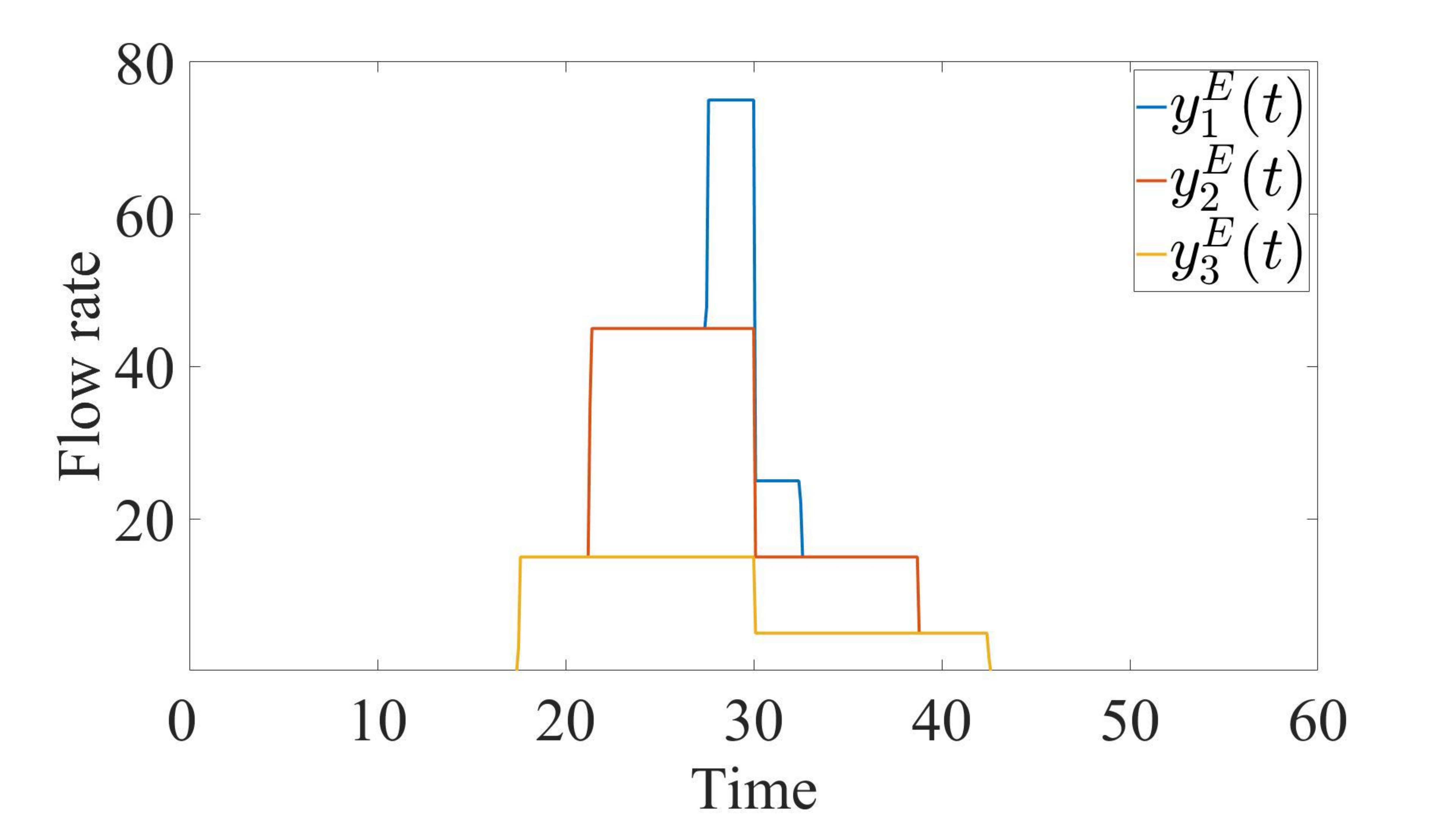}
		\subcaption{Arrival flow rate in the DUE state}
		\label{Fig:NE3_DUEFlows}
	\end{minipage}\vspace{5mm}\\
	\begin{minipage}[t]{0.45\textwidth}
	\centering
		\includegraphics[clip, width=1.0\columnwidth]{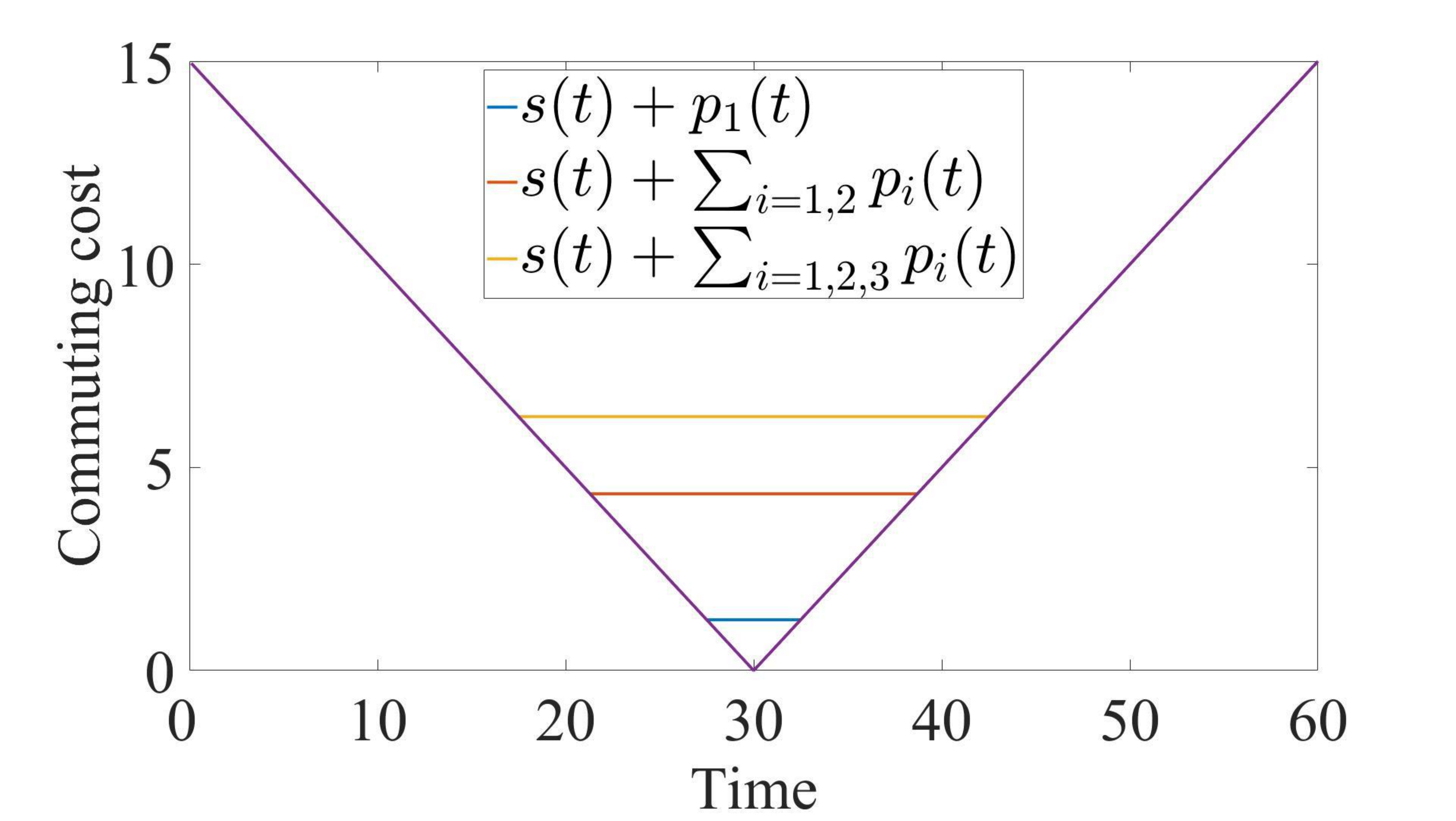}
		\subcaption{Equilibrium commuting costs in the DSO state}
		\label{Fig:NE3_DSOCosts}
	\end{minipage}
	\begin{minipage}[t]{0.45\textwidth}
	\centering
		\includegraphics[clip, width=1.0\columnwidth]{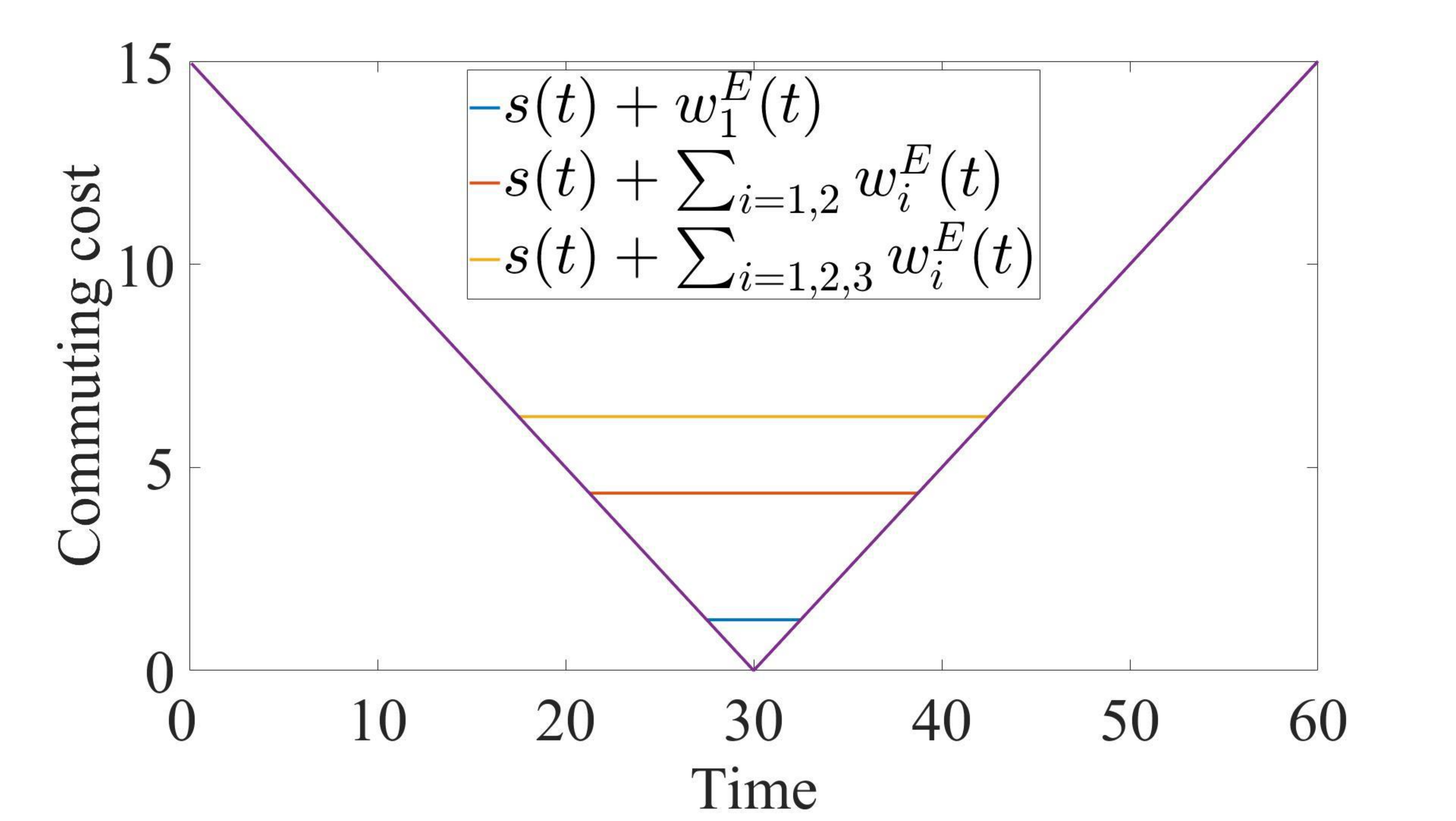}
		\subcaption{Equilibrium commuting costs in the DUE state}
		\label{Fig:NE3_DUECosts}
	\end{minipage}
	\vspace{0mm}
	\caption{Flow and cost patterns in DSO and DUE states in Example 3}
	\label{Fig:NE_Example3}
	\vspace{-2mm}
\end{figure*}

\begin{figure*}[t]
	\centering
	\begin{minipage}[t]{1.0\textwidth}
	\centering
		\includegraphics[clip, width=1.0\columnwidth]{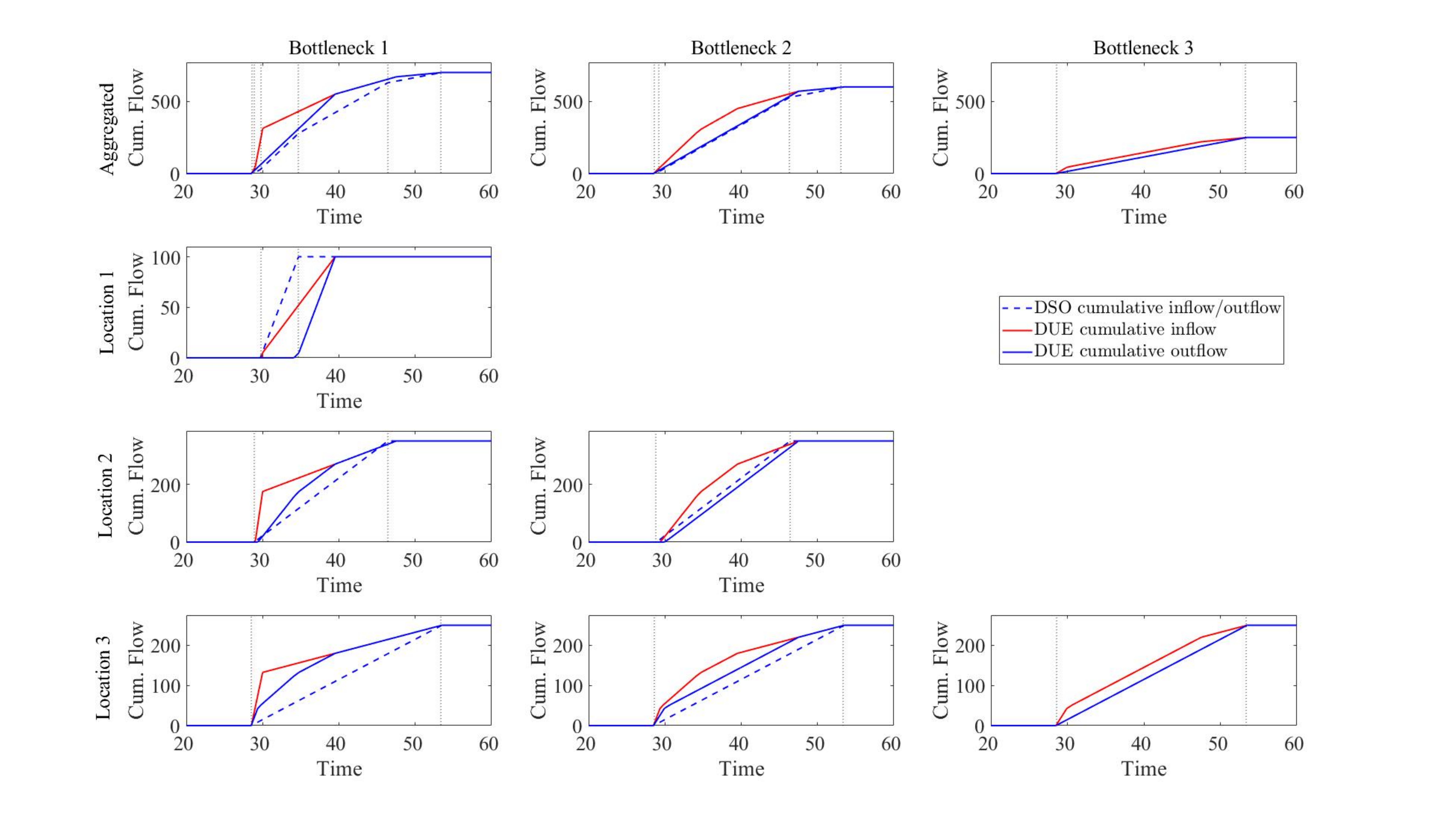}
		\vspace{-10mm}
		\subcaption{Cumulative flow curves of bottlenecks}
		\label{Fig:NE4_CumFlows}
	\end{minipage}\vspace{5mm}\\
	\begin{minipage}[t]{0.45\textwidth}
	\centering
		\includegraphics[clip, width=1.0\columnwidth]{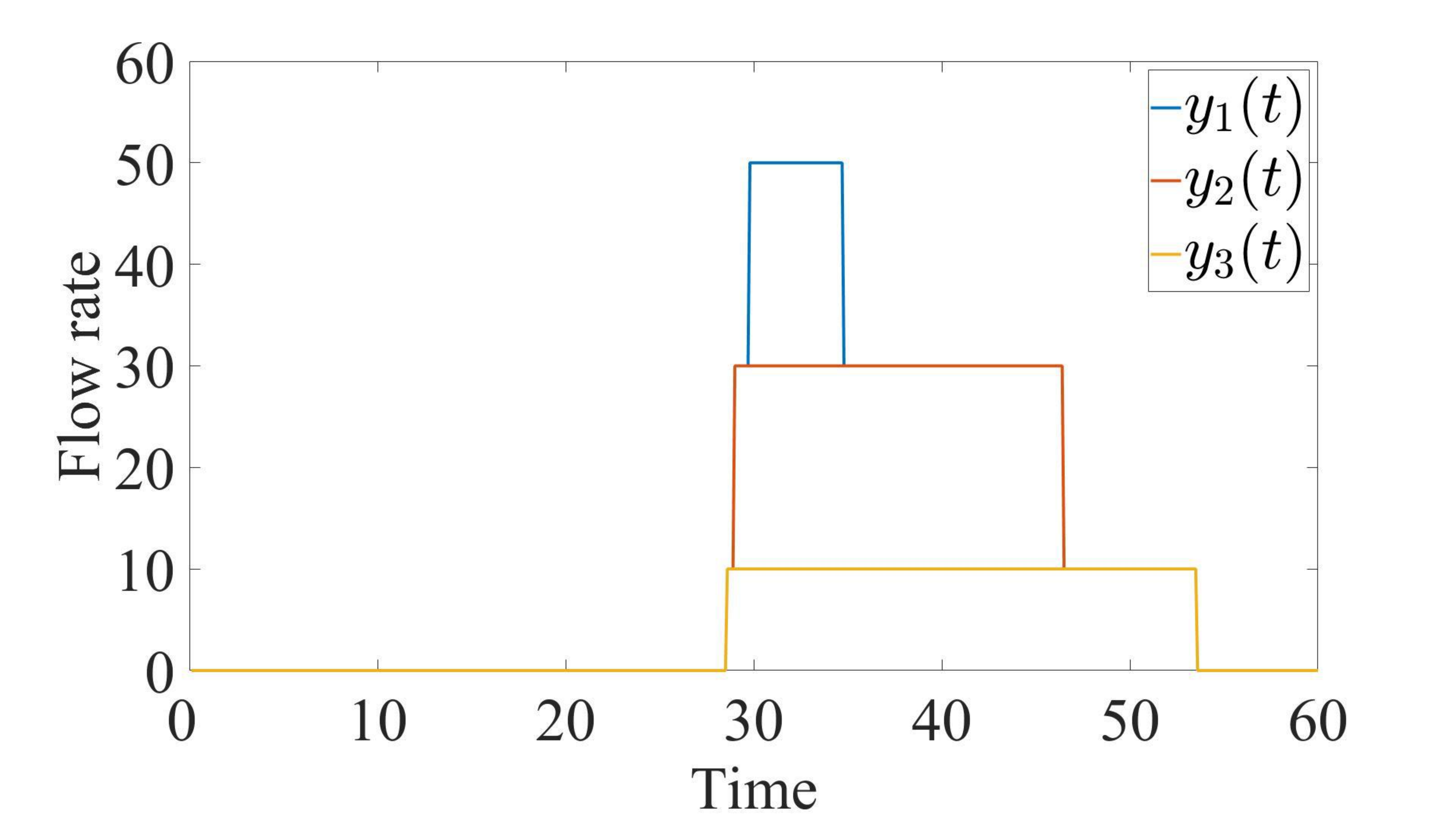}
		\subcaption{Arrival flow rate in the DSO state}
		\label{Fig:NE4_DSOFlows}
	\end{minipage}
	\begin{minipage}[t]{0.45\textwidth}
	\centering
		\includegraphics[clip, width=1.0\columnwidth]{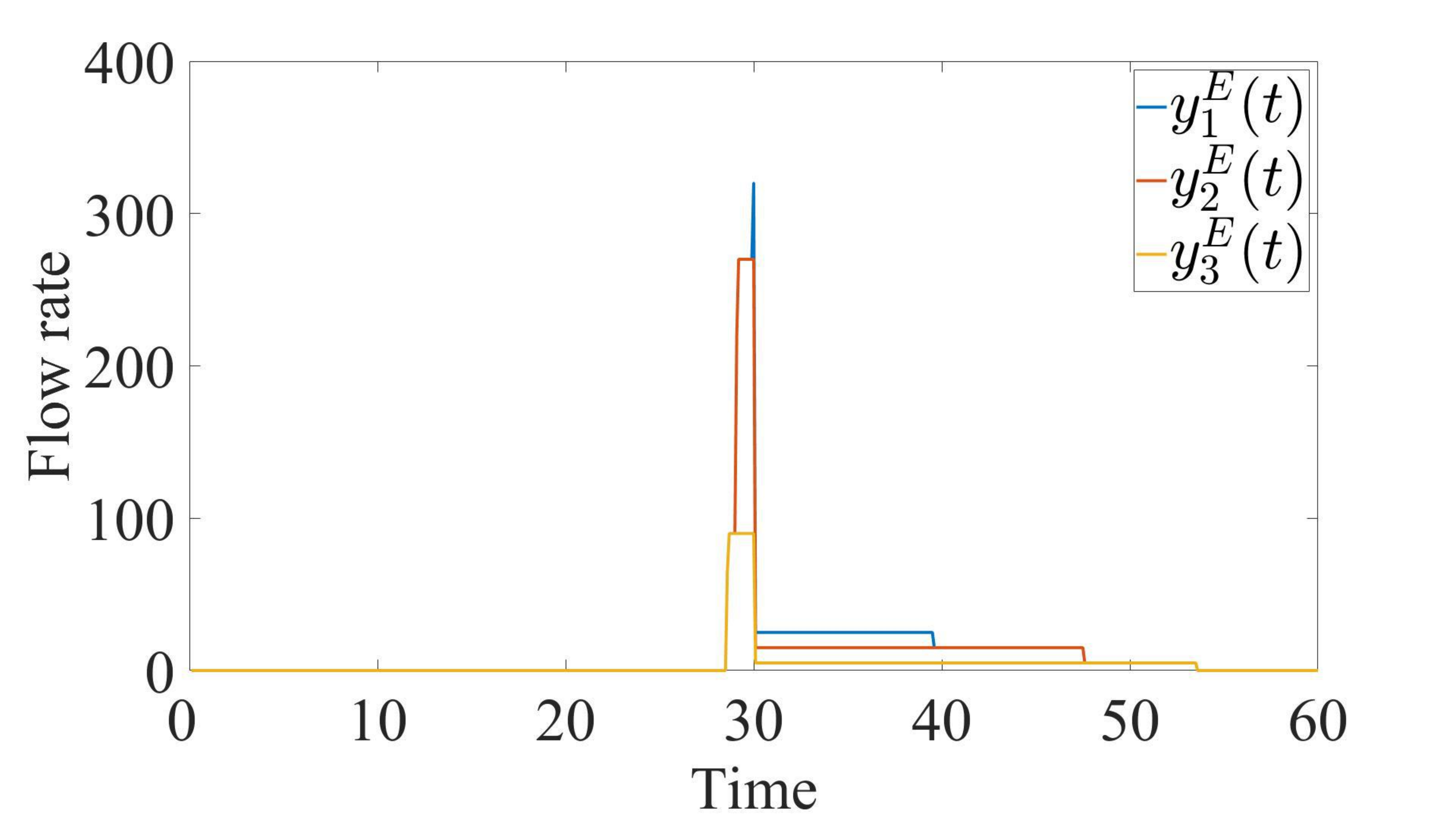}
		\subcaption{Arrival flow rate in the DUE state}
		\label{Fig:NE4_DUEFlows}
	\end{minipage}\vspace{5mm}\\
	\begin{minipage}[t]{0.45\textwidth}
	\centering
		\includegraphics[clip, width=1.0\columnwidth]{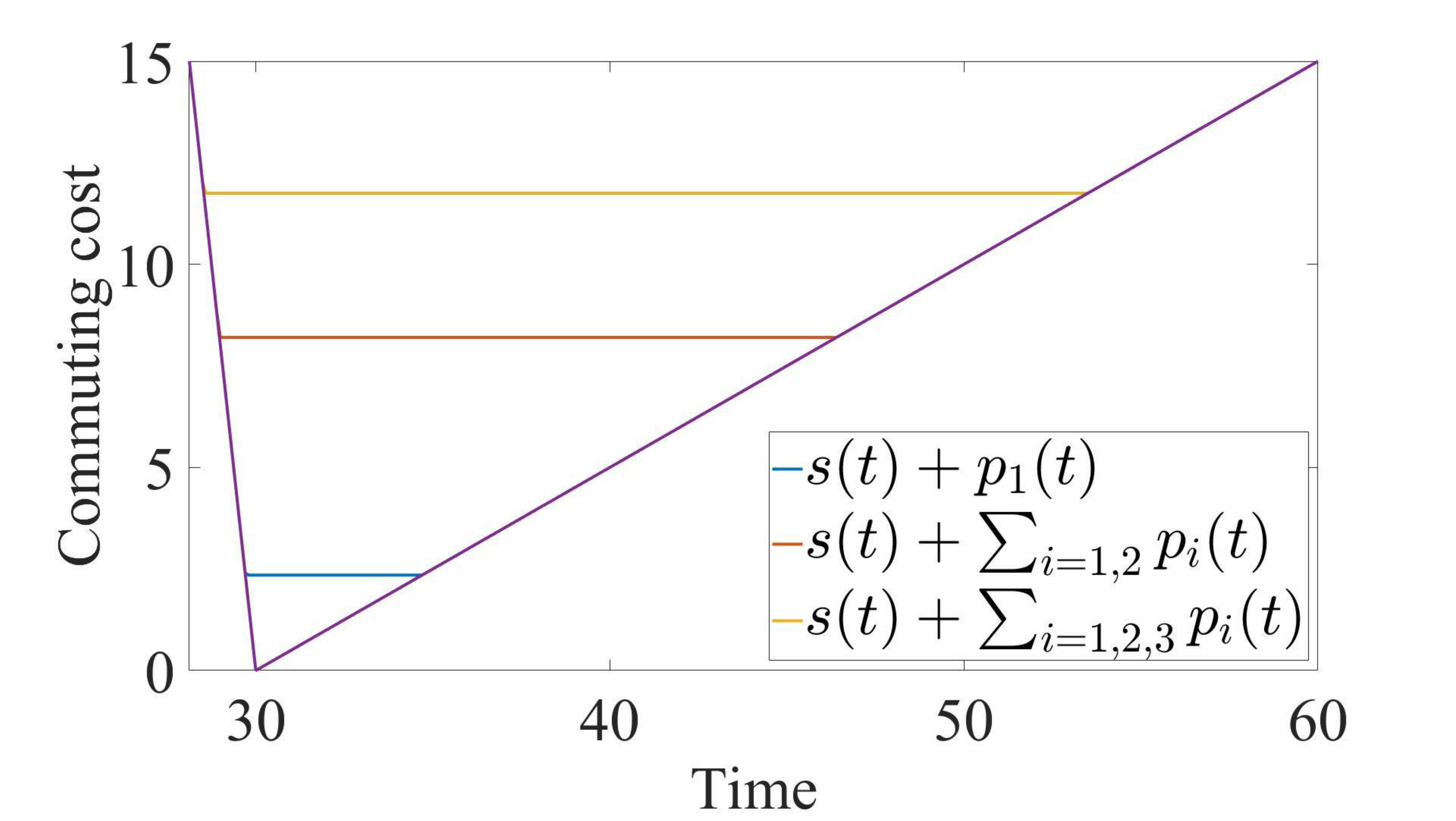}
		\subcaption{Equilibrium commuting costs in the DSO state}
		\label{Fig:NE4_DSOCosts}
	\end{minipage}
	\begin{minipage}[t]{0.45\textwidth}
	\centering
		\includegraphics[clip, width=1.0\columnwidth]{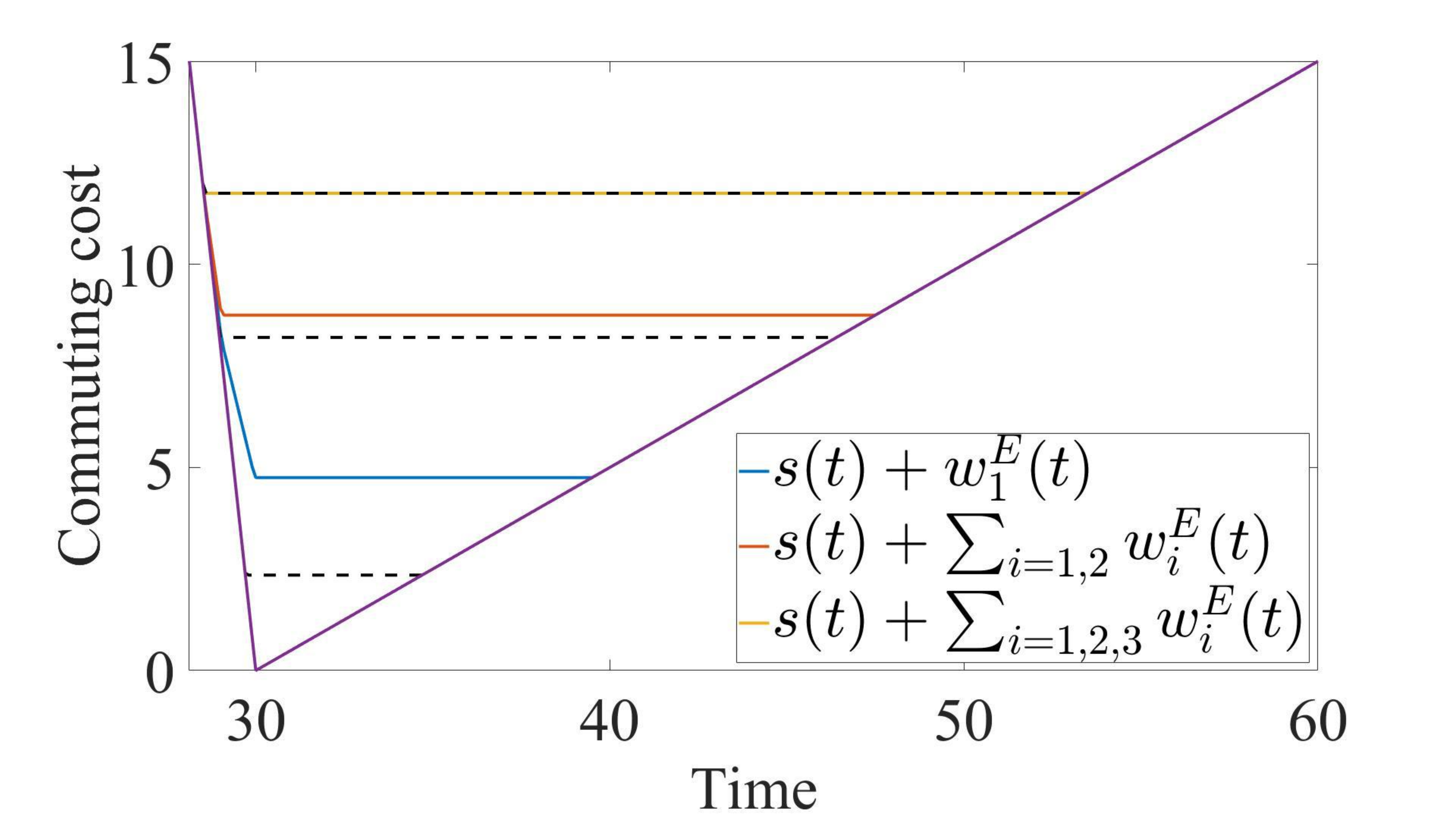}
		\subcaption{Equilibrium commuting costs in the DUE state}
		\label{Fig:NE4_DUECosts}
	\end{minipage}
	\vspace{0mm}
	\caption{Flow and cost patterns in DSO and DUE states in Example 4}
	\label{Fig:NE_Example4}
	\vspace{-2mm}
\end{figure*}

\subsection{Results: Morning commute}\label{sec:NumMor}
Figure~\ref{Fig:NE_Example1} shows the numerical results of Example 1.
Figure~\ref{Fig:NE1_CumFlows} shows the cumulative arrival and departure curves in the DSO (dashed lines) and DUE states (solid lines).
In this figure, the aggregate (first row) and disaggregate cumulative curves are illustrated.
Each disaggregate curve represents the origin-specific cumulative number of commuters at each bottleneck; 
on the other hand, the aggregate curve at each bottleneck is the sum of the disaggregate curves at the same bottleneck.
Note that in each illustration of cumulative curves, the earliest/latest destination arrival times of the corresponding aggregate/disaggregate flows in the DSO state are indicated by dotted vertical lines.
These times are the same as the times in the DUE state, as mentioned above.
Figures~\ref{Fig:NE1_DSOFlows}, \ref{Fig:NE1_DUEFlows}, \ref{Fig:NE1_DSOCosts}, and \ref{Fig:NE1_DUECosts} show the flow and cost patterns in the DSO and DUE states.

%, where the solid and dashed curves indicate the same meaning as those in Figures~\ref{fig:DUE_Example_CFC} and \ref{fig:DUEEve_Example_CFC}.
%The solid red and blue curves represent the DUE cumulative inflow and outflow curves, respectively.
%The dashed blue curves indicate the DSO cumulative inflow/outflow curves. 
%The time coordinates of the earliest/latest DSO outflow from each upstream location of bottlenecks (i.e., the time window of outflow) are marked by dotted vertical lines. 
%The first row of Figure~\ref{fig:CumFlowEx1} shows the cumulative curves of aggregate inflows/outflows for each bottleneck, and the second to fourth rows show disaggregated cumulative flows with respect to each location; top to bottom are locations 1 to 3, respectively. 

As can be seen from these figures, we can confirm that the aggregate bottleneck departure flow curves of the DSO and DUE states, namely, $D_i(\sigma_{i}(t))$ and $D^E_i(\sigma_{i}(t))$, are identical. 
However, the destination arrival flow curves with respect to each origin---represented by Eq.~\eqref{Eq:TESolution_Flow} for the DSO solutions and Eq.~\eqref{eq:DUESolution_Flow} for the DUE solution---are different.
The equilibrium commuting costs of the DSO and DUE solutions are the same, and the optimal prices $\mathbf{p}$ and queuing delay $\mathbf{w}^E$ are equal (i.e. \textbf{\prettyref{conj:EqualityQD_PP}}). 
These results verify the analysis of the DSO and DUE flow patterns, and are illustrated in Figures~\ref{Fig:NE1_DSOCosts} and \ref{Fig:NE1_DUECosts}.

The results of Example 2 are shown in Figure~\ref{Fig:NE_Example2}. 
The aggregate bottleneck departure flow curves of the DSO and DUE solutions are different, as shown in Figures~\ref{Fig:NE2_CumFlows}, \ref{Fig:NE2_DSOFlows}, and \ref{Fig:NE2_DUEFlows}. 
In addition, Figures~\ref{Fig:NE2_DSOCosts} and \ref{Fig:NE2_DUECosts} show that neither the commuting cost of the two solutions nor the optimal prices and queuing delay of each bottleneck are equal.

\subsection{Results: Evening commute}\label{sec:NumEve}

For the evening commute, the legends of the cumulative flow curves in Figures~\ref{Fig:NE_Example3} and \ref{Fig:NE_Example4} are the same as those in the morning commute. 
The time coordinates of the earliest/latest DSO departure flow heading to each downstream location of bottlenecks are indicated by dotted vertical lines.

Eqs.~\eqref{eq:DUEEve_Feasible_1} and \eqref{eq:DUEEve_Feasible_2} are satisfied in Example 3. 
As per Example 1 for the morning commute, the optimal prices and queuing delay of each bottleneck are equal (see Figures~\ref{Fig:NE3_DSOCosts} and \ref{Fig:NE3_DUECosts}). 
However, the aggregate bottleneck departure flows of the DSO and DUE solutions are different in Example 3, which is different from that observed in Example 1 (see Figure~\ref{Fig:NE3_CumFlows}). 
Instead, the two solutions exhibit the same destination arrival-flow curves, which are discerned from the subfigures of the most downstream bottleneck for each location in Figures~\ref{Fig:NE3_CumFlows} and are observed in Figure~\ref{Fig:NE3_DSOFlows} and \ref{Fig:NE3_DUEFlows}. 
These results verify the analysis of the flow patterns in the evening commute. 
As shown in Figures~\ref{Fig:NE3_DSOCosts} and \ref{Fig:NE3_DUECosts}, the equilibrium commuting costs of the DSO and DUE solutions are equal, and the optimal prices $\mathbf{p}$ and queuing delay $\mathbf{w}^E$ are equal as well (\textbf{\prettyref{conj:EveEqualityQD_PP}}).

In Example 4, the early arrival penalty is increased to violate Eq.~\eqref{eq:DUEEve_Feasible_2}. 
The results are shown in Figure~\ref{Fig:NE_Example4}. 
These figures show that both the aggregate flows and disaggregate flows of the DSO and DUE solutions are different. 
Thus, the equality between the commuting cost of the two solutions and that between dynamic toll and queuing delay is not observed under this setting. 

%These differences between the DSO and DUE solutions suggest that the analytical approach for the DUE problem is not appropriate.

%%%%%%%%%%%%%%%%%%%%%%%%%%%%%%%%%%%%%%%%%%%%%%%%%%%%%%%%%%%%%%%%%%%%
% Section 7

\section{Concluding remarks}\label{sec:Conclusion}

In this research, we studied the DSO and DUE traffic assignment problems in corridor networks, considering both morning and evening commutes. 
In the DSO problem without queues, we considered an equivalent equilibrium problem under a first-best TDM policy, and revealed an inclusion relationship of destination arrival-time (origin departure-time) windows and regularities of variables at the equilibrium in the morning (evening) commute. 
Based on these results, we derived the closed-form DSO solutions. 
In the DUE problem with queues, we first proved that, under certain conditions of a schedule delay function, the queuing delay at each bottleneck in the DUE assignment is exactly the same as an optimal toll that eliminates the queue at each bottleneck in the DSO assignment.
Based on this finding and the analytical DSO solutions, we successfully derived the closed-form DUE solutions. 
We also showed that the findings between DUE and DSO assignment provide a powerful basis for the welfare analysis of first-best/second-best TDM policies.
Furthermore, we identified the (a)symmetric properties of the assignment in the morning and evening commute problems.

The theory of the dynamic traffic assignment in the present study would be a building block for developing a theory for more general models. 
As a generalization of network topology,  dynamic assignment models in a tree network consisting of multiple corridors would be a suitable next step. 
It is also important to analyze the relationship between the DSO and DUE problems with heterogeneous commuters whose value of travel times is different.
Recently, \cite{akamatsu2021new} has revealed this relationship in a single bottleneck, and insights from this study would assist in the analyses.
Furthermore, incorporating the complex cost function for chains of activity (i.e., activity-based modeling) such as \cite{Li2014l} and/or more realistic traffic dynamics such as \cite{Lago2007} are also interesting directions.

% As a part of future works, it would be important to analyze the relationship between the DSO and DUE problems with heterogeneous commuters whose value of travel times are different.
% Recently, \cite{akamatsu2021new} has revealed this relationship in a single bottleneck, and insights from this study would assist in the analyses.
% we are interested in extending the theory to DUE assignment models in networks with a more general topology (e.g., a tree network). 
% Furthermore, incorporating the complex cost function for chains of activity (i.e., activity-based modeling) such as \cite{Li2014l} is interesting direction. 
% %analyzing the relationship between problems with commuters' route choices for further knowledge.


%%%%%%%%%%%%%%%%%%%%%%%%%%%%%%%%%%%%%%%%%%%%%%%%%%%%%%%%%%%%%%%%%%%%
% Acknowledgement

\section*{Acknowledgement}

The authors acknowledge Professor Shunsuke Hayashi and Takamasa Iryo for their comments on an early draft. 
The authors also express their gratitude to anonymous referees for their careful reading of the manuscript and useful suggestions. 
The present research was partially supported by funding from Japan Society for the Promotion of Science (KAKENHI 15H04053, 18K18916 and 20H02267).

%%%%%%%%%%%%%%%%%%%%%%%%%%%%%%%%%%%%%%%%%%%%%%%%%%%%%%%%%%%%%%%%%%%%
% Appendix
\appendix

%%%%%%%%%%%%%%%%%%%%%%%%%%%%%%%%%%%%%%%%%%%%%%%%%%%%%%%%%%%%%

\section{Algorithm}\label{app:algo}

\begin{algorithm}[h]
\caption{Constructing a reduced corridor network}
\begin{algorithmic}[1]
\REQUIRE{
Original corridor network, $\VtI\equiv[1,2,\cdots,I]$; Travel demand, $\VtQ$; Bottleneck capacities, $\Vt\mu$. 
}
\ENSURE{
Reduced corridor network, $\VtN$. %\equiv[1,2,\cdots,N]$.
}
\STATE{$\VtN\leftarrow\VtI$} (i.e. let the $i$th element of $\VtN$, $\VtN(i)$, equals to $i$);
\FOR{$i:=I-1$ to $1$} 
\WHILE{$i\neq~ length(\VtN)$ and $ Q_{\VtN(i)}/\hat\mu_{\VtN(i)} \geq Q_{\VtN(i+1)}/\hat\mu_{\VtN(i+1)}$}
\STATE{
Delete $\VtN(i+1)$ from $\VtN$;
}
\ENDWHILE
\ENDFOR
\end{algorithmic}
\end{algorithm}

The time complexity of this algorithm is $O(n)$, where $n$ is the number of locations in the original corridor network.

%%%%%%%%%%%%%%%%%%%%%%%%%%%%%%%%%%%%%%%%%%%%%%%%%%%%%%%%%%%%%
\renewcommand{\thelemma}{\Alph{section}.\arabic{lemma}}

\section{Proofs of lemmas, propositions and theorems}\label{app:proof}

%%%%%%%%%%%%%%%%%%%%%%%%%%%%%%%%%%%%%%%%%%%%%%%%%%%%%%
\subsection{Lemma~\ref{Lemma:DetectFBN}}\label{app:ProofBNP}

To prove \textbf{Lemma~\ref{Lemma:DetectFBN}}, we first show the following three lemmas as necessary components:

\begin{lemma}\label{lemma:CDBT}
If bottleneck $i\in\mathcal{N}\setminus\{1\}$ is \textit{non-false}, then for any downstream bottleneck $m~(m<i)$ of bottleneck $i$, $\ClT_m\subset\ClT_i$ is true.
\end{lemma}

\begin{lemma}\label{lemma:mui+1>=mui}
If bottleneck $i\in\mathcal{N}\setminus\{1\}$ is non-false, then $\mu_{i} < \mu_{i-1}$.
\end{lemma}

\begin{lemma}\label{lemma:PositivePrice}
If bottleneck $i\in\mathcal{N}$ is non-false, then  
\begin{align}
&p_i(t)>0, &&\forall t\in(t_i^-,t_i^+).
\end{align}
\end{lemma}

\noindent The proofs of these lemmas are provided in \ref{sec:ProofCDBT}, \ref{sec:Proofmui+1>=mui} and \ref{sec:ProofPositivePrice}, respectively. 
Given these lemmas, we prove \textbf{Lemma~\ref{Lemma:DetectFBN}}.

\subsubsection*{(1) Proof of necessity}

The assertion for necessity of ~\textbf{Lemma~\ref{Lemma:DetectFBN}} states:
\textit{If bottleneck $i\in\mathcal{N}$ is non-false, then all downstream bottlenecks of bottleneck $i$ have a smaller normalized demand than that of bottleneck $i$, i.e., $\forall m<i,m\inN$, subject to
\begin{align}
&\bar Q_m< \bar Q_i.\label{eq:ProofnonFBcriterion}
\end{align}}

 We prove the following two propositions (i) and (ii), because (ii) recursively proves this assertion for a given (i).
\begin{enumerate}[(i)]
\item Bottleneck 1 is \textit{non-false}.

\item Let the bottleneck $m~(m<i)$ be the closest downstream \textit{non-false} bottleneck of bottleneck $i$.
If bottleneck $i~ (i>1)$ is non-false, then $\forall l\geq m, l < i, l\inN$, subject to 
\begin{align}
    \bar Q_l<\bar Q_i.
\end{align}

\end{enumerate}

(i) If bottleneck 1 is false, i.e., $\rho_1=0$, all the commuters of location 1 should arrive at the destination at the desired arrival time and thus the capacity $\mu_1$ must be infinitely large. 
Therefore, bottleneck 1 is non-false for a finite bottleneck capacity $\mu_1$.

(ii) Because bottleneck $i$ is non-false, \textbf{Lemma~\ref{lemma:PositivePrice}} yields 
\begin{align}
&p_i(t)>0,&\forall t\in (t_i^-,t_i^+).
\end{align}
Substituting this into Eq.~\eqref{Eq:DSO_KKT2} yields
\begin{align}
&\sum_{j = i}^{N}q_j(t)=\mu_i,&\forall t\inT_i\label{eq:propBNPeq1}.
\end{align}
For the bottleneck $n(i)$ that is the closest upstream non-false bottleneck of bottleneck $i$, we also have
\begin{align}
&\sum_{j = n(i)}^{N}q_j(t)=\mu_{n(i)},&\forall t\inT_n\label{eq:propBNPeqn}
\end{align}
From \textbf{Lemma~\ref{lemma:CDBT}}, we know that $\ClT_j\subset\ClT_{n(i)},\forall j<n(i)$.
Thus, combining Eqs.~\eqref{eq:propBNPeq1} and \eqref{eq:propBNPeqn}, we obtain
\begin{align}
\sum_{j = i}^{n(i) - 1}Q_j&=\sum_{j = i}^{n(i) - 1} \int_{t\inT}q_j(t)\Rmd t &&(\because Eq.~\eqref{Eq:DSO_KKT3})\\
&\geq\sum_{j = i}^{n(i) - 1} \int_{t\inT_i}q_j(t)\Rmd t \\
&=\int_{t\inT_i}\left[\sum_{j = i}^{N} q_j(t)- \sum_{j = n(i)}^{N}q_j(t)\right]\Rmd t \\
&=\int_{t\inT_i}(\mu_i-\mu_{n(i)})\Rmd t &&(\because \prettyref{eq:propBNPeq1}, \prettyref{eq:propBNPeqn}, \ClT_i\subset\ClT_{n(i)})\\
&= (\mu_i-\mu_{n(i)}) \cdot T_i, \\
\Leftrightarrow \quad& \dfrac{\sum_{j = i}^{n(i) - 1}Q_j}{\mu_i-\mu_{n(i)}}\geq T_i.\label{eq:propBNPeq2}
\end{align}

For bottlenecks $l$, the joint arrival time window of location $\{l,l+1,\cdots,i-1\}$ is denoted by
\begin{align}
\breve\ClT_l\equiv\bigcup_{l\leq j <i}\ClT_j,
\end{align}
and $\breve T_l\equiv |\breve\ClT_l|$.
Then, Eq.~\eqref{Eq:DSO_KKT2} gives 
\begin{align}
&\sum_{j = l}^{N}q_j(t) \leq \mu_l,&\forall t\in\breve\ClT_l.\label{eq:propBNPeq3}
\end{align}
Combining \prettyref{eq:propBNPeq1} and \prettyref{eq:propBNPeq3} yields
\begin{align}
&\sum_{j = l}^{i-1}q_j(t)	\leq \mu_l-\mu_i,&\forall t\in\breve\ClT_l,\\
\Rightarrow \quad& \sum_{j = l}^{i-1} Q_j \leq (\mu_l-\mu_i)\cdot \breve T_l,\\
\Leftrightarrow \quad&  \dfrac{\sum_{j = l}^{i-1} Q_j}{\mu_l-\mu_i} \leq \breve T_l.\label{eq:propBNPeq4}
\end{align}

\noindent From \textbf{Lemma~\ref{lemma:CDBT}}, we know that $\breve\ClT_l\subset\ClT_i$. 
Thus, Eqs.~\eqref{eq:propBNPeq2} and \eqref{eq:propBNPeq4} can be combined to yield
\begin{align}
\dfrac{\sum_{j = i}^{n-1}Q_j}{\mu_i-\mu_n} \geq T_i > \breve T_l \geq \dfrac{\sum_{j = l}^{i-1}Q_j}{\mu_l-\mu_i}, \quad\Leftrightarrow\quad \bar Q_i>\bar Q_l.
\end{align}
This establishes the proof.

\subsubsection*{(2) Proof of sufficiency}

The sufficiency is proved by proving that the following contraposition of \textbf{Lemma~\ref{Lemma:DetectFBN}} is true: 
\textit{If bottleneck $i~(i>1)$ is false, then there exists a downstream bottleneck $m~(m<i)$ satisfying $\bar Q_m\geq \bar Q_i$.}

Let bottleneck $m(i)$ be the closest downstream non-false bottleneck of bottleneck $i$.
Therefore, $\bar Q_{m(i)} \geq \bar Q_i$.

The joint arrival-time window of locations $\{i,i+1,\cdots,n-1\}$ is denoted by
\begin{align}
\breve\ClT_i\equiv\bigcup_{i\leq j <n(i)}\ClT_j,
\end{align}
where bottleneck $n(i)$ is the closest upstream non-false bottleneck of bottleneck $i$. Also, $\breve T_i\equiv |\breve\ClT_i|$.

For bottleneck $i$, Eq.~\eqref{Eq:DSO_KKT2} states
\begin{align}
&\sum_{j = i}^{N}q_j(t)\leq \mu_i,&\forall t\in\breve\ClT_i.\label{eq:sumqi}
\end{align}
For bottleneck $n(i)$, combining Eq.~\eqref{Eq:DSO_KKT2} and \textbf{Lemma~\ref{lemma:PositivePrice}} yields
\begin{align}
&\sum_{j = n(i)}^{N}q_j(t)= \mu_{n(i)},&\forall t\inT_{n(i)}.\label{eq:sumqn}
\end{align}
\textbf{\prettyref{lemma:CDBT}} gives $\breve\ClT_i\subset\ClT_{n(i)}$. 
Therefore combining Eqs.~\eqref{eq:sumqi} and \eqref{eq:sumqn} yields
\begin{align}
&\sum_{j = i}^{n(i)-1} q_j(t)\leq \mu_i-\mu_{n(i)},&\forall t\in\breve\ClT_i,\\
\Rightarrow \quad& \sum_{j = i}^{n(i)-1} Q_j \leq (\mu_i-\mu_{n(i)})\cdot \breve T_i,\\
\Leftrightarrow \quad&  \dfrac{\sum_{j = i}^{n(i)-1} Q_j}{\mu_i-\mu_{n(i)}} \leq \breve T_i. \label{eq:sumQi}
\end{align}

For bottleneck $m(i)$, combining Eq.~\eqref{Eq:DSO_KKT2} and \textbf{\prettyref{lemma:PositivePrice}} provides
\begin{align}
&\sum_{j = m(i)}^{N}q_j(t) = \mu_{m(i)},&\forall t\inT_{m(i)}.\label{eq:sumqm}
\end{align}
Because \prettyref{lemma:CDBT} yields $\ClT_{m(i)}\subset\ClT_{n(i)}$, combining \prettyref{eq:sumqn} and \prettyref{eq:sumqm} provides
\begin{align}
\sum_{j = m(i)}^{n(i) - 1}Q_j&=\sum_{j = m(i)}^{n(i) - 1} \int_{t\inT}q_j(t)\Rmd t &&(\because Eq.~\eqref{Eq:DSO_KKT3})\\
&\geq\sum_{j = m(i)}^{n(i) - 1} \int_{t\inT_{m(i)}}q_j(t)\Rmd t \\
&=\int_{t\inT_{m(i)}}\left[\sum_{j = m(i)}^{N}q_j(t)- \sum_{j = n(i)}^{N}q_j(t)\right]\Rmd t \\
&=\int_{t\inT_{m(i)}}(\mu_i-\mu_{n(i)})\Rmd t &&(\because \prettyref{eq:sumqn}, \prettyref{eq:sumqm}, \ClT_{m(i)}\subset\ClT_{n(i)})\\
&= (\mu_{m(i)}-\mu_{n(i)}) \cdot T_i, \\
\Leftrightarrow \quad & \dfrac{\sum_{j = m(i)}^{n(i) - 1}Q_j}{\mu_{m(i)}-\mu_{n(i)}} \geq T_{m(i)}.\label{eq:sumQmn}
\end{align}

Combining Eqs.~\eqref{eq:sumQmn} and \eqref{eq:sumQi} yields that
\begin{align}
&\dfrac{\sum_{j = m(i)}^{n(i) - 1}Q_j}{\mu_{m(i)}-\mu_{n(i)}} 
\geq T_{m(i)} \geq \breve T_i \geq 
\dfrac{\sum_{j = m(i)}^{n(i) - 1}Q_j}{\mu_i-\mu_{n(i)}}, \quad \Leftrightarrow \quad \bar Q_{m(i)}\geq \bar Q_i.
\end{align}
This establishes the proof.

%%%%%%%%%%%%%%%%%%%%%%%%%%%%%%%%%%%%%%%%%%%%%%%%%%%%%%
\subsection{\prettyref{lemma:CDBT}}\label{sec:ProofCDBT}

Because bottleneck $i$ is non-false, we know $\rho_i-c_i > \rho_{i-1}-c_{i-1}$. 
In addition, $\rho_n - c_n \geq \rho_{n-1} - c_{n-1},\forall n\inN\setminus\{1\}$.
Thus, $\rho_i-c_i > \rho_m - c_m$, and the strictly quasi-convex schedule delay function yields that $\ClT_m\subset\ClT_i$.

%%%%%%%%%%%%%%%%%%%%%%%%%%%%%%%%%%%%%%%%%%%%%%%%%%%%%%
\subsection{\prettyref{lemma:mui+1>=mui}}\label{sec:Proofmui+1>=mui}

We prove the following contraposition: \textit{If $\mu_{i} \geq \mu_{i-1}$, then bottleneck $i~(i>1)$ is false.}

Because $Q_{i-1}>0$, there must exist some $t_0\inT_{i-1}\subset \ClT_{i}$ satisfying $q_{i-1}(t_0)>0$. 
Therefore, 
\begin{align}
\sum_{j =  i}^{N}q_j(t_0) & < \sum_{j = i-1}^{N}q_j(t_0) &&(\because  q_{i-1}(t_0)>0)\\
&\leq \mu_{i-1} &&(\because  Eq.~\eqref{Eq:DSO_KKT2})\\
&\leq \mu_{i}. 
\end{align}
Condition \eqref{Eq:DSO_KKT2} yields that $p_{i}(t_0)=0$. 
This further yields that
\begin{align}
\rho_{i-1}-c_{i-1} &= v_{i-1}(t_0) - c_{i-1}&&(\because  q_{i-1}(t_0)>0)\\
&= v_{i}(t_0) - c_{i} &&(\because  p_{i}(t_0)=0)\\
&\geq \rho_{i} - c_{i}   &&(\because  Eq.~\eqref{Eq:DSO_KKT1}).
\end{align}
Thus, we know that bottleneck $i+1$ is false, which establishes the proof.

%%%%%%%%%%%%%%%%%%%%%%%%%%%%%%%%%%%%%%%%%%%%%%%%%%%%%%
\subsection{\prettyref{lemma:PositivePrice}}\label{sec:ProofPositivePrice}

Proof by contradiction.
Suppose there exists $t_0\in (t_i^-,t_i^+)$ satisfying $p_i(t_0)=0$.
Because bottleneck $i$ is non-false, we have,
\begin{align}
    v_{i-1}(t_0) - c_{i-1} &= v_i(t_0) - c_i && (\because p_i(t_0)=0)\label{eq:B4Eq1}\\
    &= \rho_i - c_i &&(\because Eq.~\eqref{Eq:DSO_KKT1})\label{eq:B4Eq2}\\
    &> \rho_{i-1}-c_{i-1} &&(\because \textrm{bottleneck $i$ is non-false}).\label{eq:B4Eq3}
\end{align}
This yields $q_{i-1}(t_0)=0$, and
\begin{align}
\sum_{j = i-1}^{N} q_j(t_0) &= \sum_{j = i}^{N} q_j(t_0) &&(\because q_{i-1}(t_0)=0)\label{eq:B4Eq4}\\
&\leq \mu_{i} &&(\because Eq.~\eqref{Eq:DSO_KKT2})\label{eq:B4Eq5}\\
&< \mu_{i-1} &&(\because \prettyref{lemma:mui+1>=mui}).\label{eq:B4Eq6}
\end{align}
This means that $p_{i-1}(t_0)=0$. 

By recursively applying Eqs.~\eqref{eq:B4Eq1}-\eqref{eq:B4Eq3} and Eqs.~\eqref{eq:B4Eq4}-\eqref{eq:B4Eq6} for all downstream bottlenecks, we obtain
\begin{align}
&q_n(t_0)=0 \quad \textrm{and} \quad p_n(t_0)=0,&&\forall n< i.
\end{align}
Thus, commuters of location $i$ with arrival times $t\neq t_0$ would have an incentive to choose $t_0$ because $\rho_i -c_i > s(t_0)$. 
This means that it is not an equilibrium state and leads to contradiction. 
Therefore, such $t_0$ does not exist and $p_i(t)>0,\forall t\in(t_i^-,t_i^+)$.

%%%%%%%%%%%%%%%%%%%%%%%%%%%%%%%%%%%%%%%%%%%%%%%%%%%%%%
\subsection{\prettyref{lemma:CoT}}

Refer to the proof of \prettyref{lemma:CDBT} in \ref{sec:ProofCDBT} because \prettyref{lemma:CDBT} immediately states this proposition.

%%%%%%%%%%%%%%%%%%%%%%%%%%%%%%%%%%%%%%%%%%%%%%%%%%%%%%
\subsection{\prettyref{lemma:NonnegPrice}}
First, for $t\in(t_i^-,t_i^+)$, we can see that the lemma is true by referring to the proof of \textbf{\prettyref{lemma:PositivePrice}} in \ref{sec:ProofPositivePrice}.
In addition, for $t\notin(t_i^-,t_i^+)$, if $i=N$, the condition \eqref{Eq:DSO_KKT2} immediately yields $p_N(t)=0,\forall t\notin(t_N^-,t_N^+)$.
Thus, we only have to prove the lemma for arbitrary $t\notin(t_i^-,t_i^+)$ when $i<N$.

We prove this by contradiction.
Suppose that $t_0\notin(t_i^-,t_i^+)$ satisfies $p_i(t_0)>0$. 
For bottleneck $i$, Eq.~\eqref{Eq:DSO_KKT2} states that
\begin{align}
\sum_{j = i+1}^{N}q_j(t_0) &= \sum_{j = i}^{N}q_j(t_0) &&(\because q_i(t_0)=0)\\
&=\mu_i &&(\because p_i(t_0)>0).\label{eq:sumqi+1mui}
\end{align}
For bottleneck $i+1$, Eq.~\eqref{Eq:DSO_KKT2} gives that
\begin{align}
\sum_{j = i+1}^{N}q_j(t_0)\leq \mu_{i+1}.\label{eq:sumqimui+1}
\end{align}
Furthermore, combining Eqs.~\eqref{eq:sumqi+1mui} and \eqref{eq:sumqimui+1} yields that
\begin{align}    
\mu_{i+1} \geq \mu_i.\label{eq:mui+1mui}
\end{align}
Thus, \prettyref{lemma:mui+1>=mui} concludes that bottleneck $i+1$ is false, which results in a contradiction.
Therefore, such $t_0$ does not exist and $p_i(t)=0,\forall t\notin (t_i^-,t_i^+)$.

%%%%%%%%%%%%%%%%%%%%%%%%%%%%%%%%%%%%%%%%%%%%%%%%%%%%%%
\subsection{\prettyref{lemma:Equality_OutFlow}}\label{app:proofx=x}

Consider the two of cases $t\inT_i$ and $t\notin\ClT_i$ for this lemma:

(1) If $t\inT_i$, the demand-supply equilibrium condition \eqref{eq:TNP_PPE_2} and the queuing delay condition \eqref{eq:DUECond_QueueDelay_2} for each bottleneck $i$ yields
\begin{align}
&x^E_i(\sigma^E_i(t))=\mu_i=x_i(\sigma_i(t)),& \forall i\inN,t\inT_i\label{eq:x=xTi}.
\end{align}

(2) If $t\notin\ClT_i$, for each upstream bottleneck $j~(j>i)$ of bottleneck $i$ in the DSO problem, the DSO solution yields
\begin{align}
&x_i(\sigma_i(t)) = \mu_j, &\forall t\inT_j\setminus\ClT_{j-1}.\label{eq:Equality_OutFlow_eq1}
\end{align} 
For the DUE problem, we have
\begin{align}
    \int_{t\inT_i}q_i^E(t)\Rmd t
    &=\int_{t\inT_i}[\dot\sigma^E_i(t)\cdot x_i^E(\sigma^E_i(t))-\dot\sigma^E_{i+1}(t)\cdot x_{i+1}^E(\sigma^E_{i+1}(t))]\Rmd t \quad (\because Eqs.~\eqref{eq:Definition_DUE_Y},\eqref{eq:DUECond_FC_1})\\
    &=\int_{t\inT_i}[(1-\sum_{j<i}\dot p_j(t))\cdot \mu_i-(1-\sum_{j\leq i}\dot p_j(t))\cdot \mu_{i+1}]\Rmd t \quad (\because Eqs.~\eqref{eq:Equality_QD_PP}, \eqref{eq:x=xTi})\\
    &=\int_{t\inT_i}[(1-\sum_{j<i}\dot p_j(t))\cdot \hat\mu_i + \dot p_i(t)\cdot \mu_{i+1}]\Rmd t \\
    &=\int_{t\inT_i}[(1-\sum_{j<i}\dot p_i(t))\cdot \hat\mu_i]\Rmd t\\ &=\hat\mu_i \cdot \int_{t\inT_{i}}(1+\dot s(t))\Rmd t\\
    &=\hat\mu_i\cdot T_i \\
    &=Q_i\qquad \qquad (\because \prettyref{eq:ATWLengthinCorridor}).
\end{align}
Then, from the flow conservation condition \eqref{eq:DUECond_FCMass}, we know that 
\begin{align}
&q_i^E(t)=0,&t\notin\ClT_i.
\end{align} 
Thus, for each upstream bottleneck $j~(j>i)$ of bottleneck $i$ in the DUE problem, 
\begin{align}
x^E_i(\sigma^E_i(t))\cdot\dot\sigma_i^E(t)&= x^E_j(\sigma^E_j(t))\cdot\dot \sigma_j^E(t), &\forall t\inT_j\setminus\ClT_{j-1}.
\end{align}
Because $\dot\sigma_i^E(t)=\dot \sigma_j^E(t)=1,\forall t\inT_j\setminus\ClT_{j-1}$, 
\begin{align}
&x^E_i(\sigma^E_i(t)) =x^E_j(\sigma^E_j(t)) =\mu_j, &\forall t\inT_j\setminus\ClT_{j-1}.\label{eq:Equality_OutFlow_eq2}
\end{align} 
Thus, combining Eqs.~\eqref{eq:Equality_OutFlow_eq1} and \eqref{eq:Equality_OutFlow_eq2} yields $x^E_i(\sigma^E_i(t))=x_i(\sigma_i(t)),\forall t\notin\ClT_i$.

%\end{thebibliography}
\addcontentsline{toc}{section}{References}
\bibliographystyle{elsarticle-harv} 
\bibliography{DSODUErefs}

\end{document}